\documentclass[11pt,reqno]{article}

\usepackage[utf8]{inputenc}
\usepackage[T1]{fontenc}

\usepackage{savesym}
\usepackage{hyperref}
\usepackage{amsfonts}
\usepackage{amsmath}
\usepackage{amsthm}
\usepackage{amssymb}
\usepackage{mathrsfs}
\usepackage{amstext}
\usepackage{color}
\font\msbm=msbm10

\savesymbol{urcorner}\savesymbol{not}
\usepackage{mathabx} 
\restoresymbol{abx}{urcorner}
\restoresymbol{abx}{not}

\usepackage{stmaryrd}
\usepackage{euscript}
\usepackage{dsfont}

\numberwithin{equation}{section}

\theoremstyle{plain}
\newtheorem{Theorem}{Theorem}[section]
\newtheorem{lemma}[Theorem]{Lemma}
\newtheorem{corollary}[Theorem]{Corollary}
\newtheorem{example}[Theorem]{Example}
\newtheorem{proposition}[Theorem]{Proposition}
\newtheorem{definition}[Theorem]{Definition}

\newtheorem{remark}[Theorem]{Remark}

\def\mathbb#1{\hbox{\msbm{#1}}}

\newcommand{\re}{\mathbb{R}}\newcommand{\N}{\mathbb{N}}
\newcommand{\zz}{\mathbb{Z}}\newcommand{\C}{\mathbb{C}}

\newcommand{\com}{\mathbb{C}}
\newcommand{\Z}{{\zz}^d}
\newcommand{\n}{{\N}_0}
\newcommand{\R}{{{\re}^d}}
\newcommand{\phase}{\mathbb{S}^1}

\newcommand{\bproof}{\noindent {\bf Proof.}\, }
\newcommand{\eproof}{\hspace*{\fill} \rule{3mm}{3mm}}

\newcommand{\x}{\mathbf{x}}

\newcommand{\supp}{{\rm supp \, }}

\newcommand{\esssup}[1]{\underset{#1}{\operatorname{ess\,sup}}}
\newcommand{\essinf}[1]{\underset{#1}{\operatorname{ess\,inf}}}

\DeclareMathOperator*{\spn}{span}

\newcommand{\cp}{{\mathcal P}}

\newcommand{\ch}{{\mathcal H}}

\newcommand{\cf}{{\mathcal F}}

\newcommand{\cg}{{\mathcal G}}

\newcommand{\osc}{{\operatorname{osc}}}

\newcommand{\norm}[2]{\left\|\left.{#1}\right|{#2}\right\|}

\newcommand{\punkt}{\text{ .}}
\newcommand{\komma}{\text{ , }}
\newcommand{\e}{\textrm{e}}
\newcommand{\HLM}{\mathcal{M}}

\newcommand{\p}{{p(\cdot)}}
\newcommand{\q}{{q(\cdot)}}
\newcommand{{\s}}{{s(\cdot)}}

\newcommand{\qconst}{{\tilde{q}}}		

\renewcommand{\P}{\mathcal{P}(\R)}
\newcommand{\Plog}{\mathcal{P}^{\log}(\R)}
\newcommand{\vB}{\nu_{w,p(\cdot),q(\cdot)}}

\newcommand{\Co}{\mathsf{Co}}

\newcommand{\Lpp}{L_{\p}(\R)}

\newcommand{\Pe}{P^w_{p(\cdot),q(\cdot),a}}
\newcommand{\Le}{L^w_{p(\cdot),\tilde{q},a}}

\newcommand{\LizTri}{F^w_{p(\cdot),q(\cdot)}(\R)}
\newcommand{\Besov}{B^w_{p(\cdot),\qconst}(\R)}

\begin{document}
\title{General coorbit space theory for quasi-Banach spaces\\
and inhomogeneous function spaces with variable smoothness and
integrability}

\author{Henning Kempka$^{1,}$\footnote{\tt Corresponding author,
Email: henning.kempka@mathematik.tu-chemnitz.de}, Martin Sch\"afer$^{2}$, Tino
Ullrich$^3$\\\\
$^{1}$Faculty of Mathematics, Technical University Chemnitz\\
Reichenhainer Stra\ss e 39, 09126 Chemnitz, Germany
\\\\
$^{2}$Mathematical Institute, Technical University of Berlin\\
Stra{\ss}e des 17. Juni 136, 10623 Berlin, Germany
\\\\
$^{3}$Hausdorff Center for Mathematics \& Institute for Numerical Simulation\\
Endenicher Allee 60, 53115 Bonn, Germany
}

\maketitle

\begin{abstract}
In this paper we propose a general coorbit space theory suitable to define
coorbits of quasi-Banach spaces using an abstract continuous frame, indexed by a
locally compact Hausdorff space, and an
associated generalized voice transform. The proposed theory realizes
a further step in
the development of a universal abstract theory towards various function spaces
and their atomic decompositions which has been initiated by Feichtinger and
Gr{\"o}chenig in the late 1980ies. We combine the recent approaches in Rauhut,
Ullrich \cite{RaUl10} and Rauhut \cite{ra05-3} to identify, in particular,
various inhomogeneous (quasi-Banach) spaces of Besov-Lizorkin-Triebel
type. To prove the potential of our new theory we apply it to spaces
with variable smoothness and integrability which have attracted significant
interest in the last 10 years. From the abstract discretization machinery we
obtain atomic decompositions as well as wavelet frame isomorphisms for these
spaces.
\end{abstract}

\noindent
{\bf Key Words:} Coorbit space theory, continuous wavelet transform,
Besov-Lizorkin-Triebel type spaces,
variable smoothness, variable integrability,
2-microlocal spaces, Peetre maximal function, atomic decomposition, wavelet
bases

\medskip

\noindent
{\bf AMS Subject classification:} 42B25, 42B35, 46E35, 46F05.

\medskip


\section{Introduction}


The birth of coorbit theory dates back to the 1980ies,
starting with a series of papers by Feichtinger and Gröchenig~\cite{FeGr86, Gr88, Gr91}.
The main intention was to characterize function spaces via an abstract
transform, the so-called voice transform.
In the original setup, this transform is determined by an integrable irreducible representation
of a locally compact group on a Hilbert space $\mathcal{H}$ unifying e.g.\ the continuous wavelet transform, the
short-time Fourier transform, and the recent shearlet transform, to mention just
a few. More recently, representations which are not necessarily irreducible nor integrable have been considered~\cite{dadelasttevi14}. They allow
to treat, for instance, Paley-Wiener spaces and spaces related to Shannon wavelets and Schr\"odingerlets.

Classical examples of coorbit spaces associated to the continuous
wavelet transform on the $ax+b$-group are the homogeneous Besov-Lizorkin-Triebel
spaces \cite{Tr83,Tr88,Tr92}, identified rigorously as coorbits in Ullrich
\cite{T10}. What concerns further extensions of these spaces and interpretations
as coorbits we refer to Liang et al.\ \cite{LiSaUlYaYu11, LiSaUlYaYu12}.
More general wavelet coorbit spaces associated to a semidirect product $G=\R\rtimes H$, with a suitable subgroup $H$ of $GL(\R)$ as dilation group,
have been studied in \cite{fu13a,fu13b,furaitou15} and could recently be identified with certain decomposition spaces on the Fourier domain \cite{fuvoigt14}.
A specific example of this general setup is the shearlet transform, where $G$ is the shearlet group. The
associated shearlet spaces have first been studied in \cite{dakustte09}.
Other coorbit spaces, based
on a voice transform different from the wavelet transform, are e.g.\ modulation spaces \cite{gr01,fe83-4} and Bergman spaces \cite{FeGr86}.

Coorbit theory thus covers a great variety of different function spaces.
The underlying group structure however turns out to be a severe restriction for the theory since the identification of, e.g., inhomogeneous spaces of the above
type was long time not possible, however desirable. For that reason the theory
has evolved and several subsequent contributions have weakened among others the
assumption that the voice transform is supported on a locally compact
group. For instance, Dahlke, Steidl, and Teschke replaced it by a homogeneous
space, i.e.,\ a quotient of a group with a subgroup, with the aim to treat
functions on manifolds \cite{dastte04,dastte04-1,daforastte08}.

The starting point for the general coorbit space theory presented in this
paper is the approach used by Fornasier and Rauhut~\cite{fora05}, which was later revised and extended in~\cite{balhol10} and further expanded in~\cite{RaUl10}. There, the
group structure is abandoned completely and the voice transform is determined
solely by an abstract continuous frame $\mathcal{F} = \{\varphi_x\}_{x\in X}$ in $\mathcal{H}$
indexed by a locally compact Hausdorff space $X$ (not necessarily a group), i.e., $X$ is equipped with a Radon measure $\mu$
such that the map $x\mapsto\varphi_x$ is weakly measurable and that with constants $0<C_1,C_2<\infty$
\begin{equation}\label{eq:stab}
    C_1\|f|\mathcal{H}\|^2 \leq \int_{X} |\langle f,\varphi_x \rangle|^2 d\mu(x) \leq C_2\|f|\mathcal{H}\|^2\quad \mbox{for all }
    f\in \mathcal{H}\,.
\end{equation}
(Note that weak measurability of $x\mapsto\varphi_x$ in $\mathcal{H}$ implies that the integral in \eqref{eq:stab} is well-defined.)

We combine the approach in \cite{RaUl10} with ideas from \cite{ra05-3} to
define even coorbits
$$
  \Co(\cf,Y) := \{f~:~\langle f, \varphi_x \rangle \in Y\}
$$
of quasi-Banach spaces $Y$ using the general voice transform associated to
$\cf$. We thereby also recall the
relevant details of the existing theory, especially from \cite{fora05,RaUl10}
and fix some earlier inaccuracies. The developed theory yields noteworthy
generalizations even for the Banach case,
e.g.\
some assumptions made in \cite{fora05,RaUl10} can be weakened, such as the
uniform boundedness of the analyzing frame $\cf$ or some technical restrictions
on the weights and the admissible coverings. Most notably however, we can
generalize the main results of the discretization theory, which is possible
since we take a different -- more direct -- route to establish them.
It turns out that the three essential Lemmas \ref{auxlem:mainanalysis},
\ref{auxlem:mainsynthesis2}, and \ref{auxlem:Uinvert} below constitute the
technical foundation for the proof of the general abstract discretization
results in Theorems \ref{thm:atomicdec} and \ref{thm:frameexp}. Putting these
lemmas at the center of the exposition simplifies many arguments and
allows for a systematic approach towards new abstract discretization results. In
fact, we obtain
discrete characterizations of coorbit spaces by ``sampling'' the function using
a sub-sampled discrete frame $\cf_d = \{\varphi_{x_i}\}_{i\in I}$ on a
suitable index set $I$. Of course, as usual in coorbit space theory, there are
several technical assumptions to check. However, a great advantage of the
presented discretization machinery is
the fact that it provides a straight path towards discretization, where matters essentially reduce to checking properties (associated to $Y$) of
the analyzing frame $\cf$. This is in contrast to the usual approach where atomic decompositions and wavelet
characterizations, useful to study embeddings, $s-$numbers, interpolation
properties etc., are often developed from scratch for different related function
spaces.

To prove the potential of the theory presented here, we apply it to
identify spaces with variable smoothness and integrability, so-called variable
exponent spaces, as coorbits. Triebel-Lizorkin spaces of this kind are defined via the
quasi-norm
\begin{equation}\label{f000}
   \|f|F^w_{\p,\q}(\R)\| = \Big\|\Big(\sum\limits_{j=0}^{\infty}
    |w_j(\cdot)(\Phi_j \ast f)(\cdot)|^{q(\cdot)}\Big)^{1/q(\cdot)}|L_\p(\R)\Big\|\,,
\end{equation}
where the functions $w_j$ are weights and $\Phi_j$ are frequency filters corresponding to a dyadic decomposition of the frequency plane.
For the precise formulation see Definition \ref{inhom} below.
The functions $p(\cdot), q(\cdot)$ represent certain integrability
parameters, which may vary in the spatial variable $x$ of the space.
The
$2$-microlocal weight sequence $w_j(\cdot)$ determines the variable smoothness,
see \cite{KeVybDiff} for details. Function spaces with variable exponents are a
fast developing field thanks
to its many applications in stochastics, fluid dynamics and image processing,
see \cite{DieningHastoRoudenko2009} and \cite{DieningHastoBuch2011} and
references therein. The Lebesgue spaces $L_{p(\cdot)}(\R)$ with variable
integrability, see Definition \ref{Lppunkt} below, were already used by Orlicz
\cite{Orlicz31}. Recent contributions by Diening \cite{Diening2004} on the
boundedness of the Hardy-Littlewood maximal operator on $L_\p(\R)$ make them
accessible for harmonic analysis issues.

Surprisingly, the spaces \eqref{f000} can be handled within the generalized
coorbit space theory presented in this paper. In fact, due to unbounded left
and right translation operators (within the $ax+b$-group) a coorbit
characterization of homogeneous spaces of the above type already seems to be
rather impossible
at first glance. However, we are able to identify them as coorbits
$\Co(\cf,Y)$ of, what we call, Peetre-Wiener type spaces $Y$ by using a
suitable continuous frame $\cf = \{\varphi_x\}_{x\in X}$ with the index set $X =
\R \times [(0,1) \cup \{\infty\}]$. These spaces $Y$ are solid quasi-Banach
function spaces (QBF) defined on $X$, see Section \ref{ssec:PeetSp} below.
Peetre-Wiener type spaces can be
seen as a mixture of the Peetre type spaces introduced in \cite{T10}, and
certain Wiener amalgam spaces, see \cite{feGr89a}, \cite{ra05-4}. They appear
naturally when dealing with continuous local mean characterizations, a strategy
developed in \cite{T10} and \cite{LiSaUlYaYu11}. In fact, we show in Subsection
\ref{clm} below  that with large enough $a>0$ the quantity
\begin{equation}\label{coorb001}
\begin{split}
\|f|F^w_{\p,\q}(\R)\|_3 &= \|w(\cdot,\infty)\langle \Phi^{\ast}_0
f\rangle_{a}(\cdot)|L_\p(\R)\|\\
&+ \Big\|\Big(\int_{0}^1
|w(\cdot,t)\langle\Phi^{\ast}_t f\rangle_{a}(\cdot)|^{q(\cdot)}\frac{dt}{t}\Big)^{1/q(\cdot)}
|L_\p(\R)\Big\|\,
\end{split}
\end{equation}
represents an equivalent characterization for $F^w_{\p,\q}(\R)$. Here
$$
    \langle\Phi^{\ast}_t f\rangle_{a}(x):= \sup_{\substack{z\in \R
\\t/2\le\tau\le 2t, \tau<1}} \frac{|(\Phi_{\tau} \ast f)(x+z)|}
    {(1+|z|/\tau)^{a}}
$$
denotes the corresponding maximal function, which is essentially a
modification of the widely used Peetre maximal function, see
\eqref{Peemax} below, and is used in the definition of the Peetre-Wiener type
spaces, see Definition \ref{PeetreWiener}. Now the representation
\eqref{coorb001} is actually the identification of $F^w_{\p,\q}(\R)$
as a coorbit space of a Peetre-Wiener type space. Applying the abstract theory, in
particular Theorem \ref{thm:frameexp}, we obtain biorthogonal wavelet
expansions \cite{CoDaFe92} of the respective coorbit spaces. We describe the
application of the machinery for the rather simple (orthogonal) Meyer
wavelets, see Appendix \ref{sect:OWT}. Due to its generality, a straightforward
modification of Theorem \ref{thm:frameexp} leads to general (biorthogonal)
wavelet expansions and other tight discrete wavelet frames.

Let us mention that the continuous local mean characterizations
\eqref{coorb001} of spaces with variable exponents, see also Theorem
\ref{thm:contchar}, are new and interesting for their own
sake. In fact, one has to deal with additional difficulties
since a version of the classical Fefferman-Stein maximal
inequality, a crucial tool in this respect, is in general not true in
$L_{\p}(\ell_{\q})$ if $q(\cdot)$ is non-constant.

Finally, the provided discretizations of such spaces are not
entirely new. In
\cite{Ke11} the author used a different technique in order to obtain
discretizations with Meyer and Daubechies wavelets. However, let us mention that
the abstract Theorems \ref{thm:atomicdec}, \ref{thm:frameexp} below neither
restrict to orthonormal wavelets nor compactly supported atoms.

\subsection{Outline}

The paper is structured as follows. The abstract theory is established in
Section~\ref{sec:abstrth}. It generalizes earlier contributions, especially
\cite{fora05,RaUl10}, and in particular now includes the quasi-Banach case. In
Section~\ref{sec:varint} we give a short introduction to variable exponent
spaces, which will serve as our demonstration object for a concrete application
of the theory. We will utilize a new continuous local means characterization
in Section~\ref{sec:appcoorbit} to identify them as coorbits of a new
scale of Peetre-Wiener type spaces.
The abstract theory then yields atomic decompositions as well as discrete characterizations via wavelet frames.
Some useful facts concerning the continuous and discrete (orthogonal) wavelet
transform are collected in the Appendix.  

\subsection{Notation}

The symbols $\N, \n, \zz, \re,\re_+$, and $\com$ denote
the natural numbers, the natural numbers including $0$, the integers, the real
numbers, the non-negative real numbers, and the complex numbers. For a real number $t\in\re$ we put
$(t)_+=\max\{t,0\}$ and $(t)_-=\min\{t,0\}$. The conjugation of $z\in\C$ is
denoted by $\overline{z}$.
Let us emphasize that $\R$ has the usual meaning and $d\in\N_0$ is reserved
for its dimension. The symbol $|\cdot|$ denotes the Euclidean norm on $\R$
and $|\cdot|_1$ the $\ell_1$-norm.

The space of all sequences with entries in some set $M$ over some countable index set $I$
is denoted by $M^I$ and we write $\Lambda(i)$ for the $i$-th sequence element of
a sequence $\Lambda\in M^I$.

For topological vector spaces $Y$ and $Z$ the class of linear continuous mappings from
$Y$ to $Z$ is denoted by $\mathcal{L}(Y,Z)$. The notation $\Phi:
Y\hookrightarrow Z$ indicates that $Y$ is continuously embedded into $Z$, i.e.,
$\Phi$ is an injective continuous linear map from $Y$ into $Z$.
If the embedding is canonical we simply write $Y\hookrightarrow Z$.
If $Y$ is equipped with a quasi-norm we use $\|f|Y\|$ for the quasi-norm of
$f\in Y$. The operator quasi-norm of $A \in \mathcal{L}(Y,Z)$ is denoted by
$\|A|Y\to Z\|$.

We use the notation
$a\lesssim b$ if there exists a constant $c>0$ (independent of the
context dependent relevant parameters) such that $a \le c\,b$. If
$a\lesssim b$ and $b \lesssim a$ we write $a \asymp b$. Furthermore,
we write $Y\asymp Z$ for two quasi-normed spaces $Y,Z$ which coincide as sets and whose quasi-norms are equivalent.

\section{General coorbit space theory}
\label{sec:abstrth}


Let $\mathcal{H}$ be a separable Hilbert space 
and $X$ a locally compact Hausdorff space endowed with a positive Radon measure $\mu$ with $\supp \mu = X$.
A family $\mathcal{F} = \{\varphi_x\}_{x\in X}$ of vectors in $\mathcal{H}$
is called a continuous frame (see \cite{alanga93}) if the assignment $x\mapsto\varphi_x$ is weakly measurable and if there exist constants $0<C_1,C_2<\infty$ such that \eqref{eq:stab} is satisfied.
Let us record an important property.
\begin{lemma}\label{lem:total}
Let $\mathcal{F}=\{\varphi_x \}_{x\in X}$ be a continuous frame in $\mathcal{H}$ and $N\subset X$ a set of measure zero.
Then $\{\varphi_x \}_{x\in X\backslash N}$ is total in $\mathcal{H}$.
\end{lemma}
\bproof
Let us put $X^\ast:=X\backslash N$. We have to show that
$V:=\spn \{ \varphi_x:x\in X^\ast \}$ is dense in $\mathcal{H}$. Indeed, using the frame property of $\mathcal{F}$, we can deduce for every $f\perp V$
\[
\|f|\mathcal{H}\|^2 = \int_X |\langle f,\varphi_x \rangle|^2 \,d\mu(x) = \int_{X^*} |\langle f,\varphi_x \rangle|^2 \,d\mu(x) = 0.
\]
\eproof

\noindent
To avoid technicalities, we assume throughout this paper that $X$ is $\sigma$-compact.
We further assume that the continuous frame is Parseval, i.e.\ $C_1=C_2=1$, and note that
-- apart from minor changes -- the theory presented here is valid also for general tight frames where $C_1=C_2$.
It is also possible to develop the theory in the setting
of non-tight frames, where the associated coorbit theory has been worked out
in \cite{fora05} -- at least to a significant extent.

For $0<p<\infty$ we define the Lebesgue space $L_p(X):=L_p(X,\mu)$ as usual by
\[
\|F| L_p(X,\mu) \| := \Big( \int_X |F(x)|^p \,d\mu(x) \Big)^{1/p}<\infty.
\]
A function $F$ belongs to $L_\infty(X):=L_\infty(X,\mu)$ if and only if $F$ is essentially bounded.
The corresponding sequence spaces $\ell_p(I)$ are obtained by choosing $X$ as a countable index set $I$, equipped with the discrete topology and counting measure $\mu$.

Associated to a continuous frame $\mathcal{F}$ is the voice transform $V_{\cf}: \mathcal{H} \to L_2(X,\mu)$ defined by 
$$
V_{\cf}f(x) = \langle f,\varphi_x \rangle\,,\quad f \in \ch, x \in X,
$$
and its adjoint $V^{\ast}_{\cf}:L_2(X,\mu) \to \mathcal{H}$ given in a weak sense by the integral
\begin{align}\label{eq:adjoint}
   V^{\ast}_{\cf} F = \int_X 
   F(y)\varphi_y\,d\mu(y)\,.
\end{align}

Since we assume the frame $\cf$ to be Parseval $V_{\mathcal{F}}$ is an isometry and in particular injective.
The adjoint $V_{\mathcal{F}}^*$ is surjective with $\|V_{\mathcal{F}}^*| L_2\rightarrow\mathcal{H} \|=1$ and the
associated frame operator $S_{\mathcal{F}}:=V_{\mathcal{F}}^{\ast}V_{\mathcal{F}}$ is the identity.
Hence we have
\begin{equation*}
   f = \int_{X} V_{\cf}f(y)\varphi_y\,d\mu(y)\quad\mbox{and}\quad
   V_{\cf}f(x) = \int_{X} V_{\cf}f(y)\langle \varphi_y,\varphi_x\rangle \,d\mu(y)\,.
\end{equation*}
The second identity 
is the crucial reproducing formula $R_{\mathcal{F}}(V_{\mathcal{F}}f)=V_{\mathcal{F}}f$ for $f\in\mathcal{H}$, where
\begin{align}\label{eqdef:kernfunc}
R_{\mathcal{F}}(x,y) = \langle \varphi_y, \varphi_x \rangle\,,\quad x,y \in X,
\end{align}
is an integral kernel (operator), referred to as the \emph{frame kernel} associated to $\mathcal{F}$.
It acts as a self-adjoint bounded operator $R_{\mathcal{F}}=V_{\mathcal{F}}V_{\mathcal{F}}^*:L_2(X)\rightarrow L_2(X)$, which is
an orthogonal projection with $R_{\mathcal{F}}(L_2(X))=V_{\mathcal{F}}(\mathcal{H})$.
The converse of the reproducing formula is also true, i.e.,
if $F\in L_2(X)$ satisfies $R_{\mathcal{F}}(F)=F$ then there exists a unique element $f\in\mathcal{H}$ such that $V_{\mathcal{F}}f=F$.

We remark that we use the same notation for the function $R_{\mathcal{F}}:X\times X\to\com$ given in \eqref{eqdef:kernfunc} and the associated operator $R_{\mathcal{F}}:L_2(X)\rightarrow L_2(X)$.
It is important to note that the function $R_{\mathcal{F}}$ is measurable. Indeed, utilizing an orthonormal basis $(f_n)_{n\in\N}$ of $\mathcal{H}$ we can expand $R_{\mathcal{F}}(x,y)=\sum_{n\in\N} \langle \varphi_y,f_n \rangle \langle f_n,\varphi_x \rangle$ as a point-wise limit of measurable functions.

The idea of coorbit theory is to measure ``smoothness'' of $f$ via properties of the transform $V_{\cf}f$.
Loosely speaking, the coorbit of a function space on $X$ is its retract with respect to (a suitably extended version of) the voice transform.
The classical theory and its generalizations have been developed for the case of certain Banach function spaces on $X$.
In the classical setup, where $X$ is equipped with a group structure, the extension~\cite{ra05-3}
deals with the quasi-Banach case, and our aim is to extend the generalized theory from \cite{fora05,RaUl10} analogously.

\subsection{Function spaces on $X$}
\label{ssec:QBFspaces}

We consider \emph{(quasi)-Banach function spaces}, or shortly \emph{(Q)BF-spaces},
which are linear spaces of measurable functions on $X$, equipped with
a quasi-norm under which they are complete. Hereby, functions are identified when equal almost everywhere. Hence, when speaking of a function
one often actually refers to an equivalence class. In general, this inaccuracy of language does not pose a problem. Only when it comes to point evaluations
the precise meaning must be made clear in the context.

Recall that a quasi-norm on a linear space $Y$ generalizes the concept of a norm by replacing the triangle inequality with
the more general quasi-triangle inequality
\[
\| f + g \| \le C_Y ( \| f \| + \| g \|), \quad f,g\in Y,
\]
with associated quasi-norm constant $C_Y\ge1$. Many aspects of the theory of normed spaces carry over to the quasi-norm setting,
e.g.\ boundedness and continuity coincide, all $d$-dimensional quasi-norms are equivalent, etc..
An important exception is the Hahn-Banach theory concerned with the dual spaces. Note that the (topological) dual $Y^\prime$ of a quasi-normed space $Y$, equipped with the usual operator norm,
is always a Banach space. Due to the possible non-convexity of the quasi-norm however, it may not be sufficiently large
for the Hahn-Banach theorems to hold. In fact, $Y^\prime$ may even be trivial as the example of the $L_p$-spaces in the range $0<p<1$ shows.
This fact poses a serious problem for the theory.

An important tool for dealing with quasi-norms is the
Aoki-Rolewicz theorem~\cite{Ao42,Ro57}, which states that in every quasi-normed space $Y$ there exists an equivalent $r$-norm --
in the sense of an equal topology -- where an $r$-norm, $0<r\le1$, satisfies the $r$-triangle inequality
\[
\| f + g \|^r \le \| f \|^r + \| g \|^r ,\quad f,g\in Y,
\]
and in particular is a quasi-norm with constant $C_Y=2^{1/r-1}$.
The exponent $r=1/(\log_2 C_Y + 1)$ of the equivalent $r$-norm is called the \emph{exponent of $Y$}.

For a viable theory we need to further restrict the class of function spaces.
A quasi-normed function space $Y$ on $X$ is called \emph{solid}, if the following condition is valid,
\begin{align*}
f\text{ $\mu$-measurable},\,g\in Y,\:|f(x)|\le|g(x)|\,a.e. \quad \Rightarrow \quad  f\in Y \text{ and } \|f|Y\|\le\|g|Y\|.
\end{align*}
In a solid space $Y$ we have the equality ${\|}\,|f|\,|Y\|=\| f | Y\|$ for every $f\in Y$. Moreover, there is a useful criterion for a function $f$ to belong to $Y$,
\begin{align*}
f\in Y \quad \Leftrightarrow \quad |f|\in Y \text{ and } f \text{ $\mu$-measurable}. 
\end{align*}
A function space shall be called \emph{rich}, if it contains the characteristic functions $\chi_K$ for all compact subsets $K\subset X$.
A rich solid quasi-normed function space on $X$ then contains the characteristic functions $\chi_U$
for all relatively compact, measurable subsets $U\subset X$.

We will subsequently develop coorbit theory mainly for rich solid QBF-spaces $Y$, that are continuously embedded into $L_1^{\rm loc}(X)$.
As usual, the spaces $L_p^{\rm loc}(X):=L_p^{\rm loc}(X,\mu)$, $0<p\le\infty$, consist of all functions $F$ where $\|F\chi_K | L_p(X)\|<\infty$ for every
compact subset $K\subset X$. The case, where $Y\not\hookrightarrow L_1^{\rm loc}(X)$, is shortly commented on at the end of Subsection~\ref{ssec:coorbit}

It is important to understand the relation between the quasi-norm convergence and the pointwise convergence
of a sequence of functions in $Y$. We have the following result.

\begin{lemma}\label{lem:FuncConv1}
Let $Y$ be a solid quasi-normed function space on $X$, and assume $f_n\rightarrow f$ in $Y$.
Then for arbitrary but fixed representing functions $\widetilde{f}_n,\widetilde{f}$
the following holds true. For a.e.\ $x\in X$ there is a subsequence $(f_{n_k})_{k\in\N}$, whose choice may depend on the particular $x\in X$,
such that $\widetilde{f}_{n_k}(x)\rightarrow \widetilde{f}(x)$ as $k\rightarrow \infty$.
\end{lemma}
\bproof
Assume first that $f_n\rightarrow 0$ in the quasi-norm of $Y$, which implies $\|f_n|Y\|\rightarrow 0$.
As $\inf_{m\ge n} |f_m|$ is a measurable function with $\inf_{m\ge n} |f_m| \le |f_k|$ for all $k\ge n$
we have $\inf_{m\ge n} |f_m|\in Y$ with $\| \inf_{m\ge n} |f_m| | Y \| \le \| f_k |Y\|$ for all
$k\ge n$ by solidity. It follows
$
0\le \| \inf_{m\ge n} |f_m| | Y \| \le \inf_{m\ge n} \| f_m |Y\|=0,
$
and hence $\inf_{m\ge n} |\widetilde{f}_m|(x) =0 $ for a.e.\ $x\in X$. This implies that for these $x\in X$
there is a subsequence $(f_{n_k})_{k\in\N}$ such that $\widetilde{f}_{n_k}(x)\rightarrow 0$.
Now let $f_n\rightarrow f$. Then $(f_n-f)\rightarrow 0$ and by the previous argumentation
for a.e.\ $x\in X$ there is a subsequence $(f_{n_k})_{k\in\N}$ such that $\widetilde{f}_{n_k}(x)-\widetilde{f}(x)\rightarrow 0$, whence  $\widetilde{f}_{n_k}(x)\rightarrow \widetilde{f}(x)$.
\eproof

\begin{remark}
A more thorough investigation of pointwise convergence in solid quasi-normed function spaces is carried out in~\cite{Voigt15}.
It turns out that Lemma~\ref{lem:FuncConv1} can be strengthened using \cite[Cor.~2.2.9]{Voigt15} and the fact that $X$ is $\sigma$-finite (see Step~1 in the proof of Lemma~\ref{lem:Bochner}).
In fact, there is a subsequence $(f_{n_k})_{k\in\N}$, independent of $x\in X$, with $\widetilde{f}_{n_k}(x)\rightarrow \widetilde{f}(x)$ for a.e.\ $x\in X$.
\end{remark}

\subsection{Associated sequence spaces}

Let us take a look at sequence spaces associated with a function space $Y$ on $X$.
For this we recall the notion of an admissible covering introduced in \cite{fora05,RaUl10}.
We say that a covering $\mathcal{U}=\{U_i\}_{i\in I}$ of $X$ is \emph{locally finite} if every $x\in X$ possesses
a neighborhood which intersects only a finite number of the covering sets $U_i$.

\begin{definition}\label{def:admcov}
A covering $\mathcal{U}=\{U_i\}_{i\in I}$ of $X$ is called \emph{admissible}, if it is locally finite and if it satisfies the following conditions:
\begin{enumerate}
\item[(i)]  Each $U_i$ is measurable, relatively compact and has non-void interior. 
\item[(ii)]  The \emph{intersection number} $\sigma(\mathcal{U}):= \sup_{i\in I} \sharp\{ j ~:~ U_i\cap U_j\neq\emptyset \} $ is finite.
\end{enumerate}
\end{definition}

A covering of a locally compact Hausdorff space is
locally finite if and only if every compact subset intersects only a finite number of the covering sets.
Hence, every locally finite covering of the $\sigma$-compact space $X$ is countable.
In particular, the following lemma holds true.
\begin{lemma}\label{lem:admindex}
Every admissible covering of the $\sigma$-compact space $X$ has a countable index set.
\end{lemma}


\noindent
Following \cite{fora05,RaUl10} we now define two types of sequence spaces associated to $Y$.

\begin{definition}
For a rich solid QBF-space $Y$ on $X$ and
an admissible
covering $\mathcal{U} = \{U_i\}_{i\in I}$ of $X$
the sequence spaces $Y^{\flat}$ and $Y^{\natural}$ associated to
$Y$ and $\mathcal{U}$ are defined by
\begin{equation}\nonumber
  \begin{split}
    Y^{\flat} = Y^{\flat}(\mathcal{U}) &:= \Big\{\{\lambda_i\}_{i\in I}~:~
    \|\{\lambda_i\}_{i\in I}|Y^{\flat}\| := \Big\|\sum\limits_{i\in I}
    |\lambda_i|\chi_{U_i}|Y\Big\|<\infty
    \Big\}\,,\\
    Y^{\natural} = Y^{\natural}(\mathcal{U}) &:= \Big\{\{\lambda_i\}_{i\in I}~:~
    \|\{\lambda_i\}_{i\in I}|Y^{\natural}\| := \Big\|\sum\limits_{i\in I}
    |\lambda_i|\mu(U_i)^{-1}\chi_{U_i}|Y\Big\|<\infty
    \Big\}\,.
  \end{split}
\end{equation}
\end{definition}

Note that due to Lemma~\ref{lem:admindex} the index set $I$ of these sequence spaces is necessarily countable.
Also observe that due to condition~(i) of Definition~\ref{def:admcov} and $\supp \mu=X$ we have $\mu(U_i)>0$ for every $i\in I$, and in turn $\|\chi_{U_i}|Y\|>0$.

Viewing a sequence as a function on the index set $I$, equipped with the counting measure,
we subsequently use the terminology introduced above for function spaces. For better distinction, we will speak of a quasi-Banach sequence space and use the abbreviation QBS-space.

\begin{proposition}\label{prop:ss_basic}
The sequence spaces
$Y^\flat(\mathcal{U})$ and $Y^\natural(\mathcal{U})$ are rich solid QBS-spaces with the same quasi-norm constant $C_Y$ as $Y$.
\end{proposition}

Before we give the proof of this proposition let us establish some useful embedding results.
First observe that the mapping
\begin{align}\label{eq:ss_isom}
I^\natural_\flat: Y^\flat\to Y^\natural, \lambda_i\mapsto \mu(U_i)\lambda_i
\end{align}
is an isometric isomorphism between $Y^\flat$ and $Y^\natural$, which allows to transfer statements from one space to the other.
Moreover, if $\inf_{i\in I}\mu(U_i)>0$ we have the embedding $Y^\flat\hookrightarrow Y^\natural$. Analogously, $\sup_{i\in I}\mu(U_i)<\infty$ implies $Y^\natural\hookrightarrow Y^\flat$.
Consequently, $Y^\flat\asymp Y^\natural$ if both conditions are fulfilled.

Let $\nu:I\to[0,\infty)$ be a discrete weight and define $\|\Lambda|\ell^\nu_p\|:=\| \Lambda\nu |\ell_p \|$ for $0<p\le\infty$ and $\Lambda\in\com^I$.
The space $\ell_p^{\nu}(I):=\{ \Lambda\in\com^I : \|\Lambda| \ell^\nu_p\|<\infty \}$
is a QBS-space with quasi-norm $\|\cdot|\ell^\nu_p\|$.


\begin{lemma}\label{lem:ss_embed}
Let $0<p\le1$ be the exponent of $Y$. We then have the continuous embeddings
\begin{align*}
\ell_p^{\omega^\flat}(I) \hookrightarrow Y^\flat(\mathcal{U})\hookrightarrow \ell_\infty^{\omega^\flat}(I) \quad\text{and}\quad
\ell_p^{\omega^\natural}(I) \hookrightarrow Y^\natural(\mathcal{U})\hookrightarrow \ell_\infty^{\omega^\natural}(I)
\end{align*}
with weights defined by $\omega^\flat(i):=\| \chi_{U_i} |Y\|$ and $ \omega^\natural(i):=\mu(U_i)^{-1} \|\chi_{U_i}|Y\| $ for $i\in I$.
\end{lemma}
\bproof
We have
$
\| \{\lambda_i\}_{i\in I} | Y^\flat \|^p
=\big\| \sum_{i\in I} |\lambda_i|\chi_{U_i} \Big| Y \big\|^p
\lesssim  \sum_{i\in I} |\lambda_i|^p \| \chi_{U_i}| Y \|^p
= \| \{\lambda_i\}_{i\in I} | \ell_p^{\omega^\flat} \|^p
$
for $\{\lambda_i\}_{i\in I}\in \ell_p^{\omega^\flat}$.
If $\{\lambda_i\}_{i\in I}\in Y^\flat$ we can estimate for every $j\in I$
\begin{align}\label{eq:ss_eval}
|\lambda_j|\omega^\flat(j) =|\lambda_j| \| \chi_{U_j} | Y \| =
\| |\lambda_j|\chi_{U_j} | Y  \| \le \Big\| \sum_{i\in I} |\lambda_i|\chi_{U_i} \Big| Y \Big\|
= \| \{ \lambda_i\}_{i\in I} | Y^\flat \|.
\end{align}
The embeddings for $Y^\natural$ follow with the isometry~\eqref{eq:ss_isom}.
\eproof

The weights $\omega^\flat$ and $\omega^\natural$ also occur in the following result.

\begin{corollary}
For every $j\in I$ the evaluation $E_j:\{\lambda_i\}_{i\in I}\mapsto \lambda_j$ is a bounded functional on $ Y^\flat$ and $Y^\natural$ with
$\| E_j | Y^\flat\rightarrow\C  \| \le (\omega^\flat(j))^{-1}$ and $\| E_j | Y^\natural\rightarrow\C \| \le (\omega^\natural(j))^{-1}$.
\end{corollary}
\bproof
For $Y^\flat$ this follows directly from \eqref{eq:ss_eval}. The argument for $Y^\natural$ is similar.
\eproof

Now we are ready to give the proof of Proposition~\ref{prop:ss_basic}.

\bproof[Proof of Proposition~\ref{prop:ss_basic}]
We prove the completeness of $Y^\flat$.
The result for $Y^\natural$ follows then with the isometry \eqref{eq:ss_isom}.
A Cauchy sequence $(\Lambda_n)_{n\in\N}$ in $Y^\flat$ is also
a Cauchy sequence in $\ell_\infty^{\omega^\flat}$ by Lemma~\ref{lem:ss_embed}.
Let $\Lambda$ be the limit in $\ell_\infty^{\omega^\flat}$.
We show that $\Lambda\in Y^\flat$ and $\Lambda=\lim_{n\rightarrow\infty} \Lambda_n$ in the quasi-norm
of $Y^\flat$. For this task let us introduce the auxiliary operator $A(\Lambda):=\sum_{i\in I} |\Lambda(i)| \chi_{U_i}$,
which maps $\Lambda\in\C^I$ to a nonnegative measurable function on $X$.

For $\alpha\in\C$ and $\Lambda,\Lambda_1,\Lambda_2 \in\C^I$ we have
$A( \alpha \Lambda )= |\alpha| A(\Lambda)$ and $A(\Lambda_1+\Lambda_2)\le A(\Lambda_1)+A(\Lambda_2)$.
We also have
\begin{align}\label{aux:rel1}
|A(\Lambda_1)-A(\Lambda_2)|\le \sum_{i\in I} \big||\Lambda_1(i)|- |\Lambda_2(i)|\big| \chi_{U_i}
\le \sum_{i\in I} |\Lambda_1(i)- \Lambda_2(i)| \chi_{U_i} = A(\Lambda_1-\Lambda_2).
\end{align}
A sequence $\Lambda\in\C^I$ belongs to $Y^\flat$ if and only if $A(\Lambda)\in Y$, and we have the identity
\begin{align}\label{aux:rel2}
\| \Lambda | Y^\flat \|=  \| A(\Lambda) | Y\|.
\end{align}
Since $\Lambda$ is the limit of $(\Lambda_n)_{n\in\N}$ in $\ell_\infty^{\omega^\flat}$ it holds
$\lim_{n\rightarrow\infty} | \Lambda(i)- \Lambda_n(i)|=0$ for all $i\in I$.
Considering the local finiteness of the sum in the definition of $A$ it follows
that
\begin{align}\label{aux:rel3}
\lim_{n\rightarrow\infty} A(\Lambda-\Lambda_n)(x)=0 \quad\text{ for all }x\in X.
\end{align}
The rest of the proof relies solely on Properties \eqref{aux:rel1}-\eqref{aux:rel3} of the operator $A$ and the solidity and completeness
of $Y$.
First we show $A( \Lambda)\in Y$ which is equivalent to $\Lambda\in Y^\flat$ according to \eqref{aux:rel2}.
The sequence $(A(\Lambda_n))_{n\in\N}$ is a Cauchy sequence in $Y$ because with \eqref{aux:rel1} we can estimate
$
\| A(\Lambda_n) - A(\Lambda_m) | Y \|
\le \| A(\Lambda_n-\Lambda_m) | Y \| = \| \Lambda_n-\Lambda_m | Y^\flat \|.
$
Furthermore, from \eqref{aux:rel3} and \eqref{aux:rel1} it follows
$
\lim_{n\rightarrow\infty} A(\Lambda_n)(x) = A(\Lambda)(x) 
$
for all $x\in X$.
Since $Y$ is complete we can conclude with Lemma~\ref{lem:FuncConv1} that
$A(\Lambda_n)\rightarrow A(\Lambda)$ in $Y$ and $A(\Lambda)\in Y$.
Finally we show $\Lambda=\lim_{n\rightarrow\infty} \Lambda_n$ in $Y^\flat$. The sequence $(A(\Lambda_n-\Lambda))_{n\in\N}$ is a Cauchy sequence in $Y$, because
with \eqref{aux:rel1} we get
\begin{gather*}
\| A(\Lambda_n-\Lambda)-A(\Lambda_m-\Lambda) | Y\|
\le \| A(\Lambda_n-\Lambda_m) | Y \| = \| \Lambda_n-\Lambda_m | Y^\flat \|.
\end{gather*}
Using \eqref{aux:rel3} and Lemma~\ref{lem:FuncConv1} we deduce
$
A(\Lambda_n-\Lambda)\rightarrow 0
$
in $Y$. In view of \eqref{aux:rel2} this finishes the proof.
\eproof

We finally study sequence spaces where the finite sequences are a dense subset.
Since $Y^\flat$ and $Y^\natural$ are isometrically isomorphic via the isometry $I^\natural_\flat$ from \eqref{eq:ss_isom},
and since $I^\natural_\flat$ is a bijection on the sequences with finite support,
these are dense in $Y^\flat$ if and only if they are dense in $Y^\natural$.
The next result occurs in \cite[Thm.~5.2]{fora05} in the context of Banach spaces.
However, the boundedness of the functions required there is not necessary.

\begin{lemma}
If the functions with compact support are dense in $Y$ the finite sequences are dense
in $Y^\flat(\mathcal{U})$ and $Y^\natural(\mathcal{U})$.
\end{lemma}
\bproof
Let $\Lambda=\{\lambda_i\}_{i\in I}\in Y^\flat$ and fix $\varepsilon>0$.
Then $F:=\sum_{i\in I} |\lambda_i|\chi_{U_i} \in Y$ and there exists a function $G\in Y$ with compact support
$K$ such that $\| F- G | Y \|<\varepsilon$. As the covering $\mathcal{U}=\{U_i\}_{i\in I}$ is locally finite, the index set $J:=\{ i\in I : U_i\cap K\neq\emptyset \}$ is finite.
Let $\tilde{\Lambda}$ be the sequence which coincides with $\Lambda$ on $J$ and vanishes elsewhere.
Then $\tilde{F}:=\sum_{i\in J} |\lambda_i|\chi_{U_i}\in Y$ and $|F-\tilde{F}| \le |F - G|$.
Using the solidity of $Y$ we conclude
$
\| \Lambda - \tilde{\Lambda} |Y^\flat \| 
= \| F-\tilde{F} | Y  \| \le \| F- G | Y  \| <\varepsilon.
$
\eproof

For a countably infinite sequence $\Lambda=\{\lambda_i\}_{i\in I}$, a bijection $\sigma:\N\rightarrow I$
and $n\in\N$ we define $\Lambda^\sigma_n$
as the sequence which coincides with $\Lambda$ on $\sigma(\{1,\ldots,n\})$ and is zero elsewhere.

\begin{lemma}\label{lem:ss_findens}
Let $\mathcal{U} = \{U_i\}_{i\in I}$ be an admissible covering and assume that there is a bijection $\sigma:\N\rightarrow I$.
The finite sequences are dense in $Y^\flat(\mathcal{U})$
if and only if for all $\Lambda\in Y^\flat(\mathcal{U})$
it holds $\Lambda^\sigma_n \rightarrow \Lambda$ in the quasi-norm of $Y^\flat(\mathcal{U})$ for $n\rightarrow\infty$.
\end{lemma}
\bproof
Assume that the finite sequences are dense. For $n\in\N$ we can then choose a finite sequence $\Gamma_n\in Y^\flat$ with
$\| \Gamma_n - \Lambda | Y^\flat \|<2^{-n}$. By solidity of $Y$ we get for $N\ge 1 + \max \{ k\in\N |\Gamma_n (\sigma(k))\neq 0 \}$,
with the convention $\max \emptyset =0$, the estimate
$
\| \Lambda^\sigma_{N} - \Lambda | Y^\flat \| \le \| \Gamma_n - \Lambda | Y^\flat \| < 2^{-n}.
$
The other direction is trivial.
\eproof

We end this paragraph with an illustration and examine the sequence spaces associated to the weighted Lebesgue space $L_p^\nu(X)$,
defined by $\|F| L^\nu_p(X) \|:=\| F\nu |  L_p(X) \|<\infty$,
where $\nu$ is a weight and $0<p\le\infty$. In this special case we have a stronger statement than Lemma~\ref{lem:ss_embed}.

\begin{proposition}
Let $\mathcal{U}=\{U_i\}_{i\in I}$ be an admissible covering of $X$, $\nu$ be a
weight and $0< p\le\infty$.
Then for $Y=L_p^\nu(X)$ we have $Y^\flat(\mathcal{U})\asymp\ell_p^{\nu^\flat_p}(I)$ and
$ Y^\natural(\mathcal{U})\asymp\ell_p^{\nu^\natural_p}(I)$ with weights
given by $\nu^\flat_p(i):= \| \chi_{U_i} | L_p^\nu(X) \|$ and $\nu^\natural_p(i):= \mu(U_i)^{-1}\nu^\flat_p(i)$ for $i\in I$.
\end{proposition}
\bproof
We give the proof for $0<p<\infty$ and $Y=L_p^\nu(X)$. For $\{\lambda_i\}_{i\in I}\in\C^I$ we can estimate
\begin{gather*}
\| \{\lambda_i\}_{i\in I} | Y^\flat \|^p = \Big\| \sum_{i\in I} |\lambda_i|\chi_{U_i} \Big| Y \Big\|^p
= \int_X \Big| \sum_{i\in I} |\lambda_i|\chi_{U_i}(y)\nu(y) \Big|^p \,d\mu(y) \\
\asymp \int_X  \sum_{i\in I} |\lambda_i|^p\chi_{U_i}(y)^p \nu(y)^p \,d\mu(y)
=  \sum_{i\in I} |\lambda_i|^p  \int_X \chi_{U_i}(y)^p \nu(y)^p \,d\mu(y)
= \sum_{i\in I} |\lambda_i|^p  \nu^\flat_p(i)^p,
\end{gather*}
where we used that the intersection number $\sigma(\mathcal{U})$ is finite and the equivalence of the $p$-norm and the $1$-norm on $\C^{\sigma(\mathcal{U})}$.
Applying the isometry~\eqref{eq:ss_isom} yields the result for $Y^\natural$.
\eproof


\subsection{Voice transform extension}

For the definition of the coorbit spaces, we need a sufficiently large reservoir for the voice transform. Hence we extend it in this paragraph
following~\cite{fora05}.
For a weight $\nu:X\rightarrow [1,\infty)$ we introduce the space
$
\mathcal{H}_1^\nu := \left\{ f\in\mathcal{H} ~:~ V_\mathcal{F}f\in L_1^\nu(X,\mu) \right\}.
$
Since $\mathcal{F}$ is total in $\mathcal{H}$ by Lemma~\ref{lem:total} it is easy to verify that
$\| f | \mathcal{H}_1^\nu \| :=\| V_{\mathcal{F}}f |  L_1^\nu \|$ constitutes a norm on $\mathcal{H}_1^\nu$.
Further, we define the kernel algebra
\begin{equation*}
  \mathcal{A}_1 := \{K:X \times X \to \C~:~ K \mbox{ is measurable and } \|K|\mathcal{A}_1\| < \infty\},
\end{equation*}
where
\hfill $    \|K|\mathcal{A}_1\| := \max\Big\{\esssup{x\in X}\int_{X}|K(x,y)|d\mu(y)~,~ \esssup{y\in X} \int_{X}|K(x,y)|d\mu(x)\Big\}.$ \hfill~

\vspace*{1ex}
\noindent
Associated to $\nu$ is a weight $m_\nu$ on $X \times X$ given by
\begin{equation*}
     m_\nu(x,y) = \max\Big\{\frac{\nu(x)}{\nu(y)}, \frac{\nu(y)}{\nu(x)} \Big\}\,,\quad x,y\in X.
\end{equation*}
The corresponding sub--algebra $\mathcal{A}_{m_\nu} \subset \mathcal{A}_1$ is defined as
\begin{align}\label{eqdef:Am}
    \mathcal{A}_{m_\nu} := \{K:X\times X \to \com~:~Km_\nu \in \mathcal{A}_1\}
\end{align}
and endowed with the norm $\|K|\mathcal{A}_{m_\nu}\| := \|Km_\nu|\mathcal{A}_1\|$. Note that a kernel $K\in\mathcal{A}_{m_\nu}$ operates continuously
on $L_1^\nu(X)$ and $L_\infty^{1/\nu}(X)$ with $\|K |L_1^\nu(X)\rightarrow L_1^\nu(X) \|,\,\|K |L_\infty^{1/\nu}(X)\rightarrow L_\infty^{1/\nu}(X) \|  \le \|K | \mathcal{A}_{m_\nu} \| $.
Technically, the theory rests upon (mapping) properties of certain kernel functions.
A first example of a typical result is given by the following lemma.

\begin{lemma}\label{auxlem:crossG}
Assume that for a family $\mathcal{G}=\{\psi_x\}_{x\in X}\subset\mathcal{H}$ the Gramian kernel
\begin{align}\label{eq:crossGram}
G[\mathcal{G},\mathcal{F}](x,y):=\langle\varphi_y,\psi_x\rangle \qquad x,y\in X,
\end{align}
is contained in $\mathcal{A}_{m_\nu}$. Then
$\psi_x\in\mathcal{H}_1^\nu$ with
$\| \psi_x | \mathcal{H}_1^\nu \| \le \| G[\mathcal{G},\mathcal{F}] |\mathcal{A}_{m_\nu}\|\nu(x) $ for a.e.\ $x\in X$.
\end{lemma}
\bproof
We have
$
\| G[\mathcal{G},\mathcal{F}]  | \mathcal{A}_{m_\nu}  \|
\ge \int_X |V_{\mathcal{F}}\psi_x(y)| \frac{\nu(y)}{\nu(x)} \,d\mu(y)
= \frac{\| \psi_x | \mathcal{H}_1^\nu \|}{\nu(x)}
$
for a.e.\ $x\in X$.
\eproof

The theory in \cite{fora05,RaUl10} is developed under the global assumption that
$\mathcal{F}$ is uniformly bounded, i.e.\ $\|\varphi_x\|\le C_B$ for all $x\in X$ and some $C_B>0$.
This assumption can be weakened.

\begin{lemma}\label{lem:Bochner}
Let $\nu\ge 1$ be a weight such that the analyzing frame $\mathcal{F}$ satisfies
\vspace*{-2.5ex}
\begin{flalign}\label{eq:framecond}
&\parbox{13cm}{
\begin{enumerate}
\item[(i)] $\| \varphi_x |\mathcal{H} \|\le C_B\nu(x) $ for some constant $C_B>0$ and all $x\in X$,
\item[(ii)] $R_\mathcal{F}\in \mathcal{A}_{m_\nu}$.
\end{enumerate}
}&&
\end{flalign}
\vspace*{-4ex}

\noindent
Then
$\mathcal{H}_1^\nu$
is a Banach space and the canonical embedding $\mathcal{H}_1^\nu \hookrightarrow \mathcal{H}$
is continuous and dense.
Moreover, there is a subset $X^\ast\subset X$
such that $\varphi_x\in \mathcal{H}_1^\nu$ for every $x\in X^\ast$ and $\mu(X\backslash X^\ast)=0$.
The corresponding map
$
\Psi: X^\ast\to\mathcal{H}_1^\nu,\, x\mapsto \varphi_x
$
is Bochner-measurable in $\mathcal{H}_1^\nu$.
\end{lemma}
\bproof
A Cauchy sequence $(f_n)_{n\in\N}\subset \mathcal{H}_1^\nu$ determines a Cauchy sequence $(F_n:=Vf_n)_{n\in\N}$
in $L_1^\nu$, which converges to some $F\in L_1^\nu$. Since
the kernel $R\in\mathcal{A}_{m_\nu}$ operates continuously on $L_1^\nu$, the equality $F_n=R(F_n)$ for $n\in\N$ implies $F=R(F)$. Furthermore, because of
$\| \varphi_x |\mathcal{H} \|\le C_B\nu(x) $ it holds
$
|R(x,y)|\le C_B^2\nu(x)\nu(y)
$
for all $x,y\in X$ and we can deduce
\begin{align*}
|F(x)|=\Big|\int_X R(x,y)F(y) \,d\mu(y) \Big| \le C_B^2\nu(x) \int_X |F(y)|\nu(y) \,d\mu(y)
= C_B^2\nu(x) \| F | L_1^\nu \|.
\end{align*}
This shows $F\in L_\infty^{1/\nu}$, and as $L_\infty^{1/\nu}\cap L_1^\nu \subset L_2$
even $F\in L_2$.
The reproducing formula on $\mathcal{H}$ yields $f\in\mathcal{H}$ with $Vf=F\in L_1^\nu$, which implies $f\in \mathcal{H}_1^\nu$.
Since $\| f_n - f | \mathcal{H}_1^\nu \| = \| F_n - F | L_1^\nu \|$ we obtain $f_n \rightarrow f$ in $\mathcal{H}_1^\nu$. This proves the
completeness.
To prove the continuity of the embedding we observe
$
\| h | \mathcal{H} \|^2 = \| Vh | L_2 \|^2
\le  \| Vh | L_\infty^{1/\nu} \| \| h | \mathcal{H}_1^\nu \|
$
for $h\in\mathcal{H}_1^\nu$. Together with
$
\| Vh | L_\infty^{1/\nu} \| \le \sup_{x\in X} \left\{  \frac{ \| \varphi_x|\mathcal{H}\|}{\nu(x)} \| h|\mathcal{H}\| \right\} \le C_B \| h|\mathcal{H}\|
$,
where $\|\varphi_x | \mathcal{H} \| \le C_B\nu(x)$ was used, the continuity follows.
Due to Lemma~\ref{auxlem:crossG}, applied with $\mathcal{G}=\mathcal{F}$, there is a null-set $N\subset X$ such that $\varphi_x\in \mathcal{H}_1^\nu$ for every $x\in X^\ast:=X\backslash N$.
The density of $\mathcal{H}_1^\nu\hookrightarrow\mathcal{H}$ is thus a consequence of the totality of $\{\varphi_x\}_{x\in X^\ast}$ in $\mathcal{H}$, as stated by Lemma~\ref{lem:total}.
It remains to prove the Bochner-measurability of $\Psi$.
Since $V_\mathcal{F}:\mathcal{H}_1^\nu \to V_\mathcal{F}(\mathcal{H}_1^\nu)$ is an isometric isomorphism,
it suffices to confirm that
\[
\widetilde{\Psi}:=V_\mathcal{F}\circ \Psi: X^\ast\to L_1^\nu(X), x\mapsto V_\mathcal{F}\varphi_{x} 
\]
is Bochner-measurable in $L_1^\nu(X)$. The proof of this is divided into three steps.

\noindent
Step 1: Let us first construct an adequate partition of $X$. Since $\mu$ is a Radon measure, by definition locally finite, all compact subsets of $X$ have finite measure.
As $X$ is assumed to be $\sigma$-compact, the measure $\mu$ is thus $\sigma$-finite.
Hence $X=\bigcup_{n\in\N} L_n$ for certain subsets $L_n\subset X$ of finite measure. By subdividing each of these sets further into $L_{n,m}:=\{x\in L_n : \nu(x)\le m \}$,
disjointifying these subdivided sets, and finally by renumbering the resulting countable family of sets, we obtain a sequence $(K_n)_{n\in\N}$ of pairwise disjoint sets of finite
measure with $X=\bigcup_{n\in\N} K_n$ and such that $\nu(x)\le C_n$ holds for all $x\in K_n$ and suitable constants $C_n>0$.

\noindent
Step 2: We now show that for every $n\in\N$ the function
\[
\widetilde{\Psi}_n:  X^\ast\to L_1^\nu(X), x\mapsto V_\mathcal{F}\varphi_{x}\cdot\chi_{K_n}
\]
is Bochner-measurable in $L_1^\nu(X)$. To this end, let $(f_\ell)_{\ell\in\N}$ be an orthonormal basis of $\mathcal{H}$
with $f_\ell\in \mathcal{H}_1^\nu$ for all $\ell\in\N$. Such a basis exists since $\mathcal{H}$ is separable and $\mathcal{H}_1^\nu$ is a dense subspace of $\mathcal{H}$.
Then we define the functions 
\[
\Phi_{\ell}:=\overline{V_\mathcal{F}f_{\ell}}\in L_1^\nu(X) \quad\text{and}\quad G_{n,\ell}:= V_\mathcal{F}f_{\ell} \cdot \chi_{K_n} \in L_1^\nu(X).
\]
Note that $\Phi_{\ell}(x)=\langle \varphi_x,f_\ell \rangle$ is the $\ell$-th expansion coefficient of $\varphi_x$ with respect to $(f_\ell)_{\ell\in\N}$.
Due to the measurability of $\Phi_{\ell}\in L_1^\nu(X)$ the function $x\mapsto \Phi_{\ell}(x) G_{n,\ell}$ is clearly Bochner-measurable.
Since the pointwise limit of Bochner-measurable functions is again Bochner-measurable, Step~2 is finished if we can show that for every fixed $x\in X^\ast$
\begin{align*}
\widetilde{\Psi}_n(x)=\lim_{N\to\infty} \sum_{\ell=1}^N \Phi_{\ell}(x) G_{n,\ell}  \quad\text{in } L_1^\nu(X).
\end{align*}
This follows with Lebesgue's dominated convergence theorem: For every $y\in X$ we have
\begin{align*}
\lim_{N\to\infty} \sum_{\ell=1}^N \Phi_{\ell}(x) G_{n,\ell}(y) = \lim_{N\to\infty} V_\mathcal{F} \Big( \sum_{\ell=1}^N \Phi_{\ell}(x) f_{\ell} \Big)(y) \cdot \chi_{K_n} (y)
= V_\mathcal{F}\varphi_x(y) \cdot \chi_{K_n}(y) = \widetilde{\Psi}_n(y) .
\end{align*}
Note here that $\varphi_x=\sum_{\ell=1}^\infty \Phi_{\ell}(x)f_\ell$ with convergence in $\mathcal{H}$, and in general $V_\mathcal{F}g_N(x)\to V_\mathcal{F}g(x)$ for fixed $x\in X$ if
$g_N\to g$ in $\mathcal{H}$.
Finally, we estimate using $\|\varphi_x | \mathcal{H} \| \le C_B\nu(x)$
\begin{align*}
\Big|\sum_{\ell=1}^N \Phi_{\ell}(x) G_{n,\ell}(y)\Big| &\le \Big(\sum_{\ell=1}^N |\Phi_{\ell}(x)|^2\Big)^{\frac{1}{2}}  \Big( \sum_{\ell=1}^N |G_{n,\ell}(y)|^2 \Big)^{\frac{1}{2}} \\
&\le \|\varphi_x |\mathcal{H}\| \| \varphi_y |\mathcal{H}\|  \chi_{K_n}(y) \le C_B \nu(y) \|\varphi_x|\mathcal{H}\|  \chi_{K_n}(y) \le C_B C_n \|\varphi_x|\mathcal{H}\|  \chi_{K_n}(y).
\end{align*}
Since $K_n$ is of finite measure this provides an integrable majorant (with respect to $y$).

\noindent
Step 3: Similar to Step~2 the Bochner-measurability of $\widetilde{\Psi}$ is proved by showing for $x\in X^\ast$
\[
\widetilde{\Psi}(x) = \lim_{N\to\infty} \sum_{n=1}^N \widetilde{\Psi}_{n}(x) \quad\text{in } L_1^\nu(X).
\]
The pointwise limit is obvious: For every $y\in X$ we clearly have
\begin{align*}
[\widetilde{\Psi}(x)](y)= V_\mathcal{F}\varphi_{x}(y) =  \lim_{N\to\infty} \sum_{n=1}^N \chi_{K_n}(y) \cdot V_\mathcal{F}\varphi_{x}(y) =  \lim_{N\to\infty} \sum_{n=1}^N [\widetilde{\Psi}_n(x)](y).
\end{align*}
Using Lebesgue's dominated convergence theorem with majorant $|\widetilde{\Psi}(x)|$ proves the claim.

\eproof

\noindent
Under the assumptions~\eqref{eq:framecond} we therefore have the chain of continuous embeddings
\[
\mathcal{H}_1^\nu \overset{i}{\hookrightarrow} \mathcal{H} \overset{\hspace*{+0.5em}i^*}{\hookrightarrow} (\mathcal{H}_1^\nu)^\urcorner,
\]
where $(\mathcal{H}_1^\nu)^\urcorner$ denotes the normed anti-dual of $\mathcal{H}_1^\nu$, which
plays the role of the tempered distributions in this abstract context.
Moreover, there is a subset $X^\ast\subset X$ with $\mu(X\backslash X^\ast)=0$
such that $\varphi_x\in \mathcal{H}_1^\nu$ for $x\in X^\ast$.
%
%
Hence we may extend the transform $V_{\mathcal{F}}:\mathcal{H}\rightarrow L_2(X)$ to $(\mathcal{H}_1^\nu)^\urcorner$ by
\begin{equation}\label{eqdef:Vext}
V_\mathcal{F}f(x) = \langle f,\varphi_x \rangle\,,\quad x\in X^\ast, f\in (\mathcal{H}_1^\nu)^\urcorner,
\end{equation}
where $\langle \cdot,\cdot \rangle$ denotes the duality product on $(\mathcal{H}_1^\nu)^\urcorner \times\mathcal{H}_1^\nu$. The anti-dual is used so that this product extends
the scalar product of $\mathcal{H}$.

\begin{lemma}\label{lem:Vext}
Under the assumptions \eqref{eq:framecond} the extension \eqref{eqdef:Vext} is a well-defined continuous mapping $V_{\mathcal{F}}:(\mathcal{H}_1^\nu)^\urcorner \rightarrow L_\infty^{1/\nu}(X)$.
\end{lemma}
\bproof
Let $f\in (\mathcal{H}_1^\nu)^\urcorner$. The function $V_\mathcal{F}f(x)=\langle f,\varphi_x\rangle$ is well-defined for every $x\in X^\ast$.
It determines a measurable function on $X$, in the sense of equivalence classes, due to the Bochner measurability of $x\mapsto\varphi_x$ in $\mathcal{H}_1^\nu$ proved in Lemma~\ref{lem:Bochner}.
Using Lemma~\ref{auxlem:crossG} we can estimate
\begin{align*}
| V_{\mathcal{F}}f(x) | = |\langle f,\varphi_x \rangle| \le \|f|(\mathcal{H}_1^\nu)^\urcorner\| \|\varphi_x|\mathcal{H}_1^\nu \|
\le \|f|(\mathcal{H}_1^\nu)^\urcorner\| \|\mathcal{R}_\mathcal{F}|\mathcal{A}_{m_\nu}\| \nu(x).
\end{align*}
This shows $V_{\mathcal{F}}f\in L_\infty^{1/\nu}(X)$ with $\|V_{\mathcal{F}}f|L_\infty^{1/\nu}\|\le \|f|(\mathcal{H}_1^\nu)^\urcorner\| \|R_\mathcal{F}|\mathcal{A}_{m_\nu}\|$.

\eproof

\begin{remark}\label{rem:frame}
The membership $R_\mathcal{F}\in\mathcal{A}_{m_\nu}$ does not ensure
$\cf\subset\mathcal{H}_1^\nu$, wherefore the extended voice transform~\eqref{eqdef:Vext} might not be defined at every point $x\in X$.
This detail has not been accounted for in preceding papers, and fortunately it is negligible
since functions on $X$ are only determined up to $\mu$-equivalence classes. Therefore
we -- as in \cite{fora05,RaUl10} -- will henceforth assume $\mathcal{F}\subset\mathcal{H}_1^\nu$ to simplify the exposition.
\end{remark}

We proceed to establish the injectivity of the extended voice transform.
To this end, the following characterization of the duality bracket $\langle \cdot,\cdot \rangle_{(\mathcal{H}_1^\nu)^\urcorner \times \mathcal{H}_1^\nu }$ will be useful.

\begin{lemma}\label{lem:Visom}
If $\mathcal{F}$ has properties \eqref{eq:framecond}, then
for all $f\in(\mathcal{H}_1^\nu)^\urcorner$ and $g\in \mathcal{H}_1^\nu$ it holds
\[
\langle f,g \rangle_{(\mathcal{H}_1^\nu)^\urcorner \times \mathcal{H}_1^\nu } = \int_X V_{\mathcal{F}}f(y)\overline{V_{\mathcal{F}}g(y)} \,d\mu(y)
=: \langle V_{\mathcal{F}}f,V_{\mathcal{F}}g \rangle_{L_\infty^{1/\nu} \times L_1^\nu}.
\]
\end{lemma}
\bproof
Let $f\in(\mathcal{H}_1^\nu)^\urcorner$ and $g\in \mathcal{H}_1^\nu$. Then $V_{\mathcal{F}}g\in L_2\cap L_1^\nu$ and we get
\begin{align*}
\langle f,g \rangle&=
\langle f, V_{\mathcal{F}}^*V_{\mathcal{F}}g \rangle
= \left\langle f, \int_X V_{\mathcal{F}}g(y)\varphi_y \,d\mu(y)  \right\rangle\\
&=\int_X  \overline{V_{\mathcal{F}}g(y)} \langle f,\varphi_y  \rangle \,d\mu(y)
=\langle V_{\mathcal{F}}f,V_{\mathcal{F}}g \rangle_{L_\infty^{1/\nu} \times L_1^\nu}.
\end{align*}
For this equality, it is important that the duality product
commutes with the integral. To verify this, note that since $G:=V_{\mathcal{F}}g \in L_1^\nu$ the
integral $\int_X G(y)\varphi_y \,d\mu(y)$ also exists in the Bochner sense in $\mathcal{H}_1^\nu$.
Indeed, in view of Lemma~\ref{lem:Bochner} the integrand is Bochner-measurable in $\mathcal{H}_1^\nu$. Bochner-integrability follows then from the estimate
\begin{align*}
\int_X |G(y)| \cdot \| \varphi_y | \mathcal{H}_1^\nu \| \,d\mu(y)
\le \|R_\mathcal{F}|\mathcal{A}_{m_\nu}\| \int_X |G(y)| \nu(y)  \,d\mu(y)
= \|R_\mathcal{F}|\mathcal{A}_{m_\nu}\| \| G |  L_1^\nu \|,  
\end{align*}
where Lemma~\ref{auxlem:crossG} was used.
Moreover, the value of the Bochner integral $h:=\int_X G(y)\varphi_y \,d\mu(y)$ equals $g$ since for every $\zeta\in\mathcal{H}$
\begin{align*}
\langle g,\zeta \rangle = \int_X  V_{\mathcal{F}}g(y) \cdot \overline{V_{\mathcal{F}}\zeta(y) } \,d\mu(y) = \langle h,\zeta \rangle.
\end{align*}


\eproof

Using Lemma~\ref{lem:Visom} we can simplify the proof of \cite[Lem.~3.2]{fora05}.

\begin{lemma}\label{lem:Veqivnorm}
Assume that the analyzing frame $\mathcal{F}$ has properties \eqref{eq:framecond}.
Then the expression $\|V_{\mathcal{F}}f | L_\infty^{1/\nu}\|$ is an equivalent norm on $(\mathcal{H}_1^\nu)^\urcorner$.
\end{lemma}
\bproof
We already know from Lemma~\ref{lem:Vext} that
$
\| V_{\mathcal{F}}f | L_\infty^{1/\nu} \| \lesssim \| f | (\mathcal{H}_1^\nu)^\urcorner \|.
$
For the estimate from below we argue with the help of Lemma~\ref{lem:Visom}
\begin{align*}
\| f | (\mathcal{H}_1^\nu)^\urcorner \|
&= \sup_{\| h | \mathcal{H}_1^\nu \|=1} |\langle f,h \rangle_{(\mathcal{H}_1^\nu)^\urcorner \times \mathcal{H}_1^\nu }|
= \sup_{\| h | \mathcal{H}_1^\nu \|=1} |\langle V_{\mathcal{F}}f,V_{\mathcal{F}}h \rangle_{L_\infty^{1/\nu} \times L_1^\nu}| \\
&\le \sup_{H\in L_1^\nu, \|H |L_1^\nu\|\le1} |\langle V_{\mathcal{F}}f,H \rangle_{L_\infty^{1/\nu} \times L_1^\nu}|
= \| V_{\mathcal{F}}f | L_\infty^{1/\nu} \|.
\end{align*}
\eproof

A direct consequence of this lemma is the injectivity of $V_{\mathcal{F}}$.

\begin{corollary}
The voice transform $V_{\mathcal{F}}:\left(\mathcal{H}_1^\nu\right)^\urcorner\rightarrow L_\infty^{1/\nu}(X)$ is continuous and injective.
\end{corollary}

The injectivity of $V_{\mathcal{F}}$ on $(\mathcal{H}_1^\nu)^\urcorner$ implies that $\mathcal{F}$ is total in $\mathcal{H}_1^\nu$.

\begin{corollary}\label{cor:frametotal}
Let $N\subset X$ be a set of measure zero. Then $\{\varphi_x \}_{x\in X\backslash N}$ is total in $\mathcal{H}_1^\nu$.
\end{corollary}
\bproof
If this is not the case, the closure $\mathcal{C}$ of $\spn \{\varphi_x : x\in X\backslash N \}$ in $\mathcal{H}_1^\nu$ is a
true subspace, and the Hahn-Banach extension theorem yields $f\in(\mathcal{H}_1^\nu)^\urcorner$, $f\neq 0$, with
$\langle f,\zeta \rangle=0$ for all $\zeta\in\mathcal{C}$. Hence, $V_{\mathcal{F}}f(x)=0$ for a.e.\ $x\in X$
and therefore $f=0$ by injectivity of $V_{\mathcal{F}}$, which is true even with respect to
$\mu$-equivalence classes in the image space. This is a contradiction.
\eproof

The adjoint $V_{\mathcal{F}}^{\ast}: L_\infty^{1/\nu}(X)\rightarrow (\mathcal{H}_1^\nu)^\urcorner$
of the restriction $V_{\mathcal{F}}:\mathcal{H}_1^\nu \rightarrow L_1^\nu(X)$
naturally extends the adjoint
of $V_{\mathcal{F}}:\mathcal{H}\rightarrow L_2(X)$ due to the equality $\langle F,V_{\mathcal{F}}\zeta \rangle_{L_\infty^{1/\nu}\times L_1^\nu}
= \langle F, V_{\mathcal{F}}\zeta \rangle_{L_2\times L_2}$ in case $\zeta\in \mathcal{H}_1^\nu$ and $F\in L_\infty^{1/\nu}\cap L_2$, and
it can also be represented by a weak integral of the form \eqref{eq:adjoint}.
The relations
\begin{align}\label{eq:reladjoint}
V_{\mathcal{F}}^{\ast}V_{\mathcal{F}}f=f \quad\text{ and }\quad V_{\mathcal{F}}V_{\mathcal{F}}^{\ast}(F)=R(F)
\end{align}
remain valid for the extension, i.e., they hold for $f\in (\mathcal{H}_1^\nu)^\urcorner$ and $F\in L_\infty^{1/\nu}$.
Indeed, Lemma~\ref{lem:Visom} yields
$\langle V_{\mathcal{F}}^{\ast}V_{\mathcal{F}}f ,\zeta  \rangle
= \langle V_{\mathcal{F}}f , V_{\mathcal{F}}\zeta \rangle_{L_\infty^{1/\nu}\times L_1^\nu}
= \langle f,\zeta \rangle$ for all $\zeta\in\mathcal{H}_1^\nu$.
Further, we have
$
V_{\mathcal{F}}V_{\mathcal{F}}^{\ast} F(x)  = \langle V_{\mathcal{F}}^{\ast} F,\varphi_x \rangle
= \int_X F(y) \langle\varphi_y,\varphi_x\rangle \,d\mu(y) = R(F)(x)
$
for all $x\in X$.

An easy consequence of the relations \eqref{eq:reladjoint} is the important fact that the reproducing formula extends to $\left(\mathcal{H}_1^\nu\right)^\urcorner$,
a result obtained differently in \cite[Lemma~3.6]{fora05}.

\begin{lemma}\label{lem:extreproduce}
Let $\nu\ge1$ be a weight on $X$ and assume that the analyzing frame $\mathcal{F}$ satisfies \eqref{eq:framecond}. Then $V_{\mathcal{F}}f(x)=R(V_{\mathcal{F}}f)(x)$
for every $f\in\left(\mathcal{H}_1^\nu\right)^\urcorner$ and $x\in X$.
Conversely, if $F\in L_\infty^{1/\nu}(X)$ satisfies $F=R(F)$ then there is
a unique $f\in \left(\mathcal{H}_1^\nu\right)^\urcorner$ such that $F=V_{\mathcal{F}}f$.
\end{lemma}
\bproof
According to \eqref{eq:reladjoint} we have $R(Vf)=VV^\ast V f =V f$
for $f\in (\mathcal{H}_1^\nu)^\urcorner$. For the opposite direction assume that
$F\in L_\infty^{1/\nu}$ satisfies $F=R(F)$. Then by \eqref{eq:reladjoint} the element $V^*F\in(\mathcal{H}_1^\nu)^\urcorner$
has the property $VV^*F=R(F)=F$. It is unique since $V$ is injective on $(\mathcal{H}_1^\nu)^\urcorner$.
\eproof

Finally we state the correspondence
between the weak*-convergence of a net $(f_i)_{i\in I}$ in $(\mathcal{H}_1^\nu)^\urcorner$
and the pointwise convergence of $(V_{\mathcal{F}}f_i)_{i\in I}$ (compare \cite[Lem.~3.6]{fora05}).

\begin{lemma}
Let $(f_i)_{i\in I}$ be a net in $(\mathcal{H}_1^\nu)^\urcorner$. 
\begin{enumerate}
\item[(i)]
If $(f_i)_{i\in I}$ converges
to some $f\in (\mathcal{H}_1^\nu)^\urcorner$ in the weak*-topology of $(\mathcal{H}_1^\nu)^\urcorner$,
then $(V_{\mathcal{F}}f_i)_{i\in I}$ converges pointwise to $V_{\mathcal{F}}f$ everywhere.
\item[(ii)]
If $(V_{\mathcal{F}}f_i)_{i\in I}$ converges pointwise a.e.\
to a function $F:X\rightarrow \C$ and if
$(f_i)_{i\in I}$ is uniformly bounded in $(\mathcal{H}_1^\nu)^\urcorner$,
then $(f_i)_{i\in I}$ converges to some $f\in(\mathcal{H}_1^\nu)^\urcorner$ in the weak*-topology
with $V_{\mathcal{F}}f=F$ a.e.. 
\end{enumerate}
\end{lemma}
\bproof
We give a proof for sequences $(f_n)_{n\in \N}$ which extends straightforwardly to nets.

\noindent
Part~(i):\,
The weak*-convergence implies $\langle  f_n,\varphi_x \rangle\rightarrow \langle f,\varphi_x \rangle$ for $n\rightarrow\infty$ and all $x\in X$.

\noindent
Part~(ii):\,
Let $X^*\subset X$ denote the subset where the sequence $(V_{\cf}f_n)_{n\in\N}$ converges pointwise.
The space $M=\spn\{ \varphi_x : x\in X^* \}$ lies dense in $\mathcal{H}_1^\nu$ by Corollary~\ref{cor:frametotal}.
We define a conjugate-linear functional $\tilde{f}$ on $M$ by
$\tilde{f}(h):= \lim_{n\rightarrow\infty} \langle f_n,h \rangle$ for $h\in M$.
By assumption, there is $C>0$ so that $\| f_n | (\mathcal{H}_1^\nu)^\urcorner  \|\le C$, which leads to
$
| \langle f_n,h \rangle | \le \|f_n | (\mathcal{H}_1^\nu)^\urcorner \| \|h | \mathcal{H}_1^\nu \|  \le
C \|h | \mathcal{H}_1^\nu \|
$ for all $n\in \N$ and shows that $\tilde{f}$ is bounded on $M$ with respect to $\|\cdot| \mathcal{H}_1^\nu\|$. Hence
it can be uniquely extended to some $f\in(\mathcal{H}_1^\nu)^\urcorner$.
For $\varepsilon>0$ and $\zeta \in \mathcal{H}_1^\nu$ we choose $h\in M$ such that
$\| h-\zeta | \mathcal{H}_1^\nu  \|<\varepsilon$. We get
\[
|\langle f_n - f,\zeta \rangle|
\le \| \zeta-h | \mathcal{H}_1^\nu \| \cdot \| f_n - f | (\mathcal{H}_1^\nu)^\urcorner \| + |\langle f_n - f,h \rangle|
\le \varepsilon ( C+ \| f | (\mathcal{H}_1^\nu)^\urcorner \| ) + |\langle f_n - f,h \rangle|.
\]
Letting $n\rightarrow\infty$ it follows
$\limsup_{n\rightarrow\infty} |\langle f_n - f,\zeta \rangle| \le \varepsilon ( C+ \| f | (\mathcal{H}_1^\nu)^\urcorner \| ) $.
This holds for all $\varepsilon>0$, hence, $\lim_{n\rightarrow\infty} |\langle f_n - f ,\zeta \rangle| =0 $.
This shows that $f_n \rightarrow f$
in the weak*-topology
of $(\mathcal{H}_1^\nu)^\urcorner$. As a consequence
$V_{\mathcal{F}}f(x)=  \langle f,\varphi_x \rangle =\lim_{n\rightarrow\infty} \langle  f_n,\varphi_x \rangle
= \lim_{n\rightarrow\infty} V_{\mathcal{F}}f_n(x) = F(x)$ for all $x\in X^*$.
\eproof

A direct implication is the correspondence principle with respect to sums formulated below.
\begin{corollary}\label{cor:corrpri}
If $\sum_{i\in I} f_i$ converges
unconditionally in the weak*-topology of $(\mathcal{H}_1^\nu)^\urcorner$ then the series $\sum_{i\in I} V_{\mathcal{F}}f_i(x)$  converges absolutely
for all $x\in X$.
Conversely, if $\sum_{i\in I} V_{\mathcal{F}}f_i(x)$ converges absolutely for a.e.\ $x\in X$ and if
the finite partial sums of $\sum_{i\in I} f_i$ are uniformly bounded in $(\mathcal{H}_1^\nu)^\urcorner$ then
$\sum_{i\in I} f_i$ converges unconditionally in the weak*-topology.
\end{corollary}

\subsection{Coorbit spaces}
\label{ssec:coorbit}

In this central part we introduce the notion of coorbit spaces, building upon the correspondence between elements of $(\mathcal{H}_1^\nu)^\urcorner$ and
functions on $X$ as established by the transform $V_\mathcal{F}$. The idea is to characterize $f\in(\mathcal{H}_1^\nu)^\urcorner$
by properties of the corresponding function $V_\mathcal{F}f$.
For a viable theory
the analyzing frame $\mathcal{F}$ must fulfill certain suitability conditions with respect to $Y$.

\begin{definition}\label{def:F(v,Y)}
Let $\nu\ge1$ be a weight on $X$. We say that $\mathcal{F}$ has \emph{property $F(\nu,Y)$}
if it satisfies condition \eqref{eq:framecond} and if the following holds true,
\begin{enumerate}
\item[(i)] $R_{\cf}:Y\rightarrow Y$ acts continuously on $Y$, 
\item[(ii)] $R_{\cf}(Y)\hookrightarrow L_\infty^{1/\nu}(X)$.   
\end{enumerate}
\end{definition}

\noindent
Condition~\eqref{eq:framecond} ensures that the voice transform extends to $(\mathcal{H}_1^\nu)^{\urcorner}$.
Further, conditions (i) and (ii) imply that $R_{\cf}F(x)=\int_X R_{\cf}(x,y)F(y) \,d\mu(y)$ is well-defined for a.e.\ $x\in X$ if $F\in Y$.
In addition, also due to~(i) and~(ii),
the operator $R_{\cf}:Y\rightarrow L_\infty^{1/\nu}(X)$ is continuous: 
For $F\in Y$ we have $R(F)\in L_\infty^{1/\nu}(X)$ and
\[
\| R(F) | L_\infty^{1/\nu}\| \lesssim \| R(F) | Y \| \le  \|R|Y\to Y\|\cdot \| F|Y\| .
\]

In view of Definition~\ref{def:F(v,Y)} it makes sense to introduce the following subalgebra of $\mathcal{A}_{m_\nu}$ from \eqref{eqdef:Am}
\begin{align*}
    \mathcal{B}_{Y,{m_\nu}} = \{K:X\times X \to \C~:~K\in \mathcal{A}_{m_\nu}\mbox{ and }K\mbox{ is bounded from }Y \to Y\}\,,
\end{align*}
equipped with the quasi-norm $\|K|\mathcal{B}_{Y,{m_\nu}}\| := \max\{\|K|\mathcal{A}_{m_\nu}\|, \|K|Y\to Y\|\}$.

Now we are able to give the definition of the coorbit of a rich solid QBF-space $Y$.

\begin{definition}
Let $Y$ be a rich solid QBF-space on $X$ 
and assume that the analyzing frame $\cf = \{\varphi_x\}_{x\in X}$ has property $F(\nu,Y)$ for some weight $\nu:X\rightarrow[1,\infty)$.
The coorbit of $Y$ with respect to $\mathcal{F}$ is defined by
  $$
      \Co(\nu, \cf, Y):= \{f\in (\mathcal{H}_1^\nu)^{\urcorner}~:~V_{\mathcal{F}} f \in Y\}\quad\mbox{with quasi-norm}\quad
      \|f|\Co(\nu, \cf, Y)\| := \|V_{\cf} f|Y\|\,.
  $$
\end{definition}
Since the coorbit is independent of the weight $\nu$ in the definition, as proved by the lemma below, it is omitted in the notation and we simply write
$
\Co(\mathcal{F},Y):=\Co(\nu,\mathcal{F},Y).
$
Moreover, if the analyzing frame $\mathcal{F}$ is fixed we may just write $\Co(Y)$.

\begin{lemma}\label{lem:co_indiw}
The coorbit $\Co(\nu,\mathcal{F},Y)$ does not depend on the particular weight $\nu$
chosen in the definition in the following sense.
If $\tilde{\nu}\ge1$ is another weight such that $\mathcal{F}$ has property $F(\tilde{\nu},Y)$
then we have
$
\Co(\tilde{\nu},\mathcal{F},Y)= \Co(\nu,\mathcal{F},Y).
$
\end{lemma}
\bproof
If $\mathcal{F}$ has properties $F(\nu,Y)$ and $F(\tilde{\nu},Y)$
it also has property $F(\omega,Y)$ for $\omega=\nu+\tilde{\nu}$.
We show $\Co(\omega,Y) = \Co(\nu,Y)$.
Since $\omega\ge \nu$ we have the continuous dense embedding $\mathcal{H}_1^{\omega}\hookrightarrow \mathcal{H}_1^\nu$
which implies $(\mathcal{H}_1^\nu)^\urcorner \hookrightarrow (\mathcal{H}_1^{\omega})^\urcorner$
and hence $\Co(\nu,Y)\subset \Co(\omega,Y)$. For the opposite inclusion
let $f\in \Co(\omega,Y)$. Then $f\in (\mathcal{H}_1^{\omega})^\urcorner$ and $F:=Vf\in Y$ with
$R(F) \in L_\infty^{1/\nu}$ by property $F(\nu,Y)$. Since $F=R(F)$
according to the reproducing formula on $(\mathcal{H}_1^{\omega})^\urcorner$ we thus have $F\in L_\infty^{1/\nu}$.
The inverse reproducing formula on $(\mathcal{H}_1^\nu)^\urcorner$ then yields
$\tilde{f}\in (\mathcal{H}_1^\nu)^\urcorner \subset (\mathcal{H}_1^{\omega})^\urcorner$
with $V\tilde{f}=F$, which due to the injectivity
of $V$ is equal to $f$. This shows $f\in (\mathcal{H}_1^\nu)^\urcorner$,
and as $Vf\in Y$ even
$f\in \Co(\nu,Y)$. Finally note that
the quasi-norms on $\Co(\omega,Y)$ and $\Co(\nu,Y)$ are equal.
Analogously it follows $\Co(\omega,Y)=\Co(\tilde{\nu},Y)$.
\eproof
\begin{remark}
The claim of Lemma \ref{lem:co_indiw} has to be understood in the sense
\begin{align*}
\left\{f|_{\langle\varphi_x :x\in X\rangle}:f\in \Co(\tilde{\nu},\mathcal{F},Y)\right\}=\left\{f|_{\langle\varphi_x :x\in X\rangle}:f\in \Co({\nu},\mathcal{F},Y)\right\}
\end{align*}
since the two spaces are not strictly speaking equal. Further, the span $\langle\varphi_x :x\in X\rangle$ is dense in $\mathcal{H}_1^\nu$ and $\mathcal{H}_1^{\tilde{\nu}}$, thus the notation  $\Co(\tilde{\nu},\mathcal{F},Y)=\Co(\nu,\mathcal{F},Y)$ is justified.
\end{remark}

Regarding the applicability of the theory, it is important to decide whether a given analyzing frame $\cf=\{ \varphi_x \}_{x\in X}$ has property $F(\nu,Y)$.
In the classical theory, where $X$ is a group, the frame is of the special form $\varphi_x=\pi(x)g$, where $\pi$ is a group representation
and $g\in\mathcal{H}$ a suitable vector. In this case properties of $\cf$ break down to properties of the analyzing vector $g$,
and it suffices to check admissibility of $g$, see \cite{FeGr86,feGr89a,Gr91}.
For the continuous wavelet transform concrete conditions can be formulated in terms of smoothness, decay and vanishing moments,
generalized in \cite{furaitou15} to wavelets over general dilation groups.
In our general setup the algebras $\mathcal{A}_{m_\nu}$ and $\mathcal{B}_{Y,{m_\nu}}$ embody the concept of admissibility and
for the (inhomogeneous) wavelet transform utilized in Section~\ref{sec:appcoorbit} also concrete conditions can be deduced, see e.g.\ \cite{RaUl10}.

Concerning the independence of $\Co (\cf,Y)$
on the reservoir $(\mathcal{H}_1^\nu)^{\urcorner}$ we state \cite[Lem.~3.7]{RaUl10}, whose proof carries over directly.

\begin{lemma}\label{lem:co_indires} Assume that the analyzing frame $\cf$ satisfies $F(\nu,Y)$ and let $S$ be a topological vector
space such that $\cf\subset S \hookrightarrow \mathcal{H}_1^\nu$.
In case $\cf$ is total in $S$ and the reproducing formula $V_{\cf}f = R_{\cf}(V_{\cf}f)$ extends to all $f\in S^{\urcorner}$ (the topological anti-dual
of $S$) then
$$
    \Co(\mathcal{F},Y) = \{f\in S^{\urcorner}~:~V_{\cf}f \in Y\}\,.
$$
\end{lemma}

We have the following result concerning the coincidence of the two
spaces $\Co(\cf,Y)$ and $\Co(\cg,Y)$, where $\cf$ and $\cg$ are two different
continuous frames. 

\begin{lemma}
Assume that the frames $\cg = \{g_x\}_{x\in X}$ and $\cf = \{f_x\}_{x\in X}$
satisfy $F(\nu,Y)$. If the Gramian kernels $G[\cf, \cg]$ and $G[\cg, \cf]$ defined in \eqref{eq:crossGram}
are both contained in $\mathcal{B}_{Y,m_\nu}$, we have
$
    \Co(\cf,Y) = \Co(\cg, Y)\,
$
in the sense of equivalent quasi-norms.
\end{lemma}
\bproof
This is a consequence of the relations $V_{\mathcal{F}}=G[\cf,\cg] V_{\mathcal{G}}$ and $V_{\mathcal{G}}=G[\cg,\cf] V_{\mathcal{F}}$.
In view of Lemma~\ref{lem:Bochner} and Lemma~\ref{lem:co_examples} we have $V_{\mathcal{G}}f_x\in  L_1^\nu(X)$ for a.e.\ $x\in X$.
Further, $V_{\mathcal{G}}f\in L_\infty^{1/\nu}(X)$ for $f\in (\mathcal{H}_1^\nu)^\urcorner$ and hence with Lemma~\ref{lem:Visom}
\[
V_{\mathcal{F}}f(x)=\langle f,f_x \rangle_{(\mathcal{H}_1^\nu)^\urcorner \times \mathcal{H}_1^\nu } = \langle V_{\mathcal{G}}f,V_{\mathcal{G}}f_x \rangle_{L_\infty^{1/\nu} \times L_1^\nu}
=\int_X \langle f,g_y \rangle \overline{\langle g_y,f_x \rangle}\,d\mu(y) = G[\cf,\cg] V_{\mathcal{G}}f(x).
\]
This proves $V_{\mathcal{F}}=G[\cf,\cg] V_{\mathcal{G}}$, and by symmetry also $V_{\mathcal{G}}=G[\cg,\cf] V_{\mathcal{F}}$.
\eproof

It is essential for the theory that the reproducing formula carries over to $\Co(Y)$, which
is an immediate consequence of Lemma~\ref{lem:extreproduce}.

\begin{lemma}\label{lem:reproform3}
A function $F\in Y$ is of the form $Vf$ for some $f\in \Co(Y)$
if and only if $F=R(F)$.
\end{lemma}

The reproducing formula is the key to prove the main theorem of this section, which corresponds to
\cite[Prop.~3.7]{fora05}. We explicitly state the continuitiy of the embedding $\Co(Y)\hookrightarrow (\mathcal{H}_1^\nu)^\urcorner$.

\begin{Theorem}\label{thm:co_main}
\begin{enumerate}
\item[(i)] The space $(\Co(Y),\|\cdot|\Co(Y)\|)$ is a quasi-Banach space with quasi-norm constant $C_Y$, which is continuously embedded
into $(\mathcal{H}_1^\nu)^\urcorner$.
\item[(ii)] The map $V:\Co(Y)\rightarrow Y$ establishes
an isometric isomorphism between $\Co(Y)$ and the closed subspace
$R(Y)$ of $Y$.
\item[(iii)] The map $R:Y\rightarrow Y$
is a projection of $Y$ onto $R(Y)=V(\Co(Y))$.
\end{enumerate}
\end{Theorem}
\bproof
In general, we refer to the proof of \cite[Prop.~3.7]{fora05}.
However, the continuity of the embedding $\Co(Y)\hookrightarrow (\mathcal{H}_1^\nu)^\urcorner$ is not proved there.
It is a consequence of the following estimate for $f\in \Co(Y)$, where Lemma~\ref{lem:Veqivnorm} is used,
\begin{align*}
\| f | (\mathcal{H}_1^\nu)^\urcorner \| \asymp \| Vf | L_\infty^{1/\nu}  \|
\le \| R | Y\rightarrow L_\infty^{1/\nu} \| \| Vf | Y \|
= \| R | Y\rightarrow L_\infty^{1/\nu} \| \| f | \Co(Y) \|.
\end{align*}
Further, the proof of \cite[Prop.~3.7]{fora05} implicitly relies on the validity of $R\circ R=R$ on $Y$, which a-priori is only clear for $L_2(X)$.
Therefore, we include a proof of this relation here. Let $F\in Y$ and choose compact subsets $(K_n)_{n\in\N}$ with $X=\bigcup_{n\in\N} K_n$ and $K_n\subset K_m$ for $n\le m$, which is possible
since $X$ is $\sigma$-compact. Then we define the sets $U_n:=\{ x\in K_n : |F(x)|\le n \}$, which are
relatively compact and thus of finite measure. As a consequence,
$F_n:= \chi_{U_n} F \in L_2(X)$. Moreover, $F_n\in Y$ since $|F_n(x)|\le |F(x)|$ for every $x\in X$.
Since by assumption $R:Y\rightarrow Y$ is well-defined the assignment $y\mapsto |R(x,y) F(y)|$ is integrable for a.e.\ $x\in X$. As $F_n(y)\rightarrow F(y)$ pointwise, Lebesgue's dominated convergence theorem thus
yields for these $x\in X$
\begin{align*}
RF_n(x) = \int_X R(x,y) F_n(y) \,d\mu(y) \rightarrow \int_X R(x,y) F(y) \,d\mu(y) = RF(x).
\end{align*}
Next, observe that the function $|R(x,\cdot)|m_{\nu}(x,\cdot)$ is integrable for a.e.\ $x\in X$ since $R\in\mathcal{A}_{m_\nu}$.
Further, due to $R(Y)\hookrightarrow L_\infty^{1/\nu}(X)$ the following estimate holds true for a.e.\ $x,y\in X$
\[
|R(x,y) RF_n(y)| \le C |R(x,y)| \nu(y) \|F_n|Y\| \le  C |R(x,y)| m_{\nu}(x,y) \nu(x) \|F|Y\|.
\]
Another application of Lebesgue's dominated convergence therefore yields for a.e.\ $x\in X$
\begin{align*}
R(RF_n)(x) = \int_X R(x,y) RF_n(y) \,d\mu(y) \rightarrow \int_X R(x,y) RF(y) \,d\mu(y) = R(RF)(x).
\end{align*}
Since $F_n\in L_2(X)$ we have $RF_n=R(RF_n)$ for every $n\in\N$. Altogether, we obtain
\begin{align*}
R(RF)(x) \leftarrow  R(RF_n)(x) = RF_n(x) \rightarrow RF(x).
\end{align*}

\eproof

Let us finally provide some trivial examples, also given in \cite[Cor.~3.8]{fora05}.

\begin{lemma}\label{lem:co_examples}
If the analyzing frame $\mathcal{F}$ satisfies condition~\eqref{eq:framecond} for a weight $\nu\ge 1$,
it has properties $F(\nu,L_2)$, $F(\nu,L_\infty^{1/\nu})$, $F(\nu,L_1^\nu)$, and
it holds
\begin{align*}
(\mathcal{H}_1^\nu)^\urcorner \asymp \Co(\mathcal{F},L_\infty^{1/\nu}),
&& \mathcal{H}_1^\nu = \Co(\mathcal{F},L_1^\nu), && \mathcal{H} = \Co(\mathcal{F},L_2).
\end{align*}
\end{lemma}

Typically, the theory cannot be applied if the QBF-space $Y$ is not embedded in $L_1^{\rm loc}(X)$, since then the kernel conditions concerning operations on $Y$
can usually not be fulfilled. Let us close this paragraph with a short discussion of how to proceed in case $Y\not\hookrightarrow L_1^{\rm loc}(X)$.

\subsubsection*{The case $Y\not\hookrightarrow L_1^{\rm loc}(X)$}

The main idea is to replace $Y$ with a suitable subspace $Z$, which is embedded into $L_1^{\rm loc}(X)$ and fits into the existing theory.
The basic observation behind this is that not all the information of $Y$ is used in the definition of the coorbit. In fact, the information about $\Co(Y)$ is fully contained
in the subspace $R(Y)$, i.e., we have $\Co(\mathcal{F},Y)=\Co(\mathcal{F},R(Y))$.
Thus, we can painlessly pass over to a solid subspace $Z$ of $Y$ and regain the same coorbit if
\begin{align*}
R(Y)\hookrightarrow Z \hookrightarrow Y.
\end{align*}

This observation motivates the idea to substitute $Y$ --  in case $Y$ is not embedded into $L_1^{\rm loc}(X)$ itself -- by a
suitable subspace $Z$ of $Y$ consisting of locally integrable functions, and then to consider the coorbit of $Z$ instead.
In the classical group setting~\cite{ra05-3} Wiener amalgams~\cite{Fe83,ra05-4} were used as suitable substitutes.
Since Wiener amalgams rely on the underlying group structure, they cannot be used in our general setup however.
Instead, it is possible to resort to the closely related
decomposition spaces due to Feichtinger and Gröbner~\cite{fegr85}, which can be viewed as discrete analoga of Wiener amalgams.
This approach has been worked out in \cite{Sch12}, where the decomposition space $\mathcal{D}(Y,\mathcal{U})$ with local component $L_\infty$ and global component $Y$
is used. It is defined as follows.

\begin{definition}[\cite{Sch12}]
The decomposition space $\mathcal{D}(Y,\mathcal{U})$ associated to a rich solid QBF-space $Y$ on $X$ and an admissible
covering $\mathcal{U}=\{U_i\}_{i\in I}$ of $X$ is defined by
\begin{equation}\nonumber
  \begin{split}
    \mathcal{D}(Y,\mathcal{U}) &:= \Big\{ f\in L_\infty^{\rm loc}(X) ~:~
    \| f |\mathcal{D}(Y,\mathcal{U})\| := \Big\|\sum\limits_{i\in I}
    \|f|L_\infty(U_i)\| \chi_{U_i} |Y \Big\|<\infty
    \Big\}.
  \end{split}
\end{equation}
Note that the sum $\sum_{i\in I} \|f\|_{L_\infty(U_i)}\chi_{U_i}$
is locally finite and defines pointwise a function on $X$.
\end{definition}

The space $\mathcal{D}(Y,\mathcal{U})$ is a subspace of $Y$, continuously embedded, and a rich solid QBF-space with the same quasi-norm constant $C_Y$ as $Y$.
Moreover, it is contained in $L_1^{\rm loc}(X)$ even if $Y$ itself is not.
In fact, we have the embedding $\mathcal{D}(Y,\mathcal{U}) \hookrightarrow L_\infty^{1/\omega}(X)$,
where $\omega:X\rightarrow(0,\infty)$ defined by $\omega(x):= \max_{i:x\in U_i} \{\|\chi_{U_i} | Y \|^{-1} \}$ is a locally bounded weight.
For a short proof, let $K\subset X$ be compact and $\{U_i\}_{i\in J}$ the finite subfamily of sets in $\mathcal{U}$ intersecting $K$. Then $\omega(x)\le \max_{i\in J} \{\|\chi_{U_i} | Y \|^{-1} \}$
for all $x\in K$.

In the spirit of \cite[Def.~4.1]{ra05-3}, we may therefore pass over to $\Co(\mathcal{F},\mathcal{D}(Y,\mathcal{U}))$, the coorbit of $\mathcal{D}(Y,\mathcal{U})$.
%
In general, one can only expect $\Co(\mathcal{F},\mathcal{D}(Y,\mathcal{U}))\subset \{ f\in (\mathcal{H}_1^\nu)^\urcorner : V_\mathcal{F}f\in Y \}$ and not equality.
In many applications however the equality can be proved by methods not available in the abstract setting.
Moreover, the choice $\mathcal{D}(Y,\mathcal{U})$ is consistent with the theory due to the result below, which is 
analogous to a result obtained for Wiener amalgams \cite[Thm.~6.1]{ra05-3}.

\begin{Theorem}[{\cite[Thm.~8.1]{Sch12}}]
Assume that $Y$ is a rich solid QBF-space
and that the analyzing frame $\mathcal{F}$ has property $F(\nu,Y)$.
If $\mathcal{U}$ is an admissible covering of $X$ such that the kernel $M^*_{\mathcal{U}}=K^*_\mathcal{U}[\mathcal{F},\mathcal{F}]$ (defined in \eqref{eqdef:kerK} below) operates continuously
on $Y$, then the frame $\mathcal{F}$ has property $F(\nu,\mathcal{D}(Y,\mathcal{U}))$
and it holds
\[
\Co(\mathcal{F},\mathcal{D}(Y,\mathcal{U})) \asymp \Co(\mathcal{F},Y)
\]
in the sense of equivalent quasi-norms.
\end{Theorem}

\begin{remark}
The condition that the kernel $M^*_{\mathcal{U}}$ operates continuously on $Y$
is fulfilled for instance in the important case when $\mathcal{F}$ has property $D(\delta,\nu,Y)$ (see Definition~\ref{def:D(d,v,Y)} below).
\end{remark}

In \cite[Thm.~8.1]{Sch12} this theorem was formulated under the additional assumption that $Y$ is continuously embedded into $L_1^{\rm loc}(X)$.
However, essential for the proof is only that the frame $\mathcal{F}$ has property $F(\nu,Y)$, wherefore we chose to omit this assumption here.

\subsection{Discretizations}

A main feature of coorbit space theory is its general abstract discretization
machinery. With a coorbit characterization of a given function space at
hand, the abstract framework (Theorems~\ref{thm:atomicdec}
and~\ref{thm:frameexp} below)  provides atomic decompositions of this space,
i.e.,\ a representation of functions using ``only'' a countable number of atoms
as building blocks.

Moreover, the function space can be characterized via an equivalent
quasi-norm on an associated sequence space.

The transition to sequence spaces bears many advantages, since those usually
have a simpler, more accessible structure than the original
spaces. For example, the investigation of embedding relations becomes much simpler
by performing them on the associated sequence spaces. In addition, atomic decompositions naturally
lend themselves to real world representations of the considered functions:
By truncation one obtains approximate expansions consisting only of a finite
number of atoms.

Our discretization results, Theorem~\ref{thm:atomicdec} and
Theorem~\ref{thm:frameexp}, transfer the results from \cite{RaUl10},
namely Theorem~3.11 and Theorem~3.14, to the general quasi-Banach setting.
Applying a different strategy for their proofs, however, we are able to
strengthen these results significantly even in the Banach space setting.

\subsubsection*{Preliminaries}

Let us introduce the kernel functions $K_\mathcal{U}[\mathcal{G},\mathcal{F}]$ and $K^*_\mathcal{U}[\mathcal{G},\mathcal{F}]$, which are related by involution
and play a prominent role in the discretization theory. For a family $\mathcal{G}=\{\psi_x\}_{x\in X}$
and an admissible covering $\mathcal{U}=\{U_i\}_{i\in I}$ they are defined by
\begin{align}\label{eqdef:kerK}
K_\mathcal{U}[\mathcal{G},\mathcal{F}](x,y):=\sup_{z\in Q_y} |\langle \varphi_x,\psi_z \rangle| \quad\text{and}\quad
K^*_\mathcal{U}[\mathcal{G},\mathcal{F}](x,y):= K_\mathcal{U}[\mathcal{G},\mathcal{F}](y,x)
\end{align}
where $x,y\in X$ and $Q_y:=\bigcup\limits_{i~:~y\in U_i} U_i$ for $y\in X$.
Their mapping properties are essential for two central results, namely Lemmas~\ref{auxlem:mainanalysis} and \ref{auxlem:mainsynthesis2},
which together with Lemma~\ref{auxlem:Uinvert} provide the technical foundation for the proofs of Theorem~\ref{thm:atomicdec} and Theorem~\ref{thm:frameexp}.

We will subsequently use the symbol $\interleave\cdot\interleave$ for the operator quasi-norm $\|\cdot |Y\rightarrow Y \|$ on $Y$.

\begin{lemma}\label{auxlem:mainanalysis}
Let $Y$ be a rich solid QBF-space on $X$ 
and let the analyzing frame $\mathcal{F}=\{\varphi_x\}_{x\in X}$ possess property $F(\nu,Y)$.
Further, let $\mathcal{G}=\{\psi_x \}_{x\in X}\subset \mathcal{H}_1^\nu$ be a
family and $\mathcal{U}=\{U_i\}_{i\in I}$ an admissible covering such that $K^*_\mathcal{U}:=K^*_\mathcal{U}[\mathcal{G},\mathcal{F}]$ defines a bounded
operator on $Y$.
Then for $f\in \Co(\mathcal{F},Y)$ the function $\sum_{i\in I} \sup_{z\in U_i} |V_\mathcal{G}f(z)|\chi_{U_i}$
belongs to $Y$ with the estimate
\begin{gather*}
\Big\|\sum_{i\in I} \sup_{z\in U_i} |V_\mathcal{G}f(z)|\chi_{U_i} \Big| Y \Big\| \le
\sigma(\mathcal{U})\interleave K^*_\mathcal{U} \interleave \| f | \Co(\mathcal{F},Y)\| \,.
\end{gather*}
\end{lemma}
\noindent
Note that the sum 
$\sum_{i\in I} \sup_{z\in U_i}|V_\mathcal{G}f(z)|\chi_{U_i}$ is locally finite and defined pointwise.

\bproof
Using $V_{\mathcal{G}}f= G[\mathcal{G},\mathcal{F}]V_\mathcal{F}f$
we can estimate for $f\in \Co(\mathcal{F},Y)$ and all $x\in X$
\begin{align*}
\sup_{z\in Q_x} |V_\mathcal{G}f(z)| &= \sup_{z\in Q_x} |G[\mathcal{G},\mathcal{F}]V_\mathcal{F}f(z)|
\le \sup_{z\in Q_x} \int |G[\mathcal{G},\mathcal{F}](z,y)||V_\mathcal{F}f(y)|\,d\mu(y) \\
&\le  \int \sup_{z\in Q_x}|G[\mathcal{G},\mathcal{F}](z,y)||V_\mathcal{F}f(y)|\,d\mu(y)\\
&= \int K^*_\mathcal{U}[\mathcal{G},\mathcal{F}](x,y)|V_\mathcal{F}f(y)|\,d\mu(y) = K^*_\mathcal{U}[\mathcal{G},\mathcal{F}](|V_\mathcal{F}f|)(x).
\end{align*}
For functions $F:X\rightarrow\C$ we further have the estimate
\begin{align}\label{ineq:aux}
\sup_{z\in Q_x} |F(z)| \le \sum_{i\in I} \sup_{z\in U_i} |F(z)| \chi_{U_i}(x) \le \sigma(\mathcal{U}) \sup_{z\in Q_x} |F(z)|,
\end{align}
where $\sigma(\mathcal{U})$ is the intersection number of $\mathcal{U}$. 
Choosing $F=V_\mathcal{G}f$ in \eqref{ineq:aux}, we can conclude
\begin{align*}
\Big\| \sum_{i\in I} \sup_{z\in U_i}|V_\mathcal{G}f(z)|\chi_{U_i} \Big| Y \Big\|
\le \sigma(\mathcal{U}) \interleave K^*_\mathcal{U} \interleave \| V_\mathcal{F}f| Y \|
= \sigma(\mathcal{U}) \interleave K^*_\mathcal{U} \interleave \| f |
\Co(\mathcal{F},Y) \|. ~\hfill\hfill
\end{align*}
\eproof

We can immediately deduce an important result, which corresponds to \cite[Lemma~3.12]{RaUl10}, concerning the sampling of $V_\mathcal{G}f$.

\begin{corollary}\label{auxcor:trafosampling}
With the same assumptions as in the previous lemma let $\{x_i\}_{i\in I}$ be a family of points
such that $x_i\in U_i$.
Then  $\{V_\mathcal{G}f(x_i)\}_{i\in I}\in Y^\flat(\mathcal{U})$ 
and it holds
\begin{align*}
\| \{V_\mathcal{G}f(x_i)\}_{i\in I} | Y^\flat \| = \Big\| \sum_{i\in I} |V_\mathcal{G}f(x_i)|\chi_{U_i} \Big| Y \Big\| \le  \sigma(\mathcal{U})\interleave K^*_\mathcal{U} \interleave \| f | \Co(\mathcal{F},Y) \|.
\end{align*}
\end{corollary}

Let us turn to the synthesis side. Here the following lemma is a key result, which generalizes \cite[Lem.~5.10]{fora05} and whose short direct proof
is new and avoids technical difficulties. In particular, it does not rely on \cite[Lem.~5.4]{fora05}.

\begin{lemma}\label{auxlem:mainsynthesis1}
Let $Y$ be a rich solid QBF-space on $X$ 
and let the analyzing frame $\mathcal{F}=\{\varphi_x\}_{x\in X}$ possess property $F(\nu,Y)$.
Further, let $\mathcal{G}=\{\psi_x \}_{x\in X}$ be a
family in $\mathcal{H}$ and $\mathcal{U}=\{U_i\}_{i\in I}$ an admissible covering such that $K_\mathcal{U}:=K_\mathcal{U}[\mathcal{G},\mathcal{F}]$ defines a bounded
operator on $Y$.
Then for $\{\lambda_i\}_{i\in I}\in Y^\natural(\mathcal{U})$ and for points $x_i\in U_i$
the series $\sum_{i\in I} \lambda_i V_\mathcal{F}\psi_{x_i}(x)$
converges absolutely for a.e.\ $x\in X$ defining a function in $Y$ with
\begin{align*}
\Big\| \sum_{i\in I} \lambda_i V_\mathcal{F}\psi_{x_i} \Big| Y \Big\| \le \interleave K_\mathcal{U} \interleave \| \{\lambda_i\}_{i\in I} | Y^\natural \|.
\end{align*}
If the finite sequences are dense in $Y^\natural(\mathcal{U})$ the series
also converges unconditionally in the quasi-norm of $Y$.
\end{lemma}
\bproof
We have for every $x\in X$ the estimate
\begin{gather*}
\sum_{i\in I} |\lambda_i||V_\mathcal{F}\psi_{x_i}(x)|
\le
\sum_{i\in I} \mu(U_i)^{-1}|\lambda_i| \int_X \chi_{U_i}(y) K_\mathcal{U}(x,y) \,d\mu(y) \\
= \int_X  \sum_{i\in I} \mu(U_i)^{-1}|\lambda_i| \chi_{U_i}(y) K_\mathcal{U}(x,y)   \,d\mu(y)
= K_\mathcal{U} \left( \sum_{i\in I} \mu(U_i)^{-1}|\lambda_i| \chi_{U_i} \right) (x),
\end{gather*}
where summation and integration can be interchanged due to monotone convergence.
Since $\{\lambda_i\}_{i\in I}\in Y^\natural$ the sum $\sum_{i\in I}\mu(U_i)^{-1}|\lambda_i|\chi_{U_i}$ defines pointwise a function in $Y$. By assumption $K_\mathcal{U}$ operates continuously on $Y$ and
hence also $K_\mathcal{U} \left( \sum_{i\in I}\mu(U_i)^{-1}|\lambda_i|\chi_{U_i} \right)\in Y$,
which implies $\left| K_\mathcal{U} \left( \sum_{i\in I} \mu(U_i)^{-1}|\lambda_i| \chi_{U_i} \right)(x) \right|<\infty$ for a.e.\ $x\in X$.
It follows that $\sum_{i\in I} \lambda_iV_\mathcal{F}\psi_{x_i}(x)$
converges absolutely at these points.
As a consequence of the solidity of $Y$ and the pointwise estimate
\begin{align}\label{auxeq:est1}
\Big| \sum_{i\in I}\lambda_iV_\mathcal{F}\psi_{x_i} \Big|
\le \sum_{i\in I} |\lambda_i||V_\mathcal{F}\psi_{x_i}|
\le K_\mathcal{U} \left( \sum_{i\in I} \mu(U_i)^{-1}|\lambda_i|  \chi_{U_i} \right) \in Y,
\end{align}
the measurable functions $\sum_{i\in I} \lambda_iV_\mathcal{F}\psi_{x_i}$ and
$\sum_{i\in I} |\lambda_i||V_\mathcal{F}\psi_{x_i}|$ belong to $Y$ with
\begin{gather}\label{auxeq:est2}
\Big\| \sum_{i\in I} \lambda_iV_\mathcal{F}\psi_{x_i} \Big| Y \Big\|
\le  \Big\| \sum_{i\in I} |\lambda_i|  |V_\mathcal{F}\psi_{x_i}|  \Big| Y \Big\|
\le \interleave K_\mathcal{U}\interleave \| \{\lambda_i \}_{i\in I} |Y^\natural \|.
\end{gather}
It remains to show that $\sum_{i\in I} \lambda_i V_\mathcal{F}\psi_{x_i}$ converges
unconditionally in $Y$ to its pointwise limit, if the finite sequences are dense in $Y^\natural(\mathcal{U})$.
For this we fix an arbitrary bijection $\sigma:\N\rightarrow I$ and obtain as in \eqref{auxeq:est2}
\begin{align}\label{auxeq:tendto0}
\Big\| \sum_{m=n+1}^\infty \lambda_{\sigma(m)} V_\mathcal{F}\psi_{x_{\sigma(m)}} \Big| Y \Big\|
\le \interleave K_\mathcal{U} \interleave \| \Lambda - \Lambda^\sigma_n | Y^\natural \|,
\end{align}
where the sequence $\Lambda^\sigma_n$ is given as in Lemma~\ref{lem:ss_findens}. According to this lemma the right-hand side of \eqref{auxeq:tendto0}
tends to zero for $n\rightarrow\infty$, which finishes the proof.
\eproof

\begin{corollary}\label{cor:familyG}
With the assumptions of the previous lemma $\mathcal{G}=\{ \psi_x\}_{x\in X}\subset \Co(\mathcal{F},Y)$.
\end{corollary}
\bproof
For every $x\in X$ there is an index $i_0\in I$ such that $x\in U_{i_0}$. Set $x_{i_0}:=x$
and choose arbitrary points $x_i\in U_i$ for $i\in I\backslash\{i_0\}$.
Let $\delta^{i_0}$ denote the sequence, which has entry 1 at position $i_0$ and is zero
elsewhere. Since $Y$ is assumed to be rich $\delta^{i_0}\in Y^\natural$ and
by the previous lemma
$V_\mathcal{F}\psi_x=\sum_{i\in I} \delta^{i_0}_i V_\mathcal{F}\psi_{x_i}\in Y$, whence
$\psi_x\in \Co(\mathcal{F},Y)$.
\eproof

The correspondence principle allows to cast Lemma~\ref{auxlem:mainsynthesis1}
in a different form, which corresponds to \cite[Lem.~3.11]{RaUl10}.
However, due to the different deduction
the technical assumption $Y^\natural\hookrightarrow (L_\infty^{1/\nu})^\natural$ is not required any more.

\begin{lemma}\label{auxlem:mainsynthesis2}
With the same assumptions as in Lemma~\ref{auxlem:mainsynthesis1} the
series $\sum_{i\in I} \lambda_i \psi_{x_i}$ converges unconditionally in the weak*-topology of
$(\mathcal{H}_1^\nu)^\urcorner$ to an element $f\in \Co(\mathcal{F},Y)$ with
\begin{align*}
V_\mathcal{F}f=V_\mathcal{F}\Big(\sum_{i\in I} \lambda_i \psi_{x_i}\Big)=\sum_{i\in I} \lambda_i V_\mathcal{F}\psi_{x_i}
\end{align*}
and the estimate \hfill
$\displaystyle{
\| f | \Co(\mathcal{F},Y) \|= \Big\|\sum_{i\in I} \lambda_i \psi_{x_i} \Big| \Co(\mathcal{F},Y) \Big\|
\le \interleave K_\mathcal{U} \interleave \| \{\lambda_i\}_{i\in I} | Y^\natural \|.
}$ \hfill\hfill~\\
Moreover, if the finite sequences are dense in $Y^\natural(\mathcal{U})$ the series
also converges unconditionally in the quasi-norm of $\Co(\mathcal{F},Y)$.
\end{lemma}
\bproof
If the subset $J\subset I$ is finite we have
$
V_\mathcal{F}\Big( \sum_{i\in J} \lambda_i\psi_{x_i} \Big)(x)= \sum_{i\in J} \lambda_iV_\mathcal{F}\psi_{x_i}(x)
$
for all $x\in X$.
Moreover, we have proved in Lemma~\ref{auxlem:mainsynthesis1} that $\sum_{i\in I} \lambda_iV_\mathcal{F}\psi_{x_i}$ converges pointwise
absolutely a.e.\ to a function in $Y$.
In order to apply the correspondence principle, Corollary~\ref{cor:corrpri}, it remains to verify that
the sums $\sum_{i\in J} \lambda_i\psi_{x_i}$ for finite subsets $J\subset I$
are uniformly bounded in $(\mathcal{H}_1^\nu)^\urcorner$.
With the continuous embedding $\Co(Y) \hookrightarrow (\mathcal{H}_1^\nu)^\urcorner$ from Theorem~\ref{thm:co_main}
we can conclude
\begin{gather*}
\Big\| \sum_{i\in J} \lambda_i\psi_{x_i} \Big| (\mathcal{H}_1^\nu)^\urcorner\Big\|\lesssim\Big\|\sum_{i\in J}\lambda_i\psi_{x_i}\Big|\Co(Y)\Big\|
= \Big\| \sum_{i\in J} \lambda_iV_\mathcal{F}\psi_{x_i} \Big| Y \Big\| \le \Big\| \sum_{i\in I} |\lambda_i||V_\mathcal{F}\psi_{x_i}| \Big| Y \Big\|
\end{gather*}
for every finite subset $J\subset I$, where we used that $\psi_{x_i}\in \Co(Y)$ for all $i\in I$ by Corollary~\ref{cor:familyG}.
We have shown in the proof of Lemma~\ref{auxlem:mainsynthesis1} that $\sum_{i\in I} |\lambda_i||V_\mathcal{F}\psi_{x_i}|$ is a function in $Y$.
Hence the sums are uniformly bounded in $(\mathcal{H}_1^\nu)^\urcorner$ and Corollary~\ref{cor:corrpri}
implies the unconditional weak*-convergence
of $\sum_{i\in I} \lambda_i\psi_{x_i}$ to an element
$f\in (\mathcal{H}_1^\nu)^\urcorner$. Moreover, $f\in \Co(Y)$ because Corollary~\ref{cor:corrpri} together with the previous lemma asserts that $V_\mathcal{F}f=\sum_{i\in I} \lambda_i V_\mathcal{F}\psi_{x_i}\in Y$.

It remains to show that $\sum_{i\in I} \lambda_i \psi_{x_i}$ converges
unconditionally in $\Co(\mathcal{F},Y)$, if the finite sequences are dense in $Y^\natural$.
For a subset $\tilde{I}\subset I$ let $\tilde{\Lambda}$ denote the sequence which coincides with $\Lambda$ on $\tilde{I}$ and
is trivial elsewhere.
By solidity $\tilde{\Lambda}\in Y^\natural$ and -- applying what we have proved so far --
the sum $\sum_{i\in \tilde{I}} \lambda_i \psi_{x_i} $ converges in the weak*-topology to an element
of $\Co(Y)$ and
$
V_\mathcal{F}\Big(\sum_{i\in \tilde{I}} \lambda_i \psi_{x_i} \Big) = \sum_{i\in \tilde{I}} \lambda_i V_\mathcal{F}\psi_{x_i}.
$
In view of \eqref{auxeq:tendto0} we conclude
\begin{align*}
\Big\| \sum_{m=n+1}^\infty \lambda_{\sigma(m)} \psi_{x_{\sigma(m)}} \Big| \Co(Y) \Big\|
= \Big\| \sum_{m=n+1}^\infty \lambda_{\sigma(m)} V_\mathcal{F}\psi_{x_{\sigma(m)}} \Big| Y \Big\|
\rightarrow 0 \quad (n\rightarrow\infty),
\end{align*}
for an arbitrary bijection $\sigma:\N\rightarrow I$, which finishes the proof.
\eproof


\subsubsection*{Atomic decompositions}

Our first goal is to obtain atomic decompositions of the coorbit $\Co(Y)$.
Since $\Co(Y)$ is isometrically isomorphic to the function space $R(Y)$
we initially focus on this space and recall from Theorem~\ref{thm:co_main} that for functions $F\in R(Y)$ the reproducing formula holds, i.e.\
\[
F= R(F)= \int_X F(y)R(\cdot,y) \,d\mu(y) ,\qquad F\in R(Y).
\]
This identity can be interpreted as a ``continuous atomic decomposition'' of $F$ with atoms $R(\cdot,y)$
indexed by $y\in X$.
The strategy is to discretize the integral, an approach which originates
from Feichtinger and Gröchenig~\cite{feGr89a} and was also used in subsequent
papers e.g. in \cite{fora05,RaUl10}.
To this end let $\mathcal{U}=\{U_i\}_{i\in I}$ be an admissible
covering of $X$ and let $\Phi=\{\Phi_i\}_{i\in I}$ be a $\mathcal{U}$-PU, i.e.\
a partition of unity subordinate to the covering $\mathcal{U}$ consisting of
measurable functions $\Phi_i$ which satisfy
\begin{enumerate}
\item[(i)] $0\le\Phi_i(x)\le 1$ for all $x\in X$ and all $i\in I$,
\item[(ii)] $\supp \Phi_i\subset U_i$ for all $i\in I$,
\item[(iii)] $\sum_{i\in I} \Phi_i(x) = 1$ for all $x\in X$. 
\end{enumerate}
We note that the construction of such a family $\Phi$ with respect to a locally finite covering is standard,
see e.g.\ \cite[p.~127]{Fol84}.
Using $\Phi$ the integral operator $R$ can be written in the form
\[
R(F)(x)=\sum_{i\in I} \int_X \Phi_i(y)F(y)R(x,y) \,d\mu(y).
\]
A formal discretization yields a discrete integral operator $U_\Phi$, called the \emph{discretization operator},
\begin{align}\label{def:DiscrOp}
U_\Phi F(x) :=\sum_{i\in I} c_i F(x_i)R(x,x_i),
\end{align}
where $c_i:=\int_X \Phi_i(y)\,d\mu(y)$ and the points $\{x_i\}_{i\in I}$ are chosen such that $x_i\in U_i$.
Here we must give meaning to the point evaluations $F(x_i)$ since in general $F\in Y$ only determines an equivalence class of functions where point evaluations are not well-defined.
However, the operator $U_\Phi$ is only applied to elements $F\in R(Y)$ and pointwise evaluation can be understood in the sense
\[
F(x_i)=(RF)(x_i)=\int_X R(x_i,y)F(y) \,d\mu(y)\,.
\]

Intuitively, $U_\Phi F$ approximates $R(F)$ because the discretization resembles a Riemannian sum of the integral.
Hence we can hope to obtain an atomic decomposition from the relation
\[
F= R(F) \approx U_\Phi F = \sum_{i\in I} c_i F(x_i)R(\cdot,x_i).
\]
So far our considerations were just formal. To make the argument precise 
we have to impose conditions on $\mathcal{F}$ so that $U_\Phi$ is a well-defined
operator.
It turns out that here mapping properties of the kernels
$M_\mathcal{U}:=K_\mathcal{U}[\mathcal{F},\mathcal{F}]$ and $M^*_\mathcal{U}:=K^*_\mathcal{U}[\mathcal{F},\mathcal{F}]$
come into play. Recalling the definition~\eqref{eqdef:kerK} of $K_\mathcal{U},\,K^*_\mathcal{U}$ we have for $x,y\in X$
\begin{align}\label{eqdef:kerM}
M_\mathcal{U}(x,y)=\sup_{z\in Q_y} |\langle \varphi_x,\varphi_z \rangle| \quad\text{and}\quad
M^*_\mathcal{U}(x,y)= M_\mathcal{U}(y,x)
\end{align}
with $Q_y=\bigcup\limits_{i~:~y\in U_i} U_i$ for the covering $\mathcal{U}=\{U_i\}_{i\in I}$.

The lemma below provides definition \eqref{def:DiscrOp} with a solid foundation.

\begin{lemma}\label{lem:DiscrOp}
If $M_{\mathcal{U}}$ and
$M^*_{\mathcal{U}}$ given in \eqref{eqdef:kerM} are bounded operators on $Y$
the discretization operator defined in \eqref{def:DiscrOp}
is a well-defined continuous operator $U_\Phi: R(Y)\rightarrow R(Y)$ with operator quasi-norm
$\| U_\Phi | R(Y)\rightarrow R(Y) \| \le \sigma(\mathcal{U}) \interleave M_\mathcal{U}\interleave
\interleave M^*_\mathcal{U} \interleave$.
In general, the convergence of the sum in \eqref{def:DiscrOp} is pointwise absolutely a.e..
If the finite sequences are dense in $Y^\natural$ the convergence is also in the quasi-norm of $Y$.
\end{lemma}
\bproof
For $F\in R(Y)$ Lemma~\ref{lem:reproform3} gives an element $f\in \Co(Y)$ such that $F(x)=Vf(x)$ for all $x\in X$.
Thus, using Corollary~\ref{auxcor:trafosampling} with $\mathcal{G}=\mathcal{F}$, we can conclude $\{F(x_i)\}_{i\in I}\in Y^\flat$ with
$\| \{F(x_i)\}_{i\in I} | Y^\flat \| \le \sigma(\mathcal{U})\interleave M^*_\mathcal{U} \interleave \| F | Y\|$.
Since $\lambda_i\mapsto \mu(U_i)\lambda_i$ is an isometry from $Y^\flat$ to $Y^\natural$ and since $0\le c_i\le\mu(U_i)$ for all $i\in I$
it follows $\{c_iF(x_i)\}_{i\in I}\in Y^\natural(\mathcal{U})$  and $\|  \{c_iF(x_i)\}_{i\in I}   | Y^\natural \| \le \|  \{F(x_i)\}_{i\in I}   | Y^\flat \| $.
Therefore by Lemma~\ref{auxlem:mainsynthesis2} the sum $\sum_{i\in I} c_i F(x_i) \varphi_{x_i}$ converges in the weak*-topology
to an element in $\Co(Y)$ and $U_\Phi F=V\big(\sum_{i\in I} c_i F(x_i)\varphi_{x_i} \big)$. As a consequence $U_\Phi F\in R(Y)$ and again with Lemma~\ref{auxlem:mainsynthesis2}
\begin{align*}
\|U_\Phi F | Y \|&=\| \sum_{i\in I} c_i F(x_i) \varphi_{x_i} | \Co(Y) \|\\
&\le \interleave M_\mathcal{U} \interleave \| \{F(x_i)\}_{i\in I} | Y^\flat \|
\le \sigma(\mathcal{U}) \interleave M_\mathcal{U} \interleave  \interleave M^*_\mathcal{U} \interleave \| F | Y\|.
\end{align*}
\eproof

The operator $U_\Phi$ is self-adjoint in a certain sense.

\begin{lemma}\label{lem:Uselfadj}
Let $\mathcal{U}=\{U_i\}_{i\in I}$ be an admissible covering and assume
that the associated maximal kernels $M_{\mathcal{U}}$ and $M^*_{\mathcal{U}}$ of the analyzing frame
$\mathcal{F}$ belong to $\mathcal{A}_{m_\nu}$. Then $U_\Phi$ is a well-defined operator on
$R(L_\infty^{1/\nu})$ and $R(L_1^\nu)$ and
for every $F\in R(L_\infty^{1/\nu})$ and $G\in R(L_1^\nu)$ it holds
\begin{align}\label{eq:Uselfadj}
\langle U_\Phi F , G \rangle_{L_\infty^{1/\nu}\times L_1^\nu } = \langle  F , U_\Phi G \rangle_{L_\infty^{1/\nu}\times L_1^\nu }.
\end{align}
\end{lemma}
\bproof
For $F\in R(L_\infty^{1/\nu})$ we have $F(x)= \langle F, R(\cdot,x) \rangle_{L_\infty^{1/\nu}\times L_1^\nu}$ and -- by arguments in the proof of Lemma~\ref{lem:DiscrOp} for $Y=L_\infty^{1/\nu}$ --
$\{c_iF(x_i)\}_{i\in I}\in (L_\infty^{1/\nu})^\natural$.
Therefore, $\sum_{i\in I} c_i|F(x_i)||R(\cdot,x_i)|\in L_\infty^{1/\nu}$ by Lemma~\ref{auxlem:mainsynthesis1} and \eqref{auxeq:est1}.
Analogous statements hold for $G\in R(L_1^\nu)$. We conclude
\begin{align*}
\begin{aligned}
\langle U_\Phi F , G \rangle_{L_\infty^{1/\nu}\times L_1^\nu }
&= \sum_{i\in I} c_iF(x_i) \langle R(\cdot,x_i), G \rangle_{L_\infty^{1/\nu}\times L_1^\nu }
= \sum_{i\in I} c_i F(x_i) \overline{G(x_i)} \\
&= \sum_{i\in I} c_i \overline{G(x_i)} \langle F,R(\cdot,x_i) \rangle_{L_\infty^{1/\nu}\times L_1^\nu }
= \langle  F , U_\Phi G \rangle_{L_\infty^{1/\nu}\times L_1^\nu },
\end{aligned}
\end{align*}
where Lebesgue's dominated convergence theorem was used.
\eproof


Our next aim is to find suitable conditions on $\Phi$ and $\mathcal{U}$ such that the discretization operator $U_\Phi$ is invertible.
The possible expansion
\[
F= U_\Phi U_\Phi^{-1} F = \sum_{i\in I} c_i (U_\Phi^{-1}F)(x_i)R(\cdot,x_i)
\]
then yields an atomic decomposition for $F\in R(Y)$.
Intuitively, for the invertibility of $U_\Phi$ the functions $F\in R(Y)$ must be sufficiently ``smooth'', so that a
discrete sampling is possible without loss of information.
Since $R(Y)$ is the isomorphic image of $\Co(Y)$ under the voice transform, we have to ensure that the transforms
$V_\mathcal{F}f$ of elements $f\in \Co(Y)$ are smooth enough.
An appropriate tool for the control of the smoothness are the oscillation kernels, a concept originally due to Feichtinger and Gröchenig.
We use the extended definition from \cite{balhol10}, utilizing a \emph{phase function} $\Gamma:X\times X\rightarrow\phase$ where $\phase=\{z\in\C : |z|=1\}$, namely
\begin{align*}
\osc_{\mathcal{U},\Gamma}(x,y):= \sup_{z\in Q_y}  \left| R_{\cf}(x,y)-\Gamma(y,z)R_{\cf}(x,z) \right| && \text{and}
&& \osc^*_{\mathcal{U},\Gamma}(x,y):= \osc_{\mathcal{U},\Gamma}(y,x)
\end{align*}
with $x,y\in X$ and $Q_y$ as in \eqref{eqdef:kerK}. The choice $\Gamma\equiv1$ yields the kernels used in \cite{fora05,RaUl10}.

We can now formulate a condition on $\mathcal{F}$ which ensures invertibility of $U_\Phi$, but which is weaker than
the assumptions made in \cite{fora05,RaUl10} since we allow a larger class of coverings and weights.
\begin{definition}\label{def:D(d,v,Y)}
We say a tight continuous frame $\mathcal{F}=\{\varphi_x\}_{x\in X}\subset\mathcal{H}$ possesses \emph{property
$D(\delta,\nu,Y)$} for a weight $\nu\ge 1$ and some $\delta>0$ if it has property $F(\nu,Y)$ and if
there exists an admissible covering $\mathcal{U}$ and a phase function $\Gamma:X\times X\rightarrow\phase$ so that
\begin{enumerate}
\item[(i)] $|R_\mathcal{F}|,\, \osc_{\mathcal{U},\Gamma},\, \osc^*_{\mathcal{U},\Gamma}\,\in \mathcal{B}_{Y,m_\nu}$.
\item[(ii)] $\displaystyle{
\|\osc_{\mathcal{U},\Gamma}|\mathcal{B}_{Y,m_\nu}  \| < \delta }$ and $\displaystyle{ \| \osc^*_{\mathcal{U},\Gamma} | \mathcal{B}_{Y,m_\nu} \| < \delta }$.
\end{enumerate}
\end{definition}

\begin{remark}\label{rem:D(d,v,Y)}
A frame $\mathcal{F}$ with property $D(\delta,\nu,Y)$
for a covering $\mathcal{U}$ and a phase function $\Gamma$ automatically possesses properties
$D(\delta,\nu,L_\infty^{1/\nu})$ and $D(\delta,\nu,L_1^\nu)$ for the same covering $\mathcal{U}$ and the same phase function $\Gamma$.
\end{remark}
\bproof
Every $K\in\mathcal{A}_{m_\nu}$ operates continuously on $L_\infty^{1/\nu}$ and $L_1^\nu $
with $\| K |L_\infty^{1/\nu}\rightarrow L_\infty^{1/\nu} \| \le \| K | \mathcal{A}_{m_\nu} \| $ and
$\| K |L_1^\nu\rightarrow L_1^\nu \| \le \| K | \mathcal{A}_{m_\nu} \|  $.
Moreover, for
$Y=L_\infty^{1/\nu}$ or $Y=L_1^\nu$ it holds $R(Y)\hookrightarrow L_\infty^{1/\nu}$ and the algebras
$\mathcal{B}_{Y,m_\nu}$ and $\mathcal{A}_{m_\nu}$ coincide with equal norms.
\eproof

Note that for a measurable kernel function $K:X\times X\rightarrow\C$ the equality $\interleave K \interleave = {\interleave}\,|K|\,{\interleave}$ does not hold in general.
However, we have the following result.
\begin{lemma}\label{lem:keraction}
Let $K,L:X\times X\rightarrow\C$ be two measurable kernels and assume that $|K|$ acts continuously on $Y$.
Then, if $|L(x,y)|\le|K(x,y)|$ for almost all $x,y\in X$, also $L$ acts continuously on $Y$ with the estimate
$\interleave L \interleave \le {\interleave}\,|K|\,{\interleave}$.
In particular, $K$ acts continuously on $Y$ with $\interleave K \interleave \le {\interleave}\,|K|\,{\interleave}$.
\end{lemma}

Let us record an important consequence of the previous lemma.

\begin{corollary}\label{cor:D(d,v,Y)}
If the frame $\mathcal{F}$ has property $D(\delta,\nu,Y)$
the kernels $R_{\cf}$, $|R_{\cf}|$, $\osc_{\mathcal{U},\Gamma}$,
$\osc^*_{\mathcal{U},\Gamma}$, $M_\mathcal{U}$, and $M^*_\mathcal{U}$ are continuous operators on $Y$.
\end{corollary}
\bproof
For all $x,y\in X$ we have $|R_{\cf}(x,y)|\le M_\mathcal{U}(x,y)$ as well as the estimates
\begin{align*}
M_\mathcal{U}(x,y) \le \osc_{\mathcal{U},\Gamma}(x,y) + |R_{\cf}(x,y)| && \text{and} &&
\osc_{\mathcal{U},\Gamma}(x,y)\le M_\mathcal{U}(x,y) + |R_{\cf}(x,y)|.
\end{align*}
The corresponding estimates for the involuted kernels also hold true. Hence Lemma~\ref{lem:keraction} yields the result.
\eproof

The following lemma shows that $U_\Phi F$ approximates $F\in R(Y)$ if the analyzing frame possesses property
$D(\delta,\nu,Y)$ for a suitably small $\delta>0$.
It corresponds to \cite[Thm.~5.13]{fora05}
and the proof is still valid in our setting -- with the triangle inequality replaced by
the corresponding quasi-triangle inequality.

\begin{lemma}\label{auxlem:Uinvert}
Suppose that the analyzing frame $\mathcal{F}$ possesses property $D(\delta,\nu,Y)$ for some $\delta>0$ with
associated covering $\mathcal{U}=\{U_i\}_{i\in I}$ and phase function $\Gamma$.
Then the discretization operator $U_\Phi$ for some $\mathcal{U}$-PU $\Phi$
is a well-defined bounded operator $U_\Phi:R(Y)\rightarrow R(Y)$ and it holds
\begin{align}\label{eq:Uinvert}
\| Id-U_{\Phi}  ~|~ R(Y) \rightarrow R(Y)\| \le \delta (\interleave R \interleave + \interleave M^*_\mathcal{U} \interleave )C_Y .
\end{align}
\end{lemma}
\bproof
For $F\in R(Y)$ there is $f\in\Co(Y)$ with $F=Vf$. By adapting the proof of Lemma~\ref{auxlem:mainanalysis},
it can be shown that $\widetilde{H}:=\sum_{i\in I} \sup_{z\in U_i} |Vf(z)| \Phi_i \in Y$ with $\|\widetilde{H}|Y\|\le \interleave M^*_\mathcal{U} \interleave \| f | \Co(\mathcal{F},Y)\|$.
The intersection number $\sigma(\mathcal{U})$ does not come into play here, since the inequality~\eqref{ineq:aux} can be improved when using $\Phi_i$ instead of $\chi_{U_i}$.
%
A solidity argument yields
$H:=\sum_{i\in I} |F(x_i)|\Phi_i \in Y$ and also $\sum_{i\in I} F(x_i)\overline{\Gamma(\cdot,x_i)}\Phi_i \in Y$ with respective quasi-norms dominated by $\|\widetilde{H}|Y\|$.

Let us introduce the auxiliary operator $S_\Phi:R(Y)\rightarrow R(Y)$, given pointwise for $x\in X$ by
\[
S_\Phi F(x):=R\bigg( \sum_{i\in I} F(x_i)\overline{\Gamma(\cdot,x_i)}\Phi_i \bigg)(x).
\]
Since $F=R(F)$ we can estimate
\[
\|F- S_\Phi F | Y \| = \Big\| R \Big( F - \sum_{i\in I} F(x_i)\overline{\Gamma(\cdot,x_i)}\Phi_i \Big) \Big| Y \Big\|
\le \interleave R \interleave \Big\| F - \sum_{i\in I}F(x_i)\overline{\Gamma(\cdot,x_i)}\Phi_i \Big| Y \Big\|.
\]
We further obtain for every $x\in X$, because $F(x)=R(F)(x)$ even pointwise,
\begin{gather*}
\Big| F(x)- \sum_{i\in I} F(x_i)\overline{\Gamma(x,x_i)}\Phi_i(x) \Big| = \Big| \sum_{i\in I} \big(R(F)(x)-\overline{\Gamma(x,x_i)}R(F)(x_i)\big)\Phi_i(x) \Big| \\
\le \sum_{i\in I} \Phi_i(x)\int_X |R(y,x)-\Gamma(x,x_i)R(y,x_i)| |F(y)|\,d\mu(y)
\le \osc^*_{\mathcal{U},\Gamma}(|F|)(x).
\end{gather*}
We arrive at \hfill
$
\| F - S_\Phi F | Y \| \le \interleave R \interleave \| \osc^*_{\mathcal{U},\Gamma}(|F|) | Y \|
\le \interleave R \interleave \interleave \osc^*_{\mathcal{U},\Gamma} \interleave \| F|Y \|
\le \delta \interleave R \interleave \| F|Y \|.
$
\hfill
\vspace*{1.5ex}
Let us now estimate the difference of $U_\Phi$ and $S_\Phi$. First we see that for $x\in X$
\begin{align*}
S_\Phi F(x)
=\int_X R(x,y) \sum_{i\in I} F(x_i)\overline{\Gamma(y,x_i)}\Phi_i(y) \,d\mu(y)
=\sum_{i\in I} \int_X R(x,y)F(x_i)\overline{\Gamma(y,x_i)}\Phi_i(y) \,d\mu(y).
\end{align*}
Here we used Lebesgue's dominated convergence theorem, which we use again to obtain
\begin{gather*}
\left| U_\Phi F(x)- S_\Phi F(x)  \right| = \Big| \sum_{i\in I} \int_X \Phi_i(y) F(x_i) (R(x,x_i)-\overline{\Gamma(y,x_i)}R(x,y)) \,d\mu(y)  \Big| \\
 \le \sum_{i\in I} \int_X |F(x_i)|\Phi_i(y)\osc_{\mathcal{U},\Gamma}(x,y) \,d\mu(y)
 = \int_X \sum_{i\in I}  |F(x_i)|\Phi_i(y)\osc_{\mathcal{U},\Gamma}(x,y) \,d\mu(y) = \osc_{\mathcal{U},\Gamma}(H)(x),
\end{gather*}
where $H=\sum_{i\in I}  |F(x_i)|\Phi_i$ as above. We conclude
\begin{align*}
\left\| U_\Phi F- S_\Phi F | Y  \right\| \le \left\| \osc_{\mathcal{U},\Gamma}(H) | Y  \right\|
\le \interleave \osc_{\mathcal{U},\Gamma} \interleave \|H |Y\| \le \delta \interleave M^*_\mathcal{U} \interleave \| F | Y \|.
\end{align*}
Hence, altogether we have proved
\[
\| F -U_\Phi F | Y \|\le
C_Y ( \| F - S_\Phi F|Y\| +   \|S_\Phi F-U_\Phi F | Y \|) \le \delta C_Y \| F | Y \| ( \interleave M^*_\mathcal{U} \interleave + \interleave R \interleave ) .
\]
\eproof

If the righthand side of \eqref{eq:Uinvert} is less than one, $U_\Phi:R(Y)\rightarrow R(Y)$ is boundedly invertible with
the Neumann expansion $U_\Phi^{-1}=\sum_{n=0}^\infty (Id-U_\Phi)^n$, which is still valid in the quasi-Banach setting.

Finally, we are able to prove a cornerstone of the discretization theory, which generalizes \cite[Thm.~5.7]{fora05} and \cite[Thm.~3.11]{RaUl10}.
Note that the characterization via the sequence
spaces is a new result even in the Banach case and that we can drop many technical restrictions.

\begin{Theorem}\label{thm:atomicdec}
Let $Y$ be a rich solid QBF-space with quasi-norm constant $C_Y$ 
and suppose that the analyzing frame $\mathcal{F}=\{\varphi_x\}_{x\in X}$
possesses property $D(\delta,\nu,Y)$ for the covering $\mathcal{U}=\{U_i\}_{i\in I}$
and a small enough $\delta>0$ such that
\begin{align}\label{eq:dcond}
\delta\big( (1+C_Y)\big\| |R_{\mathcal{F}}| \big| \mathcal{B}_{Y,m_\nu} \big\| + \delta C_Y \big)C_Y \le 1.
\end{align}
Choosing arbitrary points $x_i\in U_i$, the sampled frame
$\mathcal{F}_d:=\{\varphi_i\}_{i\in I}:=\{\varphi_{x_i}\}_{i\in I}$
then possesses a ``dual family''
$\widehat{\mathcal{F}_d}=\{\psi_i\}_{i\in I}\subset \mathcal{H}_1^\nu\cap \Co(Y)$
such that the following holds true:
\begin{enumerate}
\item[(i)] (Analysis) An element $f\in(\mathcal{H}_1^\nu)^\urcorner$ belongs to $\Co(Y)$ if and only if
$\{ \langle f,\varphi_i \rangle \}_{i\in I}\in Y^\flat(\mathcal{U})$
\textup{(}or $\{ \langle f,\psi_i \rangle \}_{i\in I}\in Y^\natural(\mathcal{U})$\textup{)}
and we have the quasi-norm equivalences
\begin{align*}
\| f | \Co(Y) \| \asymp \| \{ \langle f,\varphi_{i} \rangle \}_{i\in I} | Y^\flat(\mathcal{U}) \|
\quad\text{ and }\quad
\| f | \Co(Y) \| \asymp \| \{ \langle f,\psi_i \rangle \}_{i\in I} | Y^\natural(\mathcal{U}) \|.
\end{align*}
\item[(ii)] (Synthesis) For every sequence $\{\lambda_i\}_{i\in I}\in Y^\natural(\mathcal{U})$
it holds $f=\sum_{i\in I} \lambda_i\varphi_i \in \Co(Y)$ with $\| f | \Co(Y) \| \lesssim \| \{\lambda_i\}_{i\in I} |Y^\natural(\mathcal{U}) \|$.
In general, the convergence of the sum is in the weak*-topology induced by $(\mathcal{H}_1^\nu)^\urcorner$. It is unconditional in the \mbox{quasi-}norm of $\Co(Y)$,
if the finite sequences are dense in $Y^\natural$.
Similarly, $f=\sum_{i\in I} \lambda_i\psi_i \in \Co(Y)$ with $\| f | \Co(Y) \| \lesssim \| \{\lambda_i\}_{i\in I} |Y^\flat(\mathcal{U}) \|$ in case $\{\lambda_i\}_{i\in I}\in Y^\flat(\mathcal{U})$.
\item[(iii)] (Reconstruction) For all $f\in \Co(Y)$ we have
$
f=\sum_{i\in I} \langle f,\psi_i \rangle\varphi_i$ 
and $f=\sum_{i\in I} \langle f,\varphi_i \rangle \psi_i$. 
\end{enumerate}
\end{Theorem}
\bproof
According to Remark~\ref{rem:D(d,v,Y)} the frame $\mathcal{F}$ has properties
$D(\delta,\nu,L_1^\nu)$ and $D(\delta,\nu,L_\infty^{1/\nu})$ with respect to the covering $\mathcal{U}$,
and by Lemma~\ref{lem:co_examples} it holds
$(\mathcal{H}_1^\nu)^\urcorner \asymp \Co(L_\infty^{1/\nu})$ and $\mathcal{H}_1^\nu= \Co(L_1^\nu)$.
In view of Theorem~\ref{thm:co_main} the voice transform
$V:(\mathcal{H}_1^\nu)^\urcorner\rightarrow R(L_\infty^{1/\nu})$ is thus
a boundedly invertible operator with isometric restrictions $V:\Co(Y)\rightarrow R(Y)$ and $V:\mathcal{H}_1^\nu\rightarrow R(L_1^\nu)$.

Let us fix a $\mathcal{U}$-PU $\Phi=\{\Phi_i\}_{i\in I}$ and put $c_i:=\int_X \Phi_i(y)\,d\mu(y)$.
According to Lemma~\ref{lem:Uselfadj} the corresponding discretization operator $U_\Phi$ is well-defined and bounded on $R(L_\infty^{1/\nu})$.
Condition~\eqref{eq:dcond} on $\delta$ further implies that $U_\Phi:R(L_\infty^{1/\nu})\rightarrow R(L_\infty^{1/\nu})$ is boundedly invertible
as a consequence of Lemma~\ref{auxlem:Uinvert}.
Indeed, using the estimates $\interleave M^*_{\mathcal{U}} \interleave \le C_Y (\interleave |R_{\cf}| \interleave + \interleave \osc^*_{\mathcal{U},\Gamma} \interleave)$ and $\interleave R_{\cf} \interleave \le \interleave |R_{\cf}| \interleave$ together with the assumption $\interleave \osc^*_{\mathcal{U},\Gamma} \interleave<\delta$ we can deduce
\begin{align*}
\delta(\interleave R_{\cf} \interleave + \interleave M^*_{\mathcal{U}} \interleave ) C_Y
&\le \delta( (1+C_Y) \interleave |R_{\cf}| \interleave + C_Y\interleave \osc^*_{\mathcal{U},\Gamma} \interleave) C_Y\\
&< \delta( (1+C_Y) \| |R_{\cf}| | \mathcal{B}_{Y,m_\nu} \| + C_Y \delta )C_Y \le 1.
\end{align*}
Analogously, it follows that
$U_\Phi:R(L_1^\nu)\rightarrow R(L_1^\nu)$ and $U_\Phi:R(Y)\rightarrow R(Y)$ are
boundedly invertible.

For the proof it is useful to note that the operator
$
T:=V^{-1}U_\Phi^{-1}V ~:~ (\mathcal{H}_1^\nu)^\urcorner \rightarrow (\mathcal{H}_1^\nu)^\urcorner
$
is a boundedly invertible isomorphism, whose restrictions $T:\mathcal{H}_1^\nu\rightarrow \mathcal{H}_1^\nu$ and $T:\Co(Y)\rightarrow \Co(Y)$ are also boundedly invertible.
Moreover, $T$ is ``self-adjoint''. For this observe that relation \eqref{eq:Uselfadj} also
holds for the inverse $U_\Phi^{-1}=\sum_{n=0}^\infty (Id-U_\Phi)^n$. 
Consequently, for $f\in(\mathcal{H}_1^\nu)^\urcorner$ and $\zeta\in\mathcal{H}_1^\nu$
\begin{align*}
\langle f,T\zeta \rangle = \langle f,V^{-1}U_\Phi^{-1}V\zeta \rangle
=\langle Vf , U_\Phi^{-1}V\zeta \rangle_{L_\infty^{1/\nu}\times L_1^\nu}
=\langle U_\Phi^{-1}Vf , V\zeta \rangle_{L_\infty^{1/\nu}\times L_1^\nu}
=\langle Tf , \zeta \rangle.
\end{align*}
It follows further that $T$ is sequentially continuous with respect to the weak*-topology of $(\mathcal{H}_1^\nu)^\urcorner$.
To see this let $f_n \rightarrow f$ in the weak*-topology. Then
$
\langle Tf_n , \zeta \rangle = \langle f_n , T\zeta \rangle \rightarrow
\langle f , T\zeta \rangle = \langle Tf , \zeta \rangle 
$
for every $\zeta\in\mathcal{H}_1^\nu$.
By Lemma~\ref{cor:familyG}, Corollary~\ref{cor:D(d,v,Y)} and Lemma~\ref{lem:co_examples} the atoms $\varphi_{x_i}$ lie in $\mathcal{H}_1^\nu\cap \Co(Y)$.
Since $T$ respects these subspaces we can define
\[
\psi_i:= c_iT\varphi_i \:\in \mathcal{H}_1^\nu\cap \Co(Y)
\]
and claim that $\widehat{\mathcal{F}_d}=\{\psi_i\}_{i\in I}$ is the desired ``dual'' of $\mathcal{F}_d=\{\varphi_i\}_{i\in I}$.

After these preliminary considerations we now turn to the proof of the assertions.

\noindent
\textit{Step 1.}\,
If $f\in \Co(Y)$ then $\{ \langle f,\varphi_i \rangle \}_{i\in I}= \{ Vf(x_i) \}_{i\in I}\in Y^\flat$
and $\| \{ \langle f,\varphi_i \rangle \}_{i\in I} | Y^\flat \| \lesssim  \| f | \Co(Y)\|$ by Corollary~\ref{auxcor:trafosampling}.
Furthermore, it holds $Tf\in \Co(Y)$ and Corollary~\ref{auxcor:trafosampling} yields 
$
\{ \langle f,\psi_i \rangle \}_{i\in I}= \{ c_i \langle Tf, \varphi_i \rangle \}_{i\in I}\in Y^\natural
$
with the estimate
$\|\{ \langle f,\psi_i \rangle \}_{i\in I} | Y^\natural \| \le
\| \{\langle Tf, \varphi_i \rangle \}_{i\in I} | Y^\flat \| \lesssim  \| Tf | \Co(Y)\|
\lesssim \|f | \Co(Y)\|$.

\noindent
\textit{Step 2.}\,
If $\{ \lambda_i \}_{i\in I}\in Y^\natural$ then by Lemma~\ref{auxlem:mainsynthesis2} the sum $\sum_{i\in I} \lambda_i\varphi_i$ converges
in the weak*-topology to an element in $\Co(Y)$ with estimate
$\| \sum_{i\in I} \lambda_i\varphi_i | \Co(Y) \|\lesssim \| \{\lambda_i\}_{i\in I} | Y^\natural \|$.
If the finite sequences are dense in $Y^\flat$ (or equivalently $Y^\natural$)
the convergence is even in the quasi-norm of $\Co(Y)$.

A similar statement holds for the dual family $\{\psi_i\}_{i\in I}$.
Indeed, for  $\{ \lambda_i \}_{i\in I}\in Y^\flat$ we have $\{ c_i\lambda_i \}_{i\in I}\in Y^\natural$
and hence $\sum_{i\in I} c_i\lambda_i\varphi_i$ converges
in the weak*-topology to an element in $\Co(Y)$. Since $T$ is sequentially continuous it follows that
\begin{align*}
\sum_{i\in I} \lambda_i\psi_i = \sum_{i\in I} c_i\lambda_iT\varphi_i = T\left( \sum_{i\in I} c_i\lambda_i\varphi_i \right)
\in \Co(Y)
\end{align*}
with weak*-convergence in the sums. The operator $T$ is also continuous on $\Co(Y)$, proving the quasi-norm convergence if the finite sequences are dense.  Moreover, we have the estimate
\begin{align*}
\Big\| \sum_{i\in I} \lambda_i\psi_i \Big| \Co(Y) \Big\| \lesssim  \Big\|  \sum_{i\in I} c_i\lambda_i\varphi_i \Big| \Co(Y) \Big\|
\lesssim  \| \{c_i\lambda_i\}_{i\in I} | Y^\natural \| \le \| \{\lambda_i\}_{i\in I} | Y^\flat \|.
\end{align*}

\noindent
\textit{Step 3.}\,
In this step we prove the expansions in (iii).
For $f\in(\mathcal{H}_1^\nu)^\urcorner$ we have the identity
\begin{align*}
Vf
=U_\Phi \left( U_\Phi^{-1}Vf \right)
= \sum_{i\in I} c_i \left(U_\Phi^{-1}  Vf\right)(x_i)  R(\cdot,x_i)
=\sum_{i\in I} \langle f,\psi_i \rangle V\varphi_i
\end{align*}
with pointwise absolute convergence a.e.\ in the sums. Since
$(\mathcal{H}_1^\nu)^\urcorner\asymp \Co(L_\infty^{1/\nu})$ the coefficients
$\{ \langle f,\psi_i \rangle \}_{i\in I}$ belong to $(L_\infty^{1/\nu})^\natural$ according to Step 1.
Hence, by Lemma~\ref{auxlem:mainsynthesis2} it holds $Vf=V(\sum_{i\in I} \langle f,\psi_i \rangle \varphi_i)$ with weak*-convergence of the sum.
The injectivity of $V$ finally yields
\begin{align} \label{proofeq:Atomic3}
f=\sum_{i\in I} \langle f,\psi_i \rangle \varphi_i.
\end{align}
Using the sequential continuity of $T$ with respect to the weak*-topology we can further deduce
\begin{align}\label{proofeq:Atomic4}
f=TT^{-1}f=\sum_{i\in I} \langle T^{-1}f,\psi_i \rangle T\varphi_i = \sum_{i\in I} \langle T^{-1}f,c_iT\varphi_i \rangle T\varphi_i
= \sum_{i\in I} \langle f, \varphi_i \rangle \psi_i.
\end{align}
In particular, these expansions are valid for $f\in \Co(Y)$ with coefficients $\{ \langle f,\psi_i \rangle \}_{i\in I}\in Y^\natural$
and $\{ \langle f,\varphi_i \rangle \}_{i\in I}\in Y^\flat$ by Step 1.

\noindent
\textit{Step 4.}\,
If $f\in(\mathcal{H}_1^\nu)^\urcorner$ and either $\{ \langle f,\varphi_i \rangle \}_{i\in I}\in Y^\flat$ or
$\{ \langle f,\psi_i \rangle \}_{i\in I} \in Y^\natural $ we can conclude from the expansions \eqref{proofeq:Atomic3}
and \eqref{proofeq:Atomic4}
together with Step 2 that $f\in \Co(Y)$. Moreover,
$ \| \{\langle f,\psi_i \rangle \}_{i\in I}  | Y^\natural \|$ and
$ \| \{\langle f,\varphi_i \rangle \}_{i\in I}  | Y^\flat \|$ are equivalent quasi-norms on $\Co(Y)$ because using Steps 1 and 2
\begin{flalign*}
&&\| f | \Co(Y) \| = \Big\| \sum_{i\in I} \langle f,\psi_i \rangle \varphi_i \Big| \Co(Y) \Big\|
\lesssim \| \{\langle f,\psi_i \rangle \}_{i\in I}  | Y^\natural \| \lesssim \| f | \Co(Y) \| &&\\
\text{and} && \| f | \Co(Y) \| = \Big\| \sum_{i\in I} \langle f,\varphi_i \rangle \psi_i \Big| \Co(Y) \Big\|
\lesssim \| \{\langle f,\varphi_i \rangle \}_{i\in I}  | Y^\flat \| \lesssim \|
f | \Co(Y) \|. &&
\end{flalign*}
\eproof

\begin{remark}
Properties (i)-(iii) in particular show that the discrete families $\mathcal{F}_d$ and $\widehat{\mathcal{F}_d}$ both constitute atomic decompositions for $\Co(Y)$, as well as quasi-Banach frames,
compare e.g.\ \cite{RaUl10,ra05-3}.
\end{remark}

\subsubsection*{Frame expansion}

Now we come to another main discretization result, which allows
to discretize the coorbit space $\Co(Y)=\Co(\cf,Y)$ by samples of a frame $\cg = \{\psi_x\}_{x\in X}$
different from the analyzing frame $\cf$.
It is a generalization of \cite[Thm.~3.14]{RaUl10},
whose original proof carries over to
the quasi-Banach setting based on
Corollary~\ref{auxcor:trafosampling} and Lemma~\ref{auxlem:mainsynthesis2}. In contrast to Theorem~\ref{thm:atomicdec}, here we require the
additional property of the covering $\mathcal{U}=\{U_i\}_{i\in I}$ that for some constant $D>0$
\begin{align}\label{eq:covbounbel}
\mu(U_i)\ge D \quad\text{for all }i\in I.
\end{align}

\begin{Theorem}\label{thm:frameexp}
Let $Y$ be a rich solid QBF-space on $X$ 
and assume that the analyzing frame
$\mathcal{F}=\{\varphi_x\}_{x\in X}$ has property $F(\nu,Y)$.
For $r\in\{1,\ldots,n\}$ let $\mathcal{G}_r=\{\psi_x^r\}_{x\in X}$ and $\tilde{\mathcal{G}}_r=\{\tilde{\psi}_x^r\}_{x\in X}$ be families in $\mathcal{H}$, and suppose
that for some admissible covering $\mathcal{U}=\{U_i\}_{i\in I}$ with the additional property \eqref{eq:covbounbel} the kernels
$K_r:=K_\mathcal{U}[\mathcal{G}_r,\mathcal{F}]$ and
$\tilde{K}^*_r:=K_\mathcal{U}^*[\tilde{\mathcal{G}}_r,\mathcal{F}]$ belong to $\mathcal{B}_{Y,{m_\nu}}$.
Then, if every $f\in\mathcal{H}$ has an expansion
\begin{align}\label{eq:frexpansion}
f=\sum_{r=1}^n \sum_{i\in I} \langle f , \tilde{\psi}^r_{x_i} \rangle \psi^r_{x_i}
\end{align}
with fixed points $x_i\in U_i$,
this expansion extends to all $f\in \Co(Y)=\Co(\mathcal{F},Y)$.
Furthermore, $f\in (H^1_\nu)^\urcorner$ belongs to $\Co(Y)$ if and only if $\{\langle f,\tilde{\psi}^r_{x_i} \rangle\}_{i\in I}\in Y^\natural(\mathcal{U})$ for each $r\in\{1,\ldots,n\}$, and in this case we have
$
\| f | \Co(Y) \| \asymp \sum_{r=1}^n \left\| \{ \langle f,\tilde{\psi}_{x_i}^r \rangle \}_{i\in I} | Y^\natural(\mathcal{U}) \right\| .
$
The convergence in \eqref{eq:frexpansion} is in the quasi-norm of $\Co(Y)$ if the finite sequences are dense in $Y^\natural(\mathcal{U})$.
In general, we have weak*-convergence induced by $(\mathcal{H}_1^\nu)^\urcorner$.
\end{Theorem}

Observe that the technical assumption $Y^\natural\hookrightarrow(L^{1/v}_\infty)^\natural$ made in \cite[Thm.~3.14]{RaUl10} is not necessary.
In view of Lemma~\ref{auxlem:crossG} it is further not necessary to require $\mathcal{G}_r, \tilde{\mathcal{G}}_r\subset \mathcal{H}_1^\nu$. In fact,
$K_r,\tilde{K}^*_r \in\mathcal{A}_{m_\nu}$ is a stronger condition than $G[\mathcal{G}_r,\mathcal{F}], G^*[\tilde{\mathcal{G}}_r,\mathcal{F}]\in\mathcal{A}_{m_\nu}$
and implies $\mathcal{G}_r, \tilde{\mathcal{G}}_r\subset \mathcal{H}_1^\nu$.

\section{Variable exponent spaces}
\label{sec:varint}

In the remainder we give a demonstration of the theory.
As an example we will show that variable exponent spaces, which have caught some attention recently, fall into
the framework of coorbit theory and can be handled conveniently within the theory.

\subsection{Spaces of variable integrability}
The spaces of variable integrability $\Lpp$ were first introduced by Orlicz~\cite{Orlicz31} in 1931 as a generalization of the Lebesgue spaces $L_p(\R)$.
Before defining them let us introduce some standard notation from \cite{KovacikRakosnik91}. For a measurable function $p:\R\to(0,\infty]$ and a set $\Omega\subset\R$
we define the quantities $p^-_\Omega=\essinf{x\in\Omega}p(x)$ and $p^+_\Omega=\esssup{x\in\Omega}\ p(x)$. Furthermore, we abbreviate $p^-=p^-_\R$ and $p^+=p^+_\R$ and
say that $\p$ belongs to the class of admissible exponents $\P$ if $p^->0$.
Having an admissible exponent $p\in\P$ we define the set $\R_{\!\!\!\infty}=\{x\in\R:p(x)=\infty\}$ and for every measurable function $f:\R\to\C$
the modular
\begin{align*}
	\varrho_\p(f)=\int_{\R\setminus\R_{\!\!\!\infty}}|f(x)|^{p(x)}dx+\esssup{x\in\R_{\!\!\!\infty}}|f(x)|\punkt
\end{align*}
\begin{definition}\label{Lppunkt} The space $\Lpp$ is the collection of all
functions $f$ such that there exists a $\lambda>0$ with $\varrho_\p(\lambda
f)<\infty$. It
is equipped with the Luxemburg quasi-norm
\begin{align*}
	\norm{f}{\Lpp}=\inf\left\{\lambda>0:\varrho_\p\left(\frac{f}{\lambda}\right)<1\right\}\punkt
\end{align*}
\end{definition}
The spaces $\Lpp$ share many properties with the constant exponent spaces $L_p(\R)$. Let us mention a few; the proofs can be found in \cite{KovacikRakosnik91} and in \cite{DieningHastoBuch2011}:
\begin{itemize}
	\item If $p(x)=p$ then $\Lpp=L_p(\R)$,
	\item if $|f(x)|\geq|g(x)|$ for a.e.\ $x\in\R$ then $\varrho_\p(f)\geq\varrho_\p(g)$ and $\norm{f}{\Lpp}\geq\norm{g}{\Lpp}$,
	\item $\varrho_\p(f)=0$ if and only if $f=0$,
	\item for $p(\cdot)\geq1$ H\"older's inequality holds \cite[Theorem 2.1]{KovacikRakosnik91}
		\begin{align*}
			\int_\R|f(x)g(x)|dx\leq 4\norm{f}{\Lpp}\norm{g}{L_{p'(\cdot)}(\R)}\komma
		\end{align*}
		where $1/\p+1/p'(\cdot)=1$ pointwise.
\end{itemize}
There are also some properties of the usual constant exponent spaces which the $\Lpp$ spaces do not share. For example in general the $\Lpp$ spaces are not translation invariant, i.e.\
$f\in\Lpp$ does not automatically imply that $f(\cdot+h)$ belongs to $\Lpp$ for $h\in\R$. As a consequence also Young's convolution inequality is not valid (see again \cite{KovacikRakosnik91} for details).

The breakthrough for $\Lpp$ spaces was made by Diening in \cite{Diening2004} when he showed that the Hardy-Littlewood maximal operator $\HLM$ is bounded on $\Lpp$ under certain regularity conditions on $\p$.
His result has been generalized in many cases (see \cite{DieningHarjulehto2009},\cite{Nekvinda2004} and \cite{CruzUribe2003}) and it turned out that logarithmic H\"older continuity classes are well adapted to the boundedness of the maximal operator.
\begin{definition}
	 Let $g\in C(\R)$. We say that $g$ is \emph{locally $\log$-H\"older continuous}, abbreviated $g\in C^{\log}_{\rm loc}(\R)$, if there exists $c_{\log}>0$ such that
    \begin{align*}
        |g(x)-g(y)|\leq\frac{c_{\log}}{\log(\e+{1}/{|x-y|})} \qquad\text{for all }x,y\in\R.
    \end{align*}
    We say that $g$ is \emph{globally $\log$-H\"older continuous}, abbreviated $g\in C^{\log}(\R)$, if $g$ is locally $\log$-H\"older continuous and there exists $g_\infty\in\re$ such that
    \begin{align*}
        |g(x)-g_\infty|\leq\frac{c_{\log}}{\log(\e+|x|)} \qquad\text{for all }x\in\R.
    \end{align*}
\end{definition}
With the help of the above logarithmic H\"older continuity the following result holds.

\begin{lemma}[{\cite[Thm.~3.6]{DieningHarjulehto2009}}]\label{lem:HLM}
	Let $p\in\P$ with $1<p^-\leq p^+\leq\infty$. If $\frac{1}{p}\in C^{\log}(\R)$, then $\HLM$ is bounded on $\Lpp$ i.e., there exists $c>0$ such that for all $f\in\Lpp$
\begin{align*}
	\norm{\HLM f}{\Lpp}\leq c\norm{f}{\Lpp}\punkt
\end{align*}	
\end{lemma}

Since logarithmic H\"older continuous exponents play an essential role we introduce the class $\Plog$ of admissible exponents $\p$ with $1/p\in C^{\log}(\R)$ and $0<p^-\leq p^+\leq\infty$.
As a consequence of Lemma~\ref{lem:HLM}, for exponents $p\in\Plog$ the maximal operator $\HLM$ is bounded on $L_{\frac\p t}(\R)$ for every $0<t<p^-$.

\subsection{2-microlocal function spaces with variable integrability}

We proceed with spaces of Besov-Triebel-Lizorkin type featuring variable integrability and smoothness.
Spaces of the form $F^{s(\cdot)}_{\p,\q}(\R)$ and $B^{s(\cdot)}_{\p,\qconst}(\R)$ have been studied in \cite{DieningHastoRoudenko2009,AlmeidaHasto2010},
where $s:\R\to\mathbb{R}$ with $s\in L_\infty(\R)\cap C^{\log}_{\rm loc}(\R)$.
A further generalization was pursued in \cite{Ke09,Ke11} replacing the smoothness parameter $s(\cdot)$ by a more general weight function $w$.
We make some reasonable restrictions on $w$ and use the class $\mathcal{W}^{\alpha_3}_{\alpha_1,\alpha_2}$ of admissible weights introduced in \cite{Ke09}.

\begin{definition}
For real numbers $\alpha_3\geq0$ and $\alpha_1\leq\alpha_2$ a weight function $w:X \to (0,\infty)$ on the index set $X=\R\times[(0,1)\cup\{\infty\}]$ belongs
to the class $\mathcal{W}^{\alpha_3}_{\alpha_1,\alpha_2}$
if and only if for $\x=(x,t) \in X$,
\begin{description}
 \item(W1) 
 $
    \left\{\begin{array}{lcl}
              \Big(\frac{s}{t}\Big)^{\alpha_1}w(x,s) \leq w(x,t) \leq \Big(\frac{s}{t}\Big)^{\alpha_2}w(x,s)&,& s \geq t\\\\
              t^{-\alpha_1}w(x,\infty) \leq w(x,t) \leq t^{-\alpha_2}w(x,\infty)&,& s = \infty\,,
          \end{array}\right.
 $
 \item(W2)
 $
    w(x,t) \leq w(y,t)\left\{\begin{array}{rcl}
                                (1 + |x-y|/t)^{\alpha_3}&,& t\in (0,1)\\
                                (1+|x-y|)^{\alpha_3}&,& t=\infty
                             \end{array}\right.\, \quad\mbox{ for all }y \in \R.
 $
\end{description}
\end{definition}
\begin{example}\label{Exampel2ml}  The main examples are weights of the form
$$
    w_{s,s'}(x,t) = \left\{\begin{array}{rcl}
                                t^{-s}\Big(1+\frac{|x-x_0|}{t}\Big)^{s'} &,& t\in (0,1)\\
                                (1+|x-x_0|)^{s'}&,& t=\infty
                    \end{array}\right.\,.
$$
where $s,s' \in \re$. These weights are continuous versions of 2-microlocal weights, used to define 2-microlocal function spaces of Besov-Lizorkin-Triebel type, see \cite{Ke09,Ke10,Ke11}.\\
By choosing $s'=0$ we get back to usual Besov-Lizorkin-Triebel spaces with smoothness $s\in\re$.
\end{example}
The special weights from this example are usually called 2-microlocal weights. Furthermore, function spaces which are defined with admissible weights $w\in \mathcal{W}^{\alpha_3}_{\alpha_1,\alpha_2}$ are usually called 2-microlocal spaces. This term was coined by Bony \cite{Bony} and Jaffard \cite{Jaffard}, who also introduced the concept of 2-microlocal analysis to study local regularity of functions.

\begin{remark}
By the conditions on admissible weights $w\in \mathcal{W}^{\alpha_3}_{\alpha_1,\alpha_2}$ we obtain the following estimates which will be useful later on:
\begin{enumerate}
	\item For $s\leq t$ we get from (W1)
	\begin{align}\label{eq_W1tilde}
		\left(\frac{s}{t}\right)^{\alpha_2}w(x,s)\leq w(x,t)\leq\left(\frac{s}{t}\right)^{\alpha_1}w(x,s).
	\end{align}
	\item For $0<c<s/t$ we have from (W1) and \eqref{eq_W1tilde}
	\begin{align}\label{eq_st1}
	 \frac{w(x,t)}{w(x,s)}\leq\max\{1,c^{\alpha_1-\alpha_2}\}\left(\frac{s}{t}\right)^{\alpha_2}.
	\end{align}
	\item For $0<c<t/s$ we obtain similarly from (W1) and \eqref{eq_W1tilde}
	\begin{align}\label{eq_st2}
	 \frac{w(x,t)}{w(x,s)}\leq\max\{1,c^{\alpha_1-\alpha_2}\}\left(\frac{s}{t}\right)^{\alpha_1}.
	\end{align}
	\item Consequently, we have for $0<c_1<s/t<c_2$ from \eqref{eq_st1} and \eqref{eq_st2}
	\begin{align*}
	 w(x,t)\asymp w(x,s)\quad\text{for all $x\in\R$.}
	\end{align*}
	\item Using (W2) and the inequalities \eqref{eq_st1} and \eqref{eq_st2} we can relate $w(x,t)$ to $w(0,1/2)$ by
	\begin{align*}
	  w(0,1/2)t^{-\alpha_1}(1+|x|)^{-\alpha_3} 
      \lesssim w(x,t) 
      \lesssim w(0,1/2)t^{-\alpha_2}(1+|x|)^{\alpha_3}.
	\end{align*}
\end{enumerate}
\end{remark}

A weight $w \in \mathcal{W}^{\alpha_3}_{\alpha_1,\alpha_2}$ gives rise to a
semi-discrete counterpart $(w_j)_{j \in \N_0}$, corresponding to
an admissible weight sequence in the sense of \cite{Ke09,Ke10,Ke11},
given by
\begin{equation}\label{eqdef:wj}
  \begin{split}
    w_j(x) = \left\{\begin{array}{rcl}
                       w(x,2^{-j})&,&j\in \N\,,\\
                       w(x,\infty)&,&j=0\,.
                    \end{array}\right.
  \end{split}
\end{equation}

In \cite[Lemma~2.6]{Ke11} it was shown that it is equivalent to consider a smoothness function $s\in L_\infty(\R)\cap C^{\log}_{\rm loc}(\R)$ or an admissible weight sequence stemming from
$w\in \mathcal{W}^{\alpha_3}_{\alpha_1,\alpha_2}$ if they are connected by $w_j(x)=2^{js(x)}$, see \eqref{eqdef:wj}. But there exist weight sequences (Example \ref{Exampel2ml} with $s'\neq0$) where it is not possible to find a smoothness function $s:\R\to\mathbb{R}$ such that the above relation holds.\\
Recently in \cite{Tu14} the concept of admissible weight sequences was extended to include more general weights. We will not follow this generalization of admissible weights, but we remark that by this definition we can have local Muckenhoupt weights as components in the sequence.\\
The spaces $B^w_{\p,\qconst}(\R)$ and $F^w_{\p,\q}(\R)$ are defined Fourier analytical as subspaces of the tempered distributions $\mathcal{S}'(\R)$. As usual the Schwartz space $\mathcal{S}(\re^d)$ denotes
the locally convex space of rapidly decreasing infinitely differentiable functions on $\re^d$. Its topology is generated by the seminorms
\begin{align*}
\|\varphi\|_{k,l}=\sup_{x\in\re^d}(1+|x|)^k\sum_{|\beta|\leq l}|D^\beta\varphi(x)|
\end{align*}
for every $k,l\in\N_0$. Its topological dual, the space of tempered distributions on $\re^d$, is denoted by $\mathcal{S}'(\re^d)$. The Fourier transform and its inverse are defined on both $\mathcal{S}(\R)$ and $\mathcal{S}'(\R$) (see Appendix A.1) and we denote them by $\hat{f}$ and $f^\vee$.
Finally, we introduce the subspace $\mathcal{S}_0(\R)$ of $\mathcal{S}(\R)$ by
\[
\mathcal{S}_0(\R):=\left\{ f\in\mathcal{S}(\R) ~:~ D^{\bar{\alpha}}\widehat{f}(0)=0 \text{ for every multi-index } \bar{\alpha}\in\N_0^d \right\}.
\]
The definition of $B^w_{\p,\qconst}(\R)$ and $F^w_{\p,\q}(\R)$ relies on a dyadic decomposition of unity, see also \cite[2.3.1]{Tr83}.

\begin{definition}
Let $\Pi(\R)$ be the collection of all systems $\{\varphi_j\}_{j\in
\n} \subset \mathcal{S}(\R)$ such that
\begin{description}
    \item(i) there is a function $\varphi\in \mathcal{S}(\R)$ with $\varphi_j(\xi) = \varphi(2^{-j}\xi)\,,\, j\in
    \N$\,,
    \item(ii) $\supp \varphi_0 \subset \{\xi\in \R~:~|\xi|\leq
    2\}\,,\quad
    \supp \varphi \subset \{\xi\in \R~:~1/2 \leq |\xi| \leq 2\}$\,, 
    \item(iii) $\sum\limits_{j=0}^{\infty} \varphi_j(\xi) = 1$ for every $\xi\in \re^d$\,.
\end{description}
\end{definition}
\noindent


\begin{definition}\label{inhom} Let 
$\{\varphi_j\}_{j= 0}^{\infty} \in \Pi(\R)$ and put 
$\widehat{\Phi}_j = \varphi_j$ for $j\in\N_0$.
Let further $w \in \mathcal{W}^{\alpha_3}_{\alpha_1,\alpha_2}$ with associated weight sequence
$\{w_j\}_{j\in \n}$ defined as in \eqref{eqdef:wj}.
\begin{description}
 \item(i) For $p\in\P$, $\qconst\in(0,\infty]$, we define 
 $B^w_{\p,\qconst}(\R) = \Big\{f\in \mathcal{S}'(\R):
    \|f|B^w_{\p,\qconst}(\R)\| <\infty\Big\} $ with
 \begin{equation}\nonumber
  \begin{split}
    \|f|B^w_{\p,\qconst}(\R)\| = \Big(\sum\limits_{j=0}^{\infty}
    \|w_j(\cdot)(\Phi_j \ast f)(\cdot)|L_\p(\R)\|^{\qconst}\Big)^{1/\qconst}. 
  \end{split}
 \end{equation}
 \item(ii) For $p,q\in\P$ we define 
   $F^w_{\p,\q}(\R) = \Big\{f\in \mathcal{S}'(\R):
    \|f|F^w_{\p,\q}(\R)\| <\infty\Big\}$ with
 \begin{equation}\nonumber
  \begin{split}
    \|f|F^w_{\p,\q}(\R)\| = \Big\|\Big(\sum\limits_{j=0}^{\infty}
    |w_j(\cdot)(\Phi_j \ast f)(\cdot)|^{q(\cdot)}\Big)^{1/q(\cdot)}|L_\p(\R)\Big\|. 
  \end{split}
 \end{equation}
\end{description}
\end{definition}
\begin{remark}
It is also possible to consider Besov spaces $B^w_{\p,\q}(\R)$ with variable
index $\q$, which were introduced and studied in \cite{AlmeidaHasto2010}. The
definition of these spaces is very technical since
they require a new modular. Surprisingly it is much harder to work with Besov
spaces with variable indices $\p$ and $\q$ than to work with variable
Triebel-Lizorkin spaces, in sharp contrast to the constant exponent case. For
example, Besov spaces with variable $\q$ are not always normed spaces for
$\min\{\p,\q\}\geq1$, even if $\p$ is a constant (see \cite{KeVybNorm} for
details). So we restrict our studies on Besov spaces to the case were the index
$\q$ remains a constant $\qconst$ and we leave the fully variable case for further
research.
\end{remark}
Formally, the definition of $F^{w}_{\p,\q}(\R)$ and $B^{w}_{\p,\qconst}(\R)$ depends on the chosen decomposition of unity $\{\varphi_j\}_{j= 0}^{\infty} \in \Pi(\R)$.
The following characterization by local means shows that under certain regularity conditions on the indices $\p,\q$ it is in fact independent, in the sense of equivalent quasi-norms.

To get useful further characterizations of the spaces defined above we need a
replacement for the classical Fefferman-Stein maximal inequality since it does
not hold in our case if $q(\cdot)$ is non-constant. We will use the following
convolution inequality.
\begin{lemma}[Theorem 3.2 in \cite{DieningHastoRoudenko2009}]\label{lem:FaltungsUnglg}
	Let $p,q\in\Plog$ with $1<p^-\leq p^+<\infty$ and $1<q^-\leq q^+<\infty$, then for $m>d$ there exists a constant $c>0$ such that
	\begin{align*}
		\norm{\norm{\left(\eta_{\nu,m}\ast f_\nu\right)_{\nu\in\N_0}}{\ell_\q}}{\Lpp}\leq c\norm{\norm{\left(f_\nu\right)_{\nu\in\N_0}}{\ell_\q}}{\Lpp},
	\end{align*}
where $\eta_{\nu,m}(x)=2^{\nu d}(1+2^\nu|x|)^{-m}$.
\end{lemma}


\subsection{Continuous local means characterization}
\label{clm}

For our purpose, it is more convenient to reformulate Definition~\ref{inhom} in terms of a continuous characterization, where the discrete dilation parameter $j\in \n$ is replaced by $t>0$ and
the sums become integrals over $t$. Characterizations of this type have some history and are usually referred to as characterizations via (continuous) local means.
For further references and some historical facts we mainly refer to \cite{Tr92, BuPaTa96, Ry99a} and in particular to the recent contribution \cite{T10}, which provides a complete and self-contained reference.

The system $\{\varphi_j\}_{j\in\n} \in\Pi(\R)$ may be replaced by a more general one. Essential are functions $\Phi_0, \Phi \in \mathcal{S}(\R)$ satisfying
the so-called Tauberian conditions
\begin{equation}\label{condphi1}
\begin{split}
  |\widehat{\Phi}_0(\xi)| > 0 \quad &\mbox{ on }\quad \{|\xi| < 2\varepsilon\}\,,\\
  |\widehat{\Phi}(\xi)| > 0 \quad&\mbox{ on }\quad \{\varepsilon/2<|\xi|< 2\varepsilon\}\,,
\end{split}
\end{equation}
for some $\varepsilon>0$, and -- for some $R+1 \in \n $ -- the moment conditions
\begin{equation}\label{condphi2}
   D^{{\beta}} \widehat{\Phi}(0) = 0\quad\mbox{for all}\quad
   |\beta|_1 \leq R\,.
\end{equation}
If $R+1 = 0$ the condition \eqref{condphi2} is void. We will call the functions
$\Phi_0$ and $\Phi$ \emph{kernels for local
means} and use the notations $\Phi_k = 2^{kd}\Phi(2^{k}\cdot)$, $k\in \N$, as
well as
$\Phi_t = \mathcal{D}_t \Phi = t^{-d}\Phi(\cdot/t)$ for $t>0$. The associated \emph{Peetre maximal function}
\begin{equation}\label{Peemax}
    (\Phi^{\ast}_t f)_{a}(x) = \sup\limits_{y \in \R}\frac{|(\Phi_t \ast f)(x+y)|}
    {(1+|y|/t)^{a}}\quad,\quad x\in \R\,,t>0\,,
\end{equation}
was introduced in \cite{Pe75} for $f\in \mathcal{S}^\prime(\R)$ and $a>0$. We also need the stronger version
$$
    \langle\Phi^{\ast}_t f\rangle_{a}(x) = 
    \sup\limits_{\substack{\frac t2\leq\tau\leq 2t\\ \tau<1}} (\Phi^{\ast}_\tau
f)_{a}(x) \quad,\quad x\in \R\,,t>0\,, \quad\text{(Convention:
$\sup\emptyset=0$)}
$$
which we will refer to as
\emph{Peetre-Wiener maximal function} and which was utilized for the coorbit characterization of the classical Besov-Lizorkin-Triebel-spaces in~\cite{Sch12}.
To adapt to the inhomogeneous setting we further put
$
\langle\Phi_0^*f\rangle_a=(\Phi_0^*f)_a= ((\Phi_0)_1^*f)_a 
$.

Using these maximal functions we now state several different characterizations.

\begin{Theorem}\label{thm:contchar} Let $w \in \mathcal{W}^{\alpha_3}_{\alpha_1,\alpha_2}$ and choose functions $\Phi_0,\Phi \in \mathcal{S}(\R)$
satisfying \eqref{condphi1} and \eqref{condphi2} with $R+1>\alpha_2$. For $x\in\R$ and $t\in(0,1)$ define
$A_1f(x,t):=(\Phi_t \ast f)(x)$, $A_2f(x,t):=(\Phi_t^\ast f)_a(x)$, and $A_3f(x,t):=\langle \Phi_t^\ast f\rangle_a(x)$, $a>0$.
Further, put $A_1f(x,\infty):=(\Phi_0 \ast f)(x)$, $A_2f(x,\infty):=(\Phi_0^\ast f)_a(x)$, and $A_3f(x,\infty):=\langle \Phi_0^\ast f\rangle_a(x)$.

\begin{description}
 \item(i) If $p\in\Plog$, $0<\qconst\leq\infty$, and $a>\frac{d}{p^-}+\alpha_3$ then
 $$
    B^w_{\p,\qconst}(\R) = \{f\in \mathcal{S}'(\R)~:~\|f|B^w_{\p,\qconst}(\R)\|_i < \infty\}\quad,\quad i=1,2,3,4,
 $$
 where for $i=1,2,3$
 \begin{flalign*}
     && \|f|B^w_{\p,\qconst}(\R)\|_i &= \|w(\cdot,\infty)A_if(\cdot,\infty)|L_\p(\R)\| &\\
     &&& + \Big(\int_{0}^1
     \|w(\cdot,t)A_if(\cdot,t)|L_\p(\R)\|^\qconst\frac{dt}{t}\Big)^{1/\qconst}\,, & \\
     \text{and} &&
      \|f|B^w_{\p,\qconst}(\R)\|_4 &= \|w(\cdot,\infty)(\Phi_0^{\ast}f)_a(\cdot)|L_\p(\R)\|&\\
     &&&+ \Big(\sum\limits_{j=1}^{\infty}
    \Big\|w_j(\cdot)(\Phi_{2^{-j}}^{\ast} f)_a(\cdot)|L_\p(\R)\Big\|^{\qconst}\Big)^{1/\qconst}\,. &
 \end{flalign*}
Moreover, $\|\cdot|B^w_{\p,\qconst}(\R)\|_i$, $i=1,2,3,4$, are equivalent quasi-norms in $B^w_{\p,\qconst}(\R)$\,.
\item(ii) If $p,q\in\Plog$ with $0<q^-\leq q^+<\infty$, $0<p^-\leq p^+ < \infty$, and $a>\max\{\frac{d}{p^-},\frac{d}{q^-}\}+\alpha_3$ then
 $$
    F^w_{\p,\q}(\R) = \{f\in \mathcal{S}'(\R)~:~\|f|F^w_{\p,\q}(\R)\|_i < \infty\}\quad,\quad i=1,2,3,4,
 $$
 where for $i=1,2,3$
 \begin{flalign*}
     &&\|f|F^w_{\p,\q}(\R)\|_i &= \|w(\cdot,\infty)A_if(\cdot,\infty)|L_\p(\R)\| &\\
     &&&+ \Big\|\Big(\int_{0}^1
|w(\cdot,t)A_if(\cdot,t)|^{q(\cdot)}\frac{dt}{t}\Big)^{1/q(\cdot)}
|L_\p(\R)\Big\|\,, &\\
\text{and} &&\|f|F^w_{\p,\q}(\R)\|_4 &= \|w(\cdot,\infty)(\Phi_0^{\ast}f)_a(\cdot)|L_\p(\R)\|&\\
     &&&+ \Big\|\Big(\sum\limits_{j=1}^{\infty}
    |w_j(\cdot)(\Phi_{2^{-j}}^{\ast} f)_a(\cdot)|^{q(\cdot)}\Big)^{1/q(\cdot)}|L_\p(\R)\Big\|. 
    \nonumber
 \end{flalign*}
Moreover, $\|\cdot|F^w_{\p,\q}(\R)\|_i$, $i=1,2,3,4$, are equivalent quasi-norms in $F^w_{\p,\q}(\R)$\,.
\end{description}

\end{Theorem}

Before we present a sketch of the proof recall an important convolution inequality from \cite{Ke09}.
\begin{lemma}
	\label{lem:ConvIneq}
    Let $0<\qconst\leq\infty$, $\delta>0$ and $p,q\in\P$. Let $(g_k)_{k\in\N_0}$
    be a sequence of non-negative measurable functions on $\R$
    and denote $G_\ell=\sum_{k=0}^\infty2^{-|\ell-k|\delta}g_k$ for $\ell\in\N_0$.
    Then there exist constants $C_1,C_2\geq0$ such that
    \begin{align*}
        \norm{\{G_\ell\}_{\ell}}{\ell_{\qconst}(L_{p(\cdot)})}\leq C_1\norm{\{g_k\}_{k}}{\ell_{\qconst}(L_{p(\cdot)})} \quad\text{and}\quad
        \norm{\{G_\ell\}_{\ell}}{L_{p(\cdot)}(\ell_\q)}\leq C_2\norm{\{g_k\}_{k}}{L_{p(\cdot)}(\ell_\q)}\punkt
    \end{align*}
\end{lemma}

\noindent{\bf Proof of Theorem \ref{thm:contchar}.} We only prove (ii) and
comment afterwards briefly on the necessary modifications
for (i). The arguments are more or less the same as in the proofs of
\cite[Thm.\ 2.6]{T10} and \cite[Thm. 9.6]{Sch12}. We remark that the
equivalences $\|\cdot|F^w_{\p,\q}(\R)\|\asymp\|\cdot|F^w_{\p,\q}(\R)\|_4$ and
$\|\cdot|B^w_{\p,\qconst}(\R)\|\asymp\|\cdot|B^w_{\p,\qconst}(\R)\|_4$ are
already known, see \cite{Ke09}.

\noindent

{\em Step 1.} First, we prove a central estimate \eqref{eq:essestimate} between different start functions $\Phi$ and $\Psi$ incorporating the different types of Peetre maximal operators. The needed norm inequalities in the theorem are consequences of this central estimate \eqref{eq:essestimate}, and are  subsequently deduced in the following steps.

Let us put $\varphi_0:=\widehat{\Phi}_0$ and $\varphi_k:=\widehat{\Phi}_k$ for $k\in\N$.
We can find a pair of functions $\lambda_0,\lambda\in\mathcal{S}(\R) $ with
$\supp \lambda_0\subset \{ \zeta\in\R : |\zeta|\le 2\varepsilon \}$ and
$\supp \lambda\subset \{ \zeta\in\R : \varepsilon/2 \le |\zeta| \le 2\varepsilon\}$ such that
$\sum_{k\in\N_0} \lambda_k \varphi_k \equiv 1$, where $\lambda_k=\lambda(2^{-k}\cdot)$ for $k\in\N$.
Let us shortly
demonstrate how to do that. We use the special dyadic
decomposition of unity given by $\eta_0(t) = 1$ if $|t|\leq 4/3$ and
$\eta_0(t) =0$ if $|t|>3/2$. We
put $\eta_k:=\eta_0(\cdot/2^k)-\eta_0(\cdot/2^{k-1})$ for $k\in\N$. Then
clearly $\eta_0 + \sum_{k=1}^\infty \eta_k \equiv 1$ and we obtain $\sum_{k\in\N_0} \lambda_k \varphi_k \equiv 1$ by
defining $\lambda_k := \eta_k(\cdot/\varepsilon)/\varphi_k$ for $k\in\n$ and $\lambda:=\lambda_1(2\cdot)$.

The support of the function $\theta:=1-\sum_{k\in\N} \lambda_k\varphi_k\in C^\infty_0(\R)$ is
fully contained in $M:=\{ |x| \le 3\varepsilon/2 \}$. Due to
the Tauberian conditions, $\varphi_0$ is positive on $M$. Inverting $\varphi_0$ on $M$ and extending appropriately outside, we
can construct a function $\gamma\in C^\infty_0(\R)$, which coincides with $1/\varphi_0$ on $M$.
Since $\lambda_0\varphi_0=\theta$ we thus have the factorization $\lambda_0=\gamma \theta$.

We now put $\lambda_{0,u}(\cdot):=\gamma(\cdot) \theta(u\cdot)$ for $u\in[1,2]$, which gives
\begin{align*}
\lambda_{0,u} \varphi_0 + \sum_{k\in\N} \lambda_k(u\cdot) \varphi_k(u\cdot) =1.
\end{align*}
We then define $\Xi$, $\Theta$, $\Lambda$, $\Lambda_{0,u}$, and $\Lambda_k$ for $k\in\n$, all elements of $\mathcal{S}(\R)$, via inverse Fourier transform of the
functions  $\gamma$, $\theta$, $\lambda$, $\lambda_{0,u}$, and $\lambda_k$, respectively.
We get $\Lambda_{0,u}=\Xi \ast \Theta_u$ and it holds
$
g= \Lambda_{0,u} \ast \Phi_0 \ast g + \sum_{k\in\N} \Lambda_{2^{-k}u} \ast \Phi_{2^{-k}u} \ast g
$
for every $g\in\mathcal{S}^\prime(\R)$.

Let $\Psi_0,\Psi\in \mathcal{S}(\R)$ be another system 
which satisfies the Tauberian conditions \eqref{condphi1} and \eqref{condphi2}.
Choosing $g=\Psi_{2^{-\ell}v}\ast f$, where $f\in\mathcal{S}^\prime(\R)$, $\ell\in\N$, and $v\in[1/2,4]$,
we get
\begin{align}\label{eq:convident}
\Psi_{2^{-\ell}v}\ast f=\sum_{k\in\N} \Psi_{2^{-\ell}v} \ast \Lambda_{2^{-k}u} \ast \Phi_{2^{-k}u} \ast f + \Psi_{2^{-\ell}v} \ast \Lambda_{0,u}\ast \Phi_0 \ast f.
\end{align}

Defining $J_{\ell,k}= \int_{\R} |\Psi_{2^{-\ell}v} \ast \Lambda_{2^{-k}u}(z)| (1+2^k|z|/u)^a \,dz$ for $k\in\N$ we have
for $y\in\R$
\begin{gather*}
\begin{aligned}
|(\Psi_{2^{-\ell}v} \ast \Lambda_{2^{-k}u} \ast \Phi_{2^{-k}u} \ast f)(y)|
&\le \int_{\R} |\Psi_{2^{-\ell}v} \ast \Lambda_{2^{-k}u}(z)| | \Phi_{2^{-k}u} \ast f(y-z)  | \,dz\\
&\le (\Phi_{2^{-k}u}^*f)_a(y) J_{\ell,k},
\end{aligned}
\end{gather*}
For $k=0$ we get with $J_{\ell,0}=\int_{\R} |\Psi_{2^{-\ell}v} \ast \Lambda_{0,u}(z)| (1+|z|)^a \,dz$
\begin{gather*}
\begin{aligned}
|(\Psi_{2^{-\ell}v} \ast \Lambda_{0,u} \ast \Phi_0 \ast f)(y)|
\le \int_{\R} |\Psi_{2^{-\ell}v} \ast \Lambda_{0,u}(z)| | \Phi_0 \ast f(y-z)  | \,dz
\le (\Phi_0^*f)_a(y) J_{\ell,0}.
\end{aligned}
\end{gather*}

To estimate $J_{\ell,k}$ the following identity for functions $\mu,\nu\in\mathcal{S}(\R)$ is used,
\begin{equation*}
(\mu_u\ast\nu_v)(x)= \frac{1}{u^d} [\mu\ast\nu_{v/u}](x/u) = \frac{1}{v^d} [\mu_{u/v}\ast\nu](x/v),
\end{equation*}
valid for $u,v>0$ and $x\in\R$. In case $\ell\ge k>0$ we obtain
\begin{align*}
J_{\ell,k}= \int_{\R} |(\Psi_{2^{k-\ell}\frac{v}{u}} \ast \Lambda)(z)|(1+|z|)^a \,dz
\lesssim \sup_{z\in\R} \big| (\Psi_{2^{k-\ell}\frac{v}{u}}\ast \Lambda)(z)
(1+|z|)^{a+d+1}  \big| \lesssim  2^{(k-\ell)(R+1)}\,,
\end{align*}
where we used \cite[Lemma~1]{Ry99a} in the last step.
In case $0<\ell< k$ we estimate similarly to obtain
\begin{align*}
J_{\ell,k}= \int_{\R} |(\Psi \ast \Lambda_{2^{-(k-\ell)}{u/v}}(z))|(1+2^{k-\ell}u|z|/v)^a \,dz
\lesssim  2^{(\ell-k)(L+1-a)},
\end{align*}
where $L$ can be chosen arbitrarily large since $\Lambda\in \mathcal{S}_0(\R)$ fulfills moment conditions for all $L\in\N_0$.

For $\ell>k=0$ we estimate as follows, taking advantage of $\Xi\in\mathcal{S}(\R)$,
\begin{align*}
J_{\ell,0}&= \int_{\R} |(\Psi_{2^{-\ell}v} \ast
\Theta_u) \ast \Xi (z)|(1+|z|)^a \,dz \\
&\lesssim \sup_{y\in\R} \big| (\Psi_{2^{-\ell}v}\ast
\Theta_u)(y) (1+|y|)^{a+d+1}  \big|  \int_{\R} \int_{\R} |\Xi(z-y)| (1+|z-y|)^a  (1+|y|)^{-d-1} \,dz dy  \\
&\lesssim \sup_{y\in\R} \big| (\Psi_{2^{-\ell}v/u}\ast
\Theta)(y) (1+|y|)^{a+d+1}  \big|  \int_{\R} \int_{\R}  (1+|z|)^{-d-1}  (1+|y|)^{-d-1} \,dz dy  \lesssim  2^{-\ell(R+1)}.
\end{align*}

Using $1+t|x|\le \max\{1,t\} (1+|x|)$ and $1+|x+y|/t\le (1+|y|/t) (1+|x|/t)$
for $t>0$ and $x,y\in\R$ we further deduce for $k\in\N$
\begin{align*}
(\Phi^*_{2^{-k}u}f)_a(y) &\le (\Phi^*_{2^{-k}u}f)_a(x)(1+2^k|x-y|/u)^a \\
&\lesssim (\Phi^*_{2^{-k}u}f)_a(x)(1+2^\ell|x-y|/v)^a \max\{1,2^{(k-\ell)}\}^a.
\end{align*}
and $(\Phi_0^*f)_a(y) \lesssim (\Phi_0^*f)_a(x) (1+2^\ell|x-y|/v)^a$.
Altogether, we arrive - for $k\ge 1$ - at
\begin{align*}
\sup_{y\in\R} \frac{|(\Psi_{2^{-\ell}v} \ast \Lambda_{2^{-k}u} \ast (\Phi_{2^{-k}u} \ast f))(y)|}{(1+2^\ell|x-y|/v)^a}
\lesssim  (\Phi^*_{2^{-k}u}f)_a(x)
\begin{cases} 2^{(k-\ell)(R+1)} ~&:~ \ell\ge k, \\
2^{(\ell-k)(L+1-2a)}  ~&:~ \ell< k,
\end{cases}
\end{align*}
with an implicit constant independent of $u\in[1,2]$ and $v\in [1/2,4]$.
For $k=0$ we obtain
\begin{align*}
\sup_{y\in\R} \frac{|(\Psi_{2^{-\ell}v} \ast \Lambda_{0,u} \ast \Phi_0 \ast f)(y)|}{(1+2^\ell|x-y|/v)^a}
\lesssim  (\Phi^*_{0}f)_a(x) 2^{-\ell(R+1)}.
\end{align*}
We thus conclude from \eqref{eq:convident} that uniformly in $t,u\in[1,2]$
\begin{align*}
\langle \Psi^*_{2^{-\ell}t}f \rangle_a(x) &= \sup_{t/2\le v \le 2t, v<1} (\Psi^*_{2^{-\ell}v}f)_a(x) \\ \notag
&\lesssim (\Phi_0^\ast f)_a(x) 2^{-\ell(R+1)}
+ \sum_{k\in\N}  (\Phi^*_{2^{-k}u}f)_a(x) \begin{cases} 2^{(k-\ell)(R+1)} \,&: \ell\ge k, \\
2^{(\ell-k)(L+1-2a)}  &: \ell< k.
\end{cases}
\end{align*}
Writing $\tilde{w}_{\ell,t}(x)=w(x,2^{-\ell}t)$ for $\ell\in\N$ and $\tilde{w}_{0,t}(x)=w(x,\infty)$ we have
\[
\tilde{w}_{\ell,t}(x) \tilde{w}_{k,u}(x)^{-1}  \lesssim  \begin{cases}  2^{(\ell-k)\alpha_2} \quad &\ell\ge k, \\
 2^{(\ell-k)\alpha_1}  &\ell< k,
\end{cases}
\]
as a consequence of $(W1)$, \eqref{eq_st1}, and \eqref{eq_st2}.
Multiplying both sides with $w(x,2^{-\ell}t)$
we finally derive with an implicit constant independent of $t,u\in[1,2]$
\begin{align*}
	w(x,2^{-\ell}t) \langle \Psi_{2^{-\ell}t}^{\ast} f \rangle_a (x)
	&\lesssim  w(x,\infty) (\Phi_0^\ast f)_a(x) 2^{-\ell(R+1-\alpha_2)} \\
	&+\sum_{k\in\N} w(x,2^{-k}u) (\Phi_{2^{-k}u}^{\ast}f)_a(x)
\begin{cases} 2^{(k-\ell)(R+1-\alpha_2)} \,&:~ \ell\ge k, \\
2^{(\ell-k)(L+1-2a+\alpha_1)}  &:~ \ell< k.
\end{cases}
\end{align*}

Choosing $L\ge 2a - \alpha_1$ 
we have with $0<\delta = \min\{1,R+1-\alpha_2\}$ the central estimate
%

\begin{align}\label{eq:essestimate}
   \langle \Psi_{2^{-\ell}t}^{\ast} f \rangle_a (x)
	\lesssim 2^{-\ell\delta} \frac{w(x,\infty)}{ w(x,2^{-\ell}t)} (\Phi_0^\ast f)_a(x) + \sum_{k\in\N}  2^{-|k-\ell|\delta} \frac{w(x,2^{-k}u)}{ w(x,2^{-\ell}t)} (\Phi_{2^{-k}u}^{\ast}f)_a(x) .
\end{align}

{\em Step 2.} We show $\|f|F^w_{\p,\q}(\R)\|_1 \asymp \|f|F^w_{\p,\q}(\R)\|_{2,3}$. The direction $\|f|F^w_{\p,\q}(\R)\|_1 \lesssim \|f|F^w_{\p,\q}(\R)\|_{2,3}$ is obvious
and it remains to verify $\|f|F^w_{\p,\q}(\R)\|_{3} \lesssim \|f|F^w_{\p,\q}(\R)\|_1$.\\
We use \eqref{eq:essestimate} with $\Psi=\Phi$. Choosing $0<\tilde{\delta} \le \delta$ we obtain for any $r>0$, using an embedding argument if $0<r\le1$ and Hölder's inequality otherwise,
\begin{align*}
   \langle \Phi_{2^{-\ell}t}^{\ast} f \rangle_a^r (x)w^r(x,2^{-\ell}t)
	\lesssim 2^{-\ell\tilde{\delta} r} w^r(x,\infty) (\Phi_0^\ast f)_a^r(x) + \sum_{k\in\N}  2^{-|k-\ell|\tilde{\delta} r} w^r(x,2^{-k}u) (\Phi_{2^{-k}u}^{\ast}f)^r_a(x) .
\end{align*}
To estimate the sum on the right hand side we use (2.66) proved in Substep~1.3 of the proof of \cite[Thm.\ 2.6]{T10}. It states that for
$x\in\R$, $f\in\mathcal{S}^\prime(\R)$, $k\in\N$, $u\in[1,2)$, $r>0$, and
$0<a\le N$ for some arbitrary but fixed $N\in\n$
\begin{equation}\label{2.30}
    (\Phi^{\ast}_{2^{-k}u}f)_{a}(x)^r\le C_N
    \sum\limits_{j\in \n} 2^{-jNr}2^{(k+j)d}
    \int_{\R}\frac{|(\Phi_{2^{-(k+j)}u}\ast
    f)(y)|^r}{(1+2^{k}|x-y|)^{a
    r}}\,dy\,,
\end{equation}
where the constant $C_N$ is independent of $x,f,k$, and $u\in[1,2)$,
but may depend on $r$, $a$ and $N$. Taking into account $(W2)$ and \eqref{eq_W1tilde}, which give the relation
$w(x,2^{-k}u) \lesssim 2^{-j\alpha_1}(1+2^{k}|x-y|)^{\alpha_3}w(y,2^{-(j+k)}u)$ and $(1+2^k|z|)^{-M}\leq2^{jM}(1+2^{k+j}|z|)^{-M}$, this leads to
\begin{align}\notag
 &\langle \Phi_{2^{-\ell}t}^{\ast} f \rangle_a^r (x)w^r(x,2^{-\ell}t)
	\lesssim 2^{-\ell\tilde{\delta} r} w^r(x,\infty) (\Phi_0^\ast f)_a^r(x)\\
	&\hspace{5em}+ \sum_{k\in\N}  2^{-|k-\ell|\tilde{\delta} r}\sum_{j\in\n}2^{-jr\tilde{N}}2^{(k+j)d}\int_{\R}\frac{|(\Phi_{2^{-(k+j)}u}\ast f)(y)w(y,2^{-(j+k)}u)|^r}{(1+2^{k+j}|x-y|)^{(a-\alpha_3)r}}\,dy\,\label{eq:step2_eq1}
\end{align}
with $\tilde{N}=N-a+\alpha_1+\alpha_3>0$. Since $x\in\R$ is fixed we can apply in $t$ the $L_{q(x)/r}([1,2);\frac{dt}{t})$ norm with $r<\min\{p^-,q^-\}$. This changes only the constant and the
left-hand side of \eqref{eq:step2_eq1}.
The $L_{q^-/r}([1,2);\frac{du}{u})$ (quasi-)norm in the variable $u$ only affects the right-hand side of \eqref{eq:step2_eq1}.
With Minkowski's integral inequality we obtain
\begin{align}
	&\left(\int_1^2|\langle\Phi_{2^{-\ell}t}^{\ast} f \rangle_a(x)w(x,2^{-\ell}t)|^{q(x)}\frac{dt}{t}\right)^{r/q(x)}-2^{-\ell\tilde{\delta} r} w^r(x,\infty) (\Phi_0^\ast f)_a^r(x)\notag\\
	&\hspace{2em}\lesssim\sum_{k\in\N}  2^{-|k-\ell|\tilde{\delta} r}\sum_{j\in\n}2^{-|j-k|\tilde{N}r} 2^{jd} \int_\R\frac{\left(\int_1^2|(\Phi_{2^{-j}u}\ast f)(y)w(y,2^{-j}u)|^{q^-}\frac{du}{u}\right)^{r/q^-}}{(1+2^{j}|x-y|)^{(a-\alpha_3)r}}dy\notag\\
	&\hspace{2em}\lesssim \sum_{k\in\N}2^{-|k-\ell|\tilde{\delta} r}\sum_{j\in\n}2^{-|j-k|\tilde{N}r}\left[\eta_{j,(a-\alpha_3)r}\ast\left(\int_1^2|(\Phi_{2^{-j}u}\ast f)(\cdot)w(\cdot,2^{-j}u)|^{{q}^-}\frac{du}{u}\right)^{r/{q}^-}\right](x)\ \notag
\end{align}
with functions $\eta_{\nu,m}(x)=2^{\nu d}(1+2^\nu|x|)^{-m}$.\\
Now we choose $r>0$ such that $\frac{d}{a-\alpha_3} < r < \min\{p^-,q^-\}$,
which is possible since $a>\alpha_3 + \frac{d}{\min\{p^-,q^-\}}$, and $N$ such that $\tilde{N}>0$.
Applying the $L_{\p/r}(\ell_{\q/r})$ norm with respect to $x\in\R$ and $\ell\in\N$ and using Lemma~\ref{lem:ConvIneq} twice together with Lemma~\ref{lem:FaltungsUnglg} (note $(a-\alpha_3)r>d$) then yields
\begin{align}
	&\norm{\left(\int_1^2|\langle\Phi_{2^{-\ell}t}^*f\rangle_a(\cdot)w(\cdot,2^{-\ell}t)|^{q(\cdot)}\frac{dt}{t}\right)^{r/q(\cdot)}}{L_{\p/r}(\ell_{\q/r})}-c\norm{w(\cdot,\infty) (\Phi_0^\ast f)_a(\cdot)}{L_\p(\R)}^r\notag\\
	&\hspace{2em}\lesssim\norm{\left[\eta_{\ell,(a-\alpha_3)r}\ast\left(\int_1^2|(\Phi_{2^{-\ell}u}\ast f)(\cdot)w(\cdot,2^{-\ell}u)|^{{q}^-}\frac{du}{u}\right)^{r/{q}^-}\right](x)}{L_{\p/r}(\ell_{\q/r})}\notag\\
	&\hspace{2em}\lesssim\norm{\left(\int_1^2|(\Phi_{2^{-\ell}u}\ast f)(\cdot)w(\cdot,2^{-\ell}u)|^{{q}^-}\frac{du}{u}\right)^{r/{q}^-}}{L_{\p/r}(\ell_{\q/r})}\ .\label{eq:step2_eq2}
\end{align}
Finally, we use H\"older's inequality to estimate the integral in the last norm. We use $0< q^-\leq q(x)$ and get
\begin{align*}
 &\left(\int_1^2|(\Phi_{2^{-\ell}u}\ast f)(x)w(x,2^{-\ell}u)|^{{q}^-}\frac{du}{u}\right)^{r/q^-}\\
 &\hspace{5em}\leq\left(\int_1^2|(\Phi_{2^{-\ell}u}\ast f)(x)w(x,2^{-\ell}u)|^{q(x)}\frac{du}{u}\right)^{r/q(x)}\left(\int_1^2\frac{du}{u}\right)^{r/{q}^-\cdot\frac{1}{\left(\frac{{q(x)}}{{q}^-}\right)'}}\\
 &\hspace{5em}\leq
\left(\int_1^2|(\Phi_{2^{-\ell}u}\ast f)(x)w(x,2^{-\ell}u)|^{q(x)}\frac{du}{u}\right)^{r/q(x)}.
\end{align*}
Using this estimate we can reformulate \eqref{eq:step2_eq2} into
\begin{align*}
 &\norm{\left(\int_0^1|\langle\Phi_\lambda^*f\rangle_a(\cdot)w(\cdot,\lambda)|^{q(\cdot)}\frac{d\lambda}{\lambda}\right)^{1/q(\cdot)}}{L_\p}\\
	&\hspace{5em}\lesssim\norm{w(\cdot,\infty) (\Phi_0^\ast f)_a(\cdot)}{L_\p(\R)}+\norm{\left(\int_0^1|(\Phi_\lambda\ast f)(\cdot)w(\cdot,\lambda)|^{q(\cdot)}\frac{d\lambda}{\lambda}\right)^{1/q(\cdot)}}{L_\p}\ .
\end{align*}
The inhomogeneous term $(\Phi_0^\ast f)_a(x)$ needs to be treated separately. The argumentation, however, is analogous to the exposition before with \eqref{2.30} replaced by the inequality
\begin{equation*}
    (\Phi^{\ast}_0f)_{a}(x)^r\lesssim
    \sum\limits_{k\in \N} 2^{-kNr}2^{kd}
    \int_{\R}\frac{
    |(\Phi_{2^{-k}u}\ast f)(y)|^r}{(1+|x-y|)^{ar}}\,dy + \int_{\R}\frac{
    |(\Phi_{0}\ast f)(y)|^r}{(1+|x-y|)^{ar}}\,dy.
\end{equation*}

In the Besov space case we do not need the functions $\eta_{\nu,m}$ and one can
work with the usual maximal operator $\HLM$ together with Lemma~\ref{lem:HLM},
see \cite{Ke09} for details.

\noindent
\indent {\em Step 3.} In the third step we show $\|f|F^w_{\p,\q}(\R)\|_2 \asymp \|f|F^w_{\p,\q}(\R)\|_3 \asymp \|f|F^w_{\p,\q}(\R)\|_4 $. We immediately observe
$\|f|F^w_{\p,\q}(\R)\|_2 \lesssim \|f|F^w_{\p,\q}(\R)\|_3$. \\
\noindent
\indent {\em Substep 3.1.} To prove $\|f|F^w_{\p,\q}(\R)\|_3 \lesssim \|f|F^w_{\p,\q}(\R)\|_4$ we apply \eqref{eq:essestimate} with $u=1$ and $\Psi=\Phi$ .
Since the inhomogeneous terms are identical, it suffices to estimate the
homogeneous part. Integration with respect to $dt/t$ yields for $\ell\in\N$
\begin{align*}
\Big(  \int_1^2 |w(x,2^{-\ell}t) \langle \Phi_{2^{-\ell}t}^{\ast} f\rangle_a (x) |^{q(x)} \frac{dt}{t} \Big)^{\frac{1}{q(x)}}
\lesssim 2^{-\ell\delta}w_0(x) (\Phi_0^{\ast}f)_a(x) + \sum_{k\in\N} 2^{-|k-\ell|\delta}   w_k(x) (\Phi_{2^{-k}}^{\ast}f)_a(x).
\end{align*}
Let us denote the function on the right-hand side of the previous estimate by $G_\ell$.
Applying the vector-valued convolution inequality of Lemma~\ref{lem:ConvIneq} then proves
\begin{align*}
&\Big\| \Big( \sum_{\ell=1}^\infty \int_1^2 | w(\cdot,2^{-\ell}t)  \langle\Phi^*_{2^{-\ell}t}f\rangle_a(\cdot)|^{q(\cdot)}
\,\frac{dt}{t}   \Big)^{1/q(\cdot)} \big| L_{\p} \Big\| \lesssim  \| \{G_\ell \}_{\ell \in \N} | L_{\p}(\ell_{\q}) \| \\
& \qquad \lesssim \| w_0 (\Phi_0^{\ast}f)_a | L_\p \| + \| \{ w_k (\Phi_{2^{-k}}^{\ast}f)_a \}_{k\in\N} | L_\p(\ell_\q) \|
=\| f | F^w_{\p,\q} \|_4.
\end{align*}

{\em Substep 3.2:} Let us finish by proving $\|f|F^w_{\p,\q}(\R)\|_4 \lesssim \|f|F^w_{\p,\q}(\R)\|_2$.
Again it suffices to estimate the homogeneous part.
For this we let $t=1$ and $\Psi=\Phi$ in \eqref{eq:essestimate}. If $q(x)\ge1$ we can use Minkowski's inequality to deduce
\begin{align*}
w_\ell(x) (\Phi_{2^{-\ell}}^{\ast}f)_a(x) &\lesssim 2^{-\ell\delta} w(x,\infty)
(\Phi_0^\ast f)_a(x)\\
&+ \sum_{k\in\N} 2^{-|k-\ell|\delta} \Big( \int_1^2
|\tilde{w}_{k,u}(x) (\Phi_{2^{-k}u}^{\ast}f)_a(x)|^{q(x)} \, \frac{du}{u}
\Big)^{1/q(x)}.
\end{align*}
Applying the $\ell_{q(x)}$-norm on both sides, Young's convolution inequality then yields
\begin{align}\label{eq:centest}
\sum_{\ell\in\N} w_\ell(x)^{q(x)} (\Phi_{2^{-\ell}}^{\ast}f)_a(x)^{q(x)}
&\lesssim  (w(x,\infty) (\Phi_0^\ast f)_a(x))^{q(x)} \\
&+ \sum_{k\in\N}  \int_1^2 |\tilde{w}_{k,u}(x)
(\Phi_{2^{-k}u}^{\ast}f)_a(x)|^{q(x)} \, \frac{du}{u}. \notag
\end{align}
If $q(x)<1$ we use the $q(x)$-triangle inequality
\begin{align*}
\Big(w_\ell(x) (\Phi_{2^{-\ell}}^{\ast}f)_a(x) \Big)^{q(x)} &\lesssim
2^{-\ell\delta q(x)} (w(x,\infty) (\Phi_0^\ast f)_a(x))^{q(x)}\\
&+ \sum_{k\in\N} 2^{-|k-\ell|q(x)\delta} \int_1^2 |\tilde{w}_{k,u}(x)
(\Phi_{2^{-k}u}^{\ast}f)_a(x)|^{q(x)} \, \frac{du}{u}.
\end{align*}
Now we take on both sides the $\ell_1$-norm with respect to the index $\ell\in\N$ and take into account $\sum_{k\in\n} 2^{-|k|q(x)\delta} \le C$.
We thus arrive at the same estimate \eqref{eq:centest}. Taking the $L_{\p}$-quasi-norm of \eqref{eq:centest} finishes the proof of Substep 3.2 and hence Step 3.

{\em Step 4:} Relation \eqref{eq:essestimate} also immediately allows to change to a different system $\Psi_0,\Psi$, however
in the discrete setting the change of systems has already been shown in \cite{Ke09}.
\eproof

\begin{remark} The previous theorem ensures in particular the
independence of Besov-Lizorkin-Triebel type spaces with variable exponents from
the chosen resolution of unity if
	$p,q\in\Plog$ with $p^+<\infty$, $q^+<\infty$ in the $F$-case and
$p\in\Plog$, $\qconst\in(0,\infty]$ in the $B$-case.
\end{remark}

\section{Variable exponent spaces as coorbits}
\label{sec:appcoorbit}

In order to treat the spaces $B^w_{\p,\qconst}(\R)$ and $F^w_{\p,\q}(\R)$ as coorbits
we utilize an inhomogeneous version of the continuous wavelet transform, which
uses high scale wavelets together with a base scale for the analysis.
The corresponding index set is $X = \R \times [(0,1) \cup \{\infty\}]$, where
$\infty$ denotes an isolated point, equipped with the Radon measure
$\mu$ defined by
$$
    \int_{X} F(\x) d\mu(\x) = \int_{\R}\int_{0}^1 F(x,s) \frac{ds}{s^{d+1}}dx +
\int_{\R} F(x,\infty) dx\,.
$$
The wavelet transform is then given by $V_{\cf}f(\x) = \langle f,
\varphi_\x\rangle$, $\x\in X$, for a continuous frame $\mathcal{F}=\{ \varphi_\x
\}_{\x\in X}$ on
$\mathcal{H}= L_2(\R)$ of the form
\begin{equation}\label{eqdef:wavefr}
      \varphi_{(x,\infty)} = T_x \Phi_0 = \Phi_0(\cdot-x) \quad \mbox{and} \quad \varphi_{(x,t)} =
T_x\mathcal{D}^{L_2}_t \Phi = t^{-d/2} \Phi((\cdot-x)/t)\,,
\end{equation}
with suitable functions $\Phi_0,\Phi\in L_2(\R)$. Such a frame
$\cf=\cf(\Phi_0,\Phi)$ will in our context be referred to as a \emph{continuous
wavelet frame} in $L_2(\R)$.

\begin{definition}\label{def:admfr} A continuous wavelet frame $\cf = \cf(\Phi_0,\Phi)$ is \emph{admissible} if
$\Phi_0\in\mathcal{S}(\R)$ and $\Phi\in \mathcal{S}_0(\R)$ are chosen such that they
satisfy the Tauberian conditions \eqref{condphi1}, \eqref{condphi2} and the condition
$$
    |\widehat{\Phi}_0(\xi)|^2 + \int_{0}^1 |\widehat{\Phi}(t\xi)|^2\frac{dt}{t} = C \quad\text{ for a.e. }\xi\in\R.
$$
\end{definition}

An admissible wavelet frame $\cf(\Phi_0,\Phi)$ represents a tight continuous frame in the sense of \eqref{eq:stab}. To see this, apply Fubini's and Plancherel's theorem to get
\begin{equation}
    C \|f|L_2(\R)\|^2
    = \int_{\R} |\widehat{f}(\xi)|^2 \Big(|\widehat{\Phi}_0(\xi)|^2 + \int_{0}^1 |\widehat{\Phi}(t\xi)|^2 \frac{dt}{t}\Big)\,d\xi
    = (2\pi)^{-d} \int_{X} |\langle f,\varphi_{\x}\rangle|^2 d \mu(\x)\,.\notag
\end{equation}

\subsection{Peetre-Wiener type spaces on $X$}
\label{ssec:PeetSp}

We intend to define two general scales of spaces on $X$, for which we need
a Peetre type maximal function, given for a measurable function $F: X \to \C$ by
\begin{align*}
    \mathcal{P}^\ast_a F(x,t) &:= \esssup{\substack{z\in \R, \tau<1\\
    \frac{t}{2}\le \tau \le 2t}}\frac{|F(x+z,\tau)|}{(1+|z|/\tau)^a}\quad,\quad x\in\R,\, 0<t<1,\\
    \mathcal{P}^\ast_a F(x, \infty) &:= \esssup{z\in \R} \frac{|F(x+z,\infty)|}{(1+|z|)^a} \quad , \quad x \in \R.\notag
\end{align*}
The operator $\mathcal{P}^\ast_a$ is a stronger version of the usual Peetre maximal operator $\mathcal{P}_a$, which does not take the supremum over $t$
and was used e.g.\ in \cite{RaUl10}.

\begin{definition}\label{PeetreWiener}
Let $p,q\in\Plog$ with $0< p^-\leq p^+<\infty$ and $0< q^-\leq q^+<\infty$ and let $0< \qconst\leq \infty$.
Further, let $a>0$ and $w \in \mathcal{W}^{\alpha_3}_{\alpha_1,\alpha_2}$.
Then we define by
\begin{equation*}
  \begin{split}
     \Pe(X) &= \{F:X \to \C~:~\|F|\Pe\| < \infty\}\,,\\
     \Le(X) &= \{F:X \to \C~:~\|F|\Le\| < \infty\}\\
  \end{split}
\end{equation*}
two scales of function spaces on $X$ with respective quasi-norms
\begin{align*}
    \|F|\Pe\| & := \Big\|w(\cdot,\infty)\cp^\ast_a F(\cdot, \infty) |L_{p(\cdot)}(\R)\Big\|\\
   &~~~~+ \Big\|\Big(\int_{0}^1 \Big[w(\cdot,t)\cp^\ast_a F(\cdot, t) \Big]^{q(\cdot)}\frac{dt}{t}\Big)^{1/q(\cdot)}|L_{p(\cdot)}(\R)\Big\|,\\
    \|F|\Le\| & := \Big\|w(\cdot,\infty)\cp^\ast_a F(\cdot, \infty) |L_{p(\cdot)}(\R)\Big\|\nonumber  \\
   &~~~~+ \Big(\int_{0}^1 \Big\|w(\cdot,t) \cp^\ast_a F(\cdot,t)|L_{p(\cdot)}(\R)\Big\|^{\qconst}\frac{dt}{t}\Big)^{1/\qconst}\,.
\end{align*}
\end{definition}
It is not hard to verify that in case $a>d/p^-+\alpha_3$ these spaces are rich solid QBF-spaces as defined and studied in Subsection~\ref{ssec:QBFspaces}. Moreover,
the utilization of the Peetre-Wiener operator $\cp^\ast_a$ ensures that they are locally integrable, even in the quasi-Banach case in contrast to the
ordinary Peetre spaces where $\cp_a$ is used instead of $\cp^\ast_a$. In fact, there is an associated locally bounded weight function given by
\begin{equation}\label{eqdef:assweight}
\vB(x,t) = \left\{\begin{array}{rcl}
            t^{\alpha_1-d/p^-}(1+|x|)^{\alpha_3}&,&x\in\R,\,0<t<1,\\
            (1+|x|)^{\alpha_3}&,&x\in\R,\,t=\infty,
         \end{array}\right.
\end{equation}
such that the following lemma holds true.

\begin{lemma}
We have the continuous embeddings
\[
\Pe(X) \hookrightarrow L_\infty^{1/\vB}(X) \quad\text{and}\quad \Le(X) \hookrightarrow L_\infty^{1/{\nu_{w,\p,\qconst}}}(X).
\]
\end{lemma}
\bproof
It is useful to interpret the component $\R\times(0,1)$ of the index $X$ as a subset of the $ax+b$ group $\mathcal{G}=\R\times(0,\infty)$
with multiplication $(x,t)(y,s)=(x+ty,ts)$ and $(x,t)^{-1}=(-x/t,1/t)$.
Let $U^{-1}$ be the inversion of $U:=[-2,2]^d\times[\frac{1}{2},2]$ and define $U_{(x,t)}:=(x,t)U^{-1}$ and $\widetilde{U}_{(x,t)}:=(x,t)U$.
Further put $Q_{(x,t)}:=x + t[-1,1]^d$ and $U_{(x,\infty)}:=\widetilde{U}_{(x,\infty)}:=Q_{(x,1)}\times\{\infty\}$. Then
we can estimate for $F:X\to \C$ and almost all $(x,t)\in X$ at every fixed $(y,s)\in X$
\begin{align*}
|F(x,t)| \chi_{U_{(x,t)}}(y,s) \le \esssup{(x,t)\in X\cap\widetilde{U}_{(y,s)}}|F(x,t)| \lesssim \mathcal{P}^\ast_a F(y,s).
\end{align*}
For convenience, let us introduce
\[
\| F | M^w_{\p,\q} \| :=\Big\|w(\cdot,\infty) F(\cdot, \infty) |L_{p(\cdot)}\Big\|
   + \Big\|\Big(\int_{0}^1 \Big| w(\cdot,t) F(\cdot, t) \Big|^{q(\cdot)}\frac{dt}{t}\Big)^{1/q(\cdot)}|L_{p(\cdot)}\Big\|.
\]
We obtain for almost all $(x,t)\in X$
\[
|F(x,t)|\cdot \| \chi_{U_{(x,t)}} |M^w_{\p,\q} \| \lesssim \| \mathcal{P}^\ast_a F |M^w_{\p,\q} \| = \| F | \Pe \|.
\]
It remains to prove $\vB(x,t)\gtrsim \| \chi_{U_{(x,t)}} |  M^w_{\p,\q} \|^{-1}$.
Since $U^{-1}\supset[-1,1]^d\times[\frac{1}{2},2]$ we have $U_{(x,t)}\supset Q_{(x,t)} \times [\frac{t}{2},2t]$. If $0<t<1$ it follows for $x,y\in\R$
\begin{align*}
\Big(\int_0^1 \big[ w(y,s) \chi_{U_{(x,t)}}(y,s) \big]^{q(y)} \frac{ds}{s} \Big)^{1/q(y)}
&\gtrsim \ln(4)^{1/q(y)} w(y,t)  \chi_{Q_{(x,t)}}(y)  \gtrsim  w(y,t) \chi_{Q_{(x,t)}}(y)
\end{align*}
and $\chi_{U_{(x,t)}}(y, \infty)=0$. The properties (W1) and (W2) of $w\in \mathcal{W}^{\alpha_3}_{\alpha_1,\alpha_2}$ further imply
\[
w(y,t)\gtrsim w(x,t)(1+|x-y|/t)^{-\alpha_3} \gtrsim   t^{-\alpha_1} (1+|x|)^{-\alpha_3} (1+|x-y|/t)^{-\alpha_3}.
\]
This leads to
\hfill$
\| \chi_{U_{(x,t)}} | M^w_{\p,\q}  \| \gtrsim
t^{-\alpha_1} (1+|x|)^{-\alpha_3} 
\|  \chi_{Q_{(x,t)}}(\cdot) (1+|x-\cdot|/t)^{-\alpha_3}  |L_\p \|.
$\hfill~

\noindent
Since $\|\chi_{Q_{(x,t)}} |L_\p \| \ge \frac{1}{2} \min\{ |Q_{(x,t)}|^{1/p^+}, |Q_{(x,t)}|^{1/p^-} \}$ by \cite[Lemma 3.2.12]{DieningHastoBuch2011} and $|Q_{(x,t)}|=(2t)^d$ we obtain
$\|\chi_{Q_{(x,t)}}|L_\p \| \gtrsim t^{d/ p^-}$
and finally arrive at
\begin{align*}
 \| \chi_{U_{(x,t)}} | M^w_{\p,\q} \| \gtrsim t^{-\alpha_1}(1+|x|)^{-\alpha_3}\|\chi_{Q_{(x,t)}}|L_\p(\R)\| \gtrsim   (\vB(x,t))^{-1}\,,
\end{align*}
where $\chi_{Q_{(x,t)}}(y)  (1+|x-y|/t)^{-\alpha_3}\asymp \chi_{Q_{(x,t)}}(y)$ was used. If $t=\infty$ we can argue analogously.
\eproof

\subsection{Coorbit identification}

As the following lemma shows,
every admissible wavelet frame $\cf=\cf(\Phi_0,\Phi)$ in the sense of
Definition~\ref{def:admfr} is suitable for the definition of coorbits of
Peetre-Wiener spaces.

\vspace*{-4ex}
\begin{flalign}\label{standassump}
\begin{minipage}[t]{0.9\textwidth}
\textbf{Standing assumptions:} For the rest of the
paper the indices fulfill $p,q\in\Plog$ with $0< p^-\leq p^+<\infty$, $0<
q^-\leq q^+<\infty$. Further $\qconst\in(0,\infty]$ and
$w\in\mathcal{W}^{\alpha_3}_{\alpha_1,\alpha_2}$ for arbitrary but fixed
$\alpha_2\ge \alpha_1$ and $\alpha_3\geq0$.
\end{minipage} &&
\end{flalign}

\begin{lemma}\label{lem:R(Y)}
An admissible continuous wavelet frame $\cf$ in the sense of \eqref{eqdef:wavefr} with generators $\Phi_0\in\mathcal{S}(\R)$ and $\Phi\in\mathcal{S}_0(\R)$
has property $F(\nu,Y)$ for $Y=\Pe(X)$ and $Y=\Le(X)$, and where $\nu=\vB$ is the corresponding weight from \eqref{eqdef:assweight}.
\end{lemma}
\bproof The proof goes along the lines of \cite[Lem.\ 4.18]{RaUl10}.
The kernel estimates in \cite[Lem.\ 4.8, 4.24]{RaUl10} have to be adapted to the
Peetre-Wiener space. This is a straight-forward procedure and allows for
treating as well the quasi-Banach situation\,.
\eproof

Now we are ready for the coorbit characterization of $B^{w}_{\p,\qconst}(\R)$ and $F^{w}_{\p,\q}(\R)$.
Note that the weight $\tilde{w}$ defined in \eqref{eqdef:wtilde} is an element of the class $\mathcal{W}^{\alpha_3}_{\alpha_1+d/2,\alpha_2+d/2}$.

\begin{Theorem}\label{thm:coident} Let $\p,\ \q,\ \qconst,\ w$ fulfill the
standing assumptions \eqref{standassump}. We choose an admissible continuous wavelet frame
$\cf=\cf(\Phi_0,\Phi)$ according to Definition \ref{def:admfr}. Putting
\begin{equation}\label{eqdef:wtilde}
      \tilde{w}(x,t) := \left\{\begin{array}{rcl}
                                  t^{-d/2}w(x,t) &,& 0<t< 1\,,\\
                                  w(x,\infty)&,& t=\infty\,,
                               \end{array}\right.
\end{equation}
we have $B^{w}_{\p,\qconst}(\R)= \Co (\cf,L^{\tilde{w}}_{\p,\qconst,a})$ if $a>\frac{d}{p^-}+\alpha_3$ and
$F^{w}_{\p,\q}(\R) = \Co (\cf,P^{\tilde{w}}_{\p,\q,a})$ if $a>\max\{\frac{d}{p^-},\frac{d}{q^-}\}+\alpha_3$
in the sense of equivalent quasi-norms.
\end{Theorem}
\bproof
By Lemma~\ref{lem:R(Y)} the coorbits exist in accordance with the theory.
Now, let $f\in\mathcal{S}(\R)$ and $F(x,t):=V_{\cf}f(x,t) = \langle f,
\varphi_{(x,t)}\rangle$ with $\varphi_{(x,t)}$ as in \eqref{eqdef:wavefr}. According to \cite[Lem.\ A.3]{T10} 
\begin{align*}
|V_{\cf}f(x,t)|\le C_N(f) G_N(x,t) \quad\text{with}\quad G_N(x,t)=\begin{cases} t^N(1+|x|)^{-N} \quad&,\, 0<t<1, \\
(1+|x|)^{-N} &,\, t=\infty, \end{cases}
\end{align*}
where $N\in\N$ is arbitrary but fixed and $C_N(f)>0$ is a constant depending on $N$ and $f$.
Choosing $N$ large, we have $G_N\in L_1^\nu(X)$ and thus $F\in L_1^\nu(X)$ with $\| F | L_1^\nu \|\le C_N(f) \| G_N | L_1^\nu \|$. This proves $f\in \mathcal{H}_1^\nu$.
Even more, given a sequence $(f_n)_{n\in\N} \subset \mathcal{S}(\R)$ we have $C_N(f_n)\rightarrow 0$ if $f_n\rightarrow 0$ in $\mathcal{S}(\R)$.
This is due to the fact that the constants $C_N(f_n)$ can be estimated by the Schwartz semi-norms of $f_n$ up to order $N$ (see proof of \cite[Lem.\ A.3]{T10}).
Hence, $\mathcal{F}\subset\mathcal{S}(\R)\hookrightarrow \mathcal{H}_1^\nu$ and
the voice transform $V_\cf$ extends to $\mathcal{S}^\prime(\R)$. Moreover, by a straight-forward modification of the argument in \cite[Cor.\ 20.0.2]{Hol95},
the reproducing formula is still valid on $\mathcal{S}^\prime(\R)$. Therefore we may apply Lemma~\ref{lem:co_indires} and use the larger reservoir $\mathcal{S}^\prime(\R)$.

To see that the coorbits coincide with $B^{w}_{\p,\qconst}(\R)$ and $F^{w}_{\p,\q}(\R)$, note that
the functions $\tilde{\Phi}=\overline{\Phi(-\cdot)}$ and $\tilde{\Phi}_0=\overline{\Phi_0(-\cdot)}$
satisfy the Tauberian conditions \eqref{condphi1}, \eqref{condphi2} and can thus be used in the continuous characterization of Theorem~\ref{thm:contchar}.
Recall the notation $\tilde{\Phi}_t= t^{-d} \tilde{\Phi}(\cdot/t)$.
The assertion is now a direct consequence of the possible reformulation $(V_{\cf} f)(\cdot,\infty) = \tilde{\Phi}_0\ast f$ and
$$
   (V_{\cf} f)(x,t) = \left(\mathcal{D}^{L_2}_t \overline{\Phi(-\cdot)} \ast  f
\right)(x) = t^{d/2} \left( \tilde{\Phi}_t \ast f \right)(x) \qquad, 0<t<1,
x\in\R.
$$
\eproof

\subsection{Atomic decompositions and quasi-Banach frames}

Based on the coorbit characterizations 
of Theorem~\ref{thm:coident} we can now apply the abstract theory from Section \ref{sec:abstrth} in our concrete setup, in particular the
discretization machinery.
We will subsequently use the following covering of the space $X$.
For $\alpha>0$ and $\beta>1$ we consider
the family $\mathcal{U}^{\alpha,\beta} = \{U_{j,k}\}_{j\in \n, k\in \zz^d}$ of subsets
\begin{align}\nonumber
    U_{0,k} &= Q_{0,k}\times\{\infty\}
    \quad,\quad k\in \Z\,,\notag\\
U_{j,k} &= Q_{j,k} \times [\beta^{-j},\beta^{-j+1})\quad , \quad j \in \N, k \in \zz^d\,,\nonumber
\end{align}
where $Q_{j,k} = \alpha \beta^{-j} k + \alpha \beta^{-j}[0,1]^d$.
Clearly, we have $X \subset \bigcup_{j\in \n, k\in \Z} U_{j,k}$ and $\mathcal{U}=\mathcal{U}^{\alpha,\beta}$ is an admissible covering of $X$.

The abstract Theorem~\ref{thm:atomicdec} provides atomic decompositions for
$B^w_{\p,\qconst}(\R)$ and $F^w_{\p,\q}(\R)$.
To apply it we need to analyze the oscillation kernels
$\osc_{\alpha,\beta}:=\osc_{\mathcal{U},\Gamma}$ and
$\osc^{\ast}_{\alpha,\beta}:=\osc^{\ast}_{\mathcal{U},\Gamma}$, where we choose the trivial phase function $\Gamma\equiv1$.
This goes along the lines of \cite[Sect. 4.4]{RaUl10}

\begin{proposition}\label{bddkernels} Let $\cf =\cf(\Phi_0,\Phi)$ be an admissible wavelet frame, $Y=\Le(X)$ or $Y=\Pe(X)$,
and $\nu=\vB$ the associated weight~\eqref{eqdef:assweight}.
\begin{description}
     \item(i) The kernels $\osc_{\alpha,\beta}$ and $\osc^{\ast}_{\alpha, \beta}$ are bounded operators on
     $Y$ and belong to $\mathcal{A}_{m_\nu}$.
     \item(ii) If $\alpha\downarrow 0$ and $\beta\downarrow 1$ then $\|\osc_{\alpha,\beta} | \mathcal{B}_{Y,{m_\nu}} \| \to 0$ and $\|\osc^{\ast}_{\alpha,\beta} | \mathcal{B}_{Y,{m_\nu}} \| \to 0$.
  \end{description}
\end{proposition}
\bproof The proof is a straight-forward modification of \cite[Lem.\
4.22]{RaUl10}. Similar as in Lemma \ref{lem:R(Y)} above we have to adapt the
kernel estimates to the Peetre-Wiener spaces.
\eproof

Finally, Theorem~\ref{thm:atomicdec}
yields the following discretization result in our concrete setting, which we
only state
for $\LizTri$ since for $\Besov$ it is essentially the same. 

\begin{Theorem}
Let $\p,\ \q,\ w$ fulfill the standing assumptions \eqref{standassump}, assume further $a>\max\{d/p^-, d/q^-\}+\alpha_3$ and let $\tilde{w}$ be given as in \eqref{eqdef:wtilde}.
For an admissible continuous wavelet frame $\cf = \{\varphi_\x\}_{\x\in X}$ there exist $\alpha_0>0$  and $\beta_0>1$,
such that for all $0<\alpha\leq \alpha_0$ and $1< \beta\leq \beta_0$ the family
$\cf_d = \{\varphi_{\x_{j,k}}\}_{j\in \N_0, k \in \Z}$ with $\x_{j,k} = (\alpha k \beta^{-j}, \beta^{-j})$ for $j\in\N$ and
$\x_{0,k} = (\alpha k , \infty) $ is a discrete wavelet frame with a corresponding dual frame $\mathcal{E}_d = \{e_{j,k}\}_{j\in \N_0, k\in \Z}$
such that
\begin{description}
    \item(a) If $f\in F^w_{\p,\q}(\R)$ we have the quasi-norm equivalence
    \begin{align*}
        \|f|F^w_{\p,\q}(\R)\| &\asymp \|\{\langle f,\varphi_{x_{j,k}}\rangle\}_{j\in \N_0, k\in \Z}|(P^{\tilde{w}}_{\p,\q,a})^{\flat}\|\\
                   &\asymp \|\{\langle f,e_{j,k} \rangle\}_{j\in \N_0, k\in \Z}|(P^{\tilde{w}}_{\p,\q,a})^{\natural}\|\,.
    \end{align*}
    \item(b) For every $f\in F^w_{\p,\q}(\R)$ the series
      \hfill$\displaystyle{  f = \sum\limits_{j\in \N_0}\sum\limits_{k\in \Z}\langle f,e_{j,k}\rangle \varphi_{x_{j,k}}
        =\sum\limits_{j\in \N_0}\sum\limits_{k\in \Z}\langle f,\varphi_{x_{j,k}}\rangle e_{j,k} } $\hfill~\\
    converge unconditionally in the quasi-norm of $F^w_{\p,\q}(\R)$.
\end{description}
\end{Theorem}
\bproof
The assertion is a consequence of the representation $F^w_{\p,\q}(\R)=\Co (\cf,P^{\tilde{w}}_{\p,\q,a})$ and Theorem~\ref{thm:atomicdec}.
In fact, Proposition~\ref{bddkernels} proves that $\mathcal{F}$ has property $D(\delta,\nu,Y)$ and $D(\delta,\nu,L_2)$ for every $\delta>0$.
Also note that $(\Pe)^{\flat}=(\Pe)^{\natural}$ with equivalent quasi-norms.
\eproof

\subsection{Wavelet bases}

According to Appendix \ref{sect:OWT} we obtain a family of
systems $\mathcal{G}_c,\,c\in E:=\{0,1\}^d$, whose union constitutes a tensor
wavelet system in $L_2(\R)$. Our aim is now to apply the abstract result in
Theorem~\ref{thm:frameexp} to achieve wavelet basis characterizations of
$B^w_{\p,\qconst}(\R)$ and $F^w_{\p,\q}(\R)$.
We have to consider the Gramian cross kernels
$K_c=K_\mathcal{U}[\mathcal{G}_c,\mathcal{F}]$ and
$K^{\ast}_c=K^\ast_\mathcal{U}[\mathcal{G}_c,\mathcal{F}]$ from
\eqref{eqdef:kerK} in our
concrete setup.

\begin{lemma}\label{lem:wavecond}
Let $Y=\Le(X)$ or $Y=\Pe(X)$ with associated weight $\nu=\vB$ given in
\eqref{eqdef:assweight}.  Assume that $a>0$ and $\p,\ \q,\ \qconst,\ w$ fulfill
the standing assumptions \eqref{standassump}. Let further $\cf=\cf(\Phi_0,\Phi)$ be an admissible
wavelet frame, $\cg_c$ be the systems from above, and
$K_c=K_\mathcal{U}[\mathcal{G}_c,\mathcal{F}]$,
$K^{\ast}_c=K^\ast_\mathcal{U}[\mathcal{G}_c,\mathcal{F}]$, $c \in E$, the
corresponding Gramian cross kernels. Then the kernels $K_c$ and $K^{\ast}_c$
define bounded operators from $Y$ to $Y$.
\end{lemma}
\bproof The proof is analogous to the treatment of the kernels $\osc$ in
Proposition \ref{bddkernels}, see also \cite[Lem.\ 4.24]{RaUl10}.\eproof

Now we are ready for the discretization of $\Besov$ and $\LizTri$ in terms of orthonormal wavelet bases.
We again only state the result for $\LizTri$ for the sake of brevity.

\begin{Theorem}
Let $\p,\ \q,\ w\in\mathcal{W}^{\alpha_3}_{\alpha_1,\alpha_2}$ fulfill the standing assumptions \eqref{standassump}, assume further
$a>\max\{d/p^-, d/q^-\} +\alpha_3$ and let $\tilde{w}$ be given as in
\eqref{eqdef:wtilde}. Let $\psi^0, \psi^1 \in L_2(\re)$ be the Meyer
scaling function and associated wavelet.
Then every $f\in\LizTri$ has the decomposition
\begin{equation*}
   \begin{split}
        f  =&  \sum\limits_{c\in E}\sum\limits_{k\in \Z}\lambda^c_{0,k}
\psi^c(\cdot-k)+\sum\limits_{c\in E\setminus \{0\}}\sum\limits_{j\in
\N}\sum\limits_{k\in \Z}  \lambda^c_{j,k}2^{\frac{jd}{2}}\psi^c(2^j\cdot-k)
   \end{split}
\end{equation*}
with quasi-norm convergence in $F^w_{\p,\q}(\R)$ and
sequences $\lambda^c = \{\lambda^c_{j,k}\}_{j\in \N_0,k\in\Z}$ defined by
$$
   \lambda^c_{j,k} = \langle
f,2^{\frac{jd}{2}}\psi^c(2^j\cdot-k)\rangle_{\mathcal{S}^\prime\times\mathcal{S}
}\quad, \quad j\in \N_0,k\in \Z\,,
$$
which belong to the sequence space $(P^{\tilde{w}}_{\p,\q,a})^\natural$
for every $c\in E$.
Conversely, an element $f\in (\mathcal{H}_{\nu_{\tilde{w},\p,\q}}^1)^{\urcorner}$
belongs to $\LizTri$ if all sequences
$\lambda^c(f)$ belong to $(P^{\tilde{w}}_{\p,\q,a})^\natural$.

\end{Theorem}
\bproof
The statement is a direct consequence of Theorem~\ref{thm:coident} and Theorem~\ref{thm:frameexp}.
The required conditions of the kernels $K_c,K^{\ast}_c$, $c\in E$, have been proved in Lemma~\ref{lem:wavecond}.
\eproof

\setcounter{section}{0}
\renewcommand{\thesection}{\Alph{section}}
\renewcommand{\theequation}{\Alph{section}.\arabic{equation}}
\renewcommand{\theTheorem}{\Alph{section}.\arabic{Theorem}}
\renewcommand{\thesubsection}{\Alph{section}.\arabic{subsection}}

\section{Appendix: Wavelet transforms}
\subsection{The continuous wavelet transform}
\label{Appendix:WaveletTransfoms}
As usual $\mathcal{S}(\re^d)$ denotes
the locally convex space of rapidly decreasing infinitely
differentiable functions on $\re^d$ and its topological dual is denoted by $\mathcal{S}'(\re^d)$.
The Fourier transform defined on both $\mathcal{S}(\R)$ and $\mathcal{S}'(\R$)
is given by 
$\widehat{f}(\varphi) := f (\widehat{\varphi})$,
where
$\,f\in \mathcal{S}'(\R), \varphi \in \mathcal{S}(\R)$, and
$$
\widehat{\varphi}(\xi):= (2\pi)^{-d/2}\int_{\R} e^{-ix\cdot
\xi}\varphi(x)\,dx.
$$
The Fourier transform is a bijection (in both cases) and its inverse is
given by $\varphi^{\vee} = \widehat{\varphi}(-\cdot)$.

Let us introduce the continuous wavelet transform. A general reference is provided by the monograph \cite[2.4]{Dau92}.
For $x\in \R$ and $t>0$ we define the unitary dilation and translation operators $\mathcal{D}^{L_2}_{t}$
and $T_x$ by
$$
  \mathcal{D}^{L_2}_{t}g = t^{-d/2}g\Big(\frac{\cdot}{t}\Big)\quad\mbox{and}\quad T_xg = g(\cdot-x)\quad,\quad g\in L_2(\R)\,.
$$
The vector $g$ is said to be the analyzing vector for a function $f\in L_2(\R)$. The
continuous wavelet transform $W_gf$ is then defined by
$$
    W_g f(x,t) = \langle T_x\mathcal{D}^{L_2}_t g, f\rangle\quad,\quad x\in \R, t>0\,,
$$
where the bracket $\langle \cdot, \cdot \rangle$ denotes the
inner product in $L_2(\R)$.
We call $g$ an admissible wavelet if
$$
    c_g:= \int_{\R} \frac{|\widehat{g}(\xi)|^2}{|\xi|^d}\,d\xi < \infty\,.
$$
If this is the case, then the family $\{T_x\mathcal{D}^{L_2}_t g\}_{t>0, x\in \R}$ represents a tight
continuous frame in $L_2(\re)$ where $C_1 = C_2 = c_g$.

Many consideration in this paper are based on decay results of the continuous wavelet transform $W_g f(x,t)$. This decay
mainly depends on moment conditions of the analyzing vector $g$ as well as on the smoothness of $g$
and the function $f$ to be analyzed, see \cite[Lem.\ A.3]{T10} which is based
on \cite[Lem.\ 1]{Ry99a}

\subsection{Orthonormal wavelets}

\label{sect:OWT}
\subsubsection*{The Meyer wavelets}
Meyer wavelets were introduced in \cite{Meyer86} and are an important example of wavelets which belong to the Schwartz class $\mathcal{S}(\re)$. The scaling function $\psi^0\in\mathcal{S}(\re)$ and the wavelet $\psi^1\in\mathcal{S}(\re)$ are real, their Fourier transforms are compactly supported and they fulfill
\begin{align*}
 \hat{\psi}^0(0)=(2\pi)^{-1/2}\quad\text{and}\quad \supp\hat{\psi}^1\subset\left[-\frac83\pi,-\frac23\pi\right]\cup\left[\frac23\pi,\frac83\pi\right].
\end{align*}
Due to the support condition we have infinitely many moment conditions \eqref{condphi2} on $\psi^1$ and both functions are fast decaying and infinitly often differentiable, see \cite[Section 3.2]{Wo97} for more properties.

\subsubsection*{Wavelets on $\R$}

In order to treat function spaces on $\R$ let us recall the construction of
a $d$-variate wavelet basis out of a resolution of unity in $\R$, see for
instance Wojtaszczyk \cite{Wo97}. It starts with a scaling function $\psi^0$ and
a wavelet $\psi^1$ belonging to $L_2(\re)$.
For $c\in E= \{0,1\}^d$ the function $\psi^c:\R \to \re$ is then defined by the
tensor product
$\psi^c = \bigotimes_{i=1}^d \psi^{c_i}$, i.e., $\psi^c(x) =  \prod_{i=1}^d
\psi^{c_i}(x_i)$, and we let $\mathcal{G}_c = \{\psi^c_{(x,t)}\}_{(x,t)\in X}$
be the system
with
\begin{align*}
\begin{aligned}
    \psi^c_{(x,t)} = \left\{\begin{array}{rcl}
                        T_x\mathcal{D}^{L_2}_t \psi^c&,&0<t<1\,,\\
                        T_x \psi^c&,& t=\infty\,,
                     \end{array}\right.
\end{aligned}
\quad \text{if $c\neq 0$} \quad\text{and}\quad
\begin{aligned}
    \psi^0_{(x,t)} = \left\{\begin{array}{rcl}
                        0&,&0<t<1\,,\\
                        T_x \psi^0&,& t=\infty\,.
                     \end{array}\right.
\end{aligned}
\end{align*}
\section*{Acknowledgement}
The authors would like to thank Felix Voigtlaender for a careful reading of the manuscript and many valuable comments and corrections.
They would further like to thank the anonymous referees for their careful proofreading and many valuable remarks. In particular, Lemma~\ref{lem:Bochner} should be pointed out, where now -- based on their input -- also questions of Bochner-measurability are discussed. Furthermore,
a serious technical issue in the proof of Theorem~\ref{thm:contchar} has been fixed.
M.S.\ would like to thank Holger Rauhut for support during his diploma studies where some ideas of this paper were developed.


\begin{thebibliography}{10}

\bibitem{alanga93}
S.~T. {A}li, J.-P. {A}ntoine, and J.-P. {G}azeau.
\newblock {C}ontinuous frames in {H}ilbert space.
\newblock {\em {A}nn. {P}hysics}, 222(1):1--37, 1993.

\bibitem{AlmeidaHasto2010}
A.~Almeida and P.~H{\"a}st{\"o}.
\newblock Besov spaces with variable smoothness and integrability.
\newblock {\em J. Funct. Anal.}, 258(5):1628--1655, 2010.

\bibitem{Ao42}
T.~Aoki.
\newblock Locally bounded linear topological spaces.
\newblock {\em Proc. Imp. Acad. Tokyo}, 18:588--594, 1942.

\bibitem{Bony}
J.-M. Bony.
\newblock  {Second Microlocalization and Propagation of Singularities for Semi-Linear Hyperbolic Equations.}
\newblock  {\em Taniguchi Symp. HERT. Katata}, 11--49, 1984.

\bibitem{BuPaTa96}
H.-Q. {B}ui, M.~{P}aluszy\'nski, and M.~H. {T}aibleson.
\newblock A maximal function characterization of weighted {B}esov-{L}ipschitz
  and {T}riebel-{L}izorkin spaces.
\newblock {\em Stud. Math.}, 119(3):219--246, 1996.

\bibitem{CoDaFe92}
A.~Cohen, I.~Daubechies, and J.-C. Feauveau.
\newblock Biorthogonal bases of compactly supported wavelets.
\newblock {\em Comm. Pure Appl. Math.}, 45(5):485--560, 1992.

\bibitem{CruzUribe2003}
D.~Cruz-Uribe, A.~Fiorenza, and C.~J. Neugebauer.
\newblock The maximal function on variable {$L^p$} spaces.
\newblock {\em Ann. Acad. Sci. Fenn. Math.}, 28(1):223--238, 2003.

\bibitem{daforastte08}
S.~{D}ahlke, M.~{F}ornasier, H.~{R}auhut, G.~{S}teidl, and G.~{T}eschke.
\newblock {G}eneralized coorbit theory, {B}anach frames, and the relation to
  alpha-modulation spaces.
\newblock {\em {P}roc. {L}ondon {M}ath. {S}oc. (3)}, 96:464--506, 2008.

\bibitem{dakustte09}
S.~{D}ahlke, G.~{K}utyniok, G.~{S}teidl, and G.~{T}eschke.
\newblock {S}hearlet coorbit spaces and associated {B}anach frames.
\newblock {\em {A}ppl. {C}omput. {H}armon. {A}nal.}, 27(2):195--214, 2009.

\bibitem{dastte04}
S.~{D}ahlke, G.~{S}teidl, and G.~{T}eschke.
\newblock {C}oorbit spaces and {B}anach frames on homogeneous spaces with
  applications to the sphere.
\newblock {\em {A}dv. {C}omput. {M}ath.}, 21(1-2):147--180, 2004.

\bibitem{dastte04-1}
S.~{D}ahlke, G.~{S}teidl, and G.~{T}eschke.
\newblock {W}eighted coorbit spaces and {B}anach frames on homogeneous spaces.
\newblock {\em {J}. {F}ourier {A}nal. {A}ppl.}, 10(5):507--539, 2004.


\bibitem{dadelasttevi14}
S. {D}ahlke, F. {D}e {M}ari, E. {D}e {V}ito, D. {L}abate, G. {S}teidl, G. {T}eschke, S. {V}igogna.
\newblock {C}oorbit spaces with voice in a Fr\'{e}chet space.
\newblock {\em {J}. {F}ourier {A}nal. {A}ppl.}, DOI 10.1007/s00041-016-9466-x.


\bibitem{Dau92}
I.~Daubechies.
\newblock {\em {T}en {L}ectures on {W}avelets}, volume~61 of {\em CBMS-NSF
  Regional Conference Series in Applied Mathematics}.
\newblock Society for Industrial and Applied Mathematics (SIAM), Philadelphia,
  PA, 1992.

\bibitem{Diening2004}
L.~Diening.
\newblock Maximal function on generalized {L}ebesgue spaces {$L^{p(\cdot)}$}.
\newblock {\em Math. Inequal. Appl.}, 7(2):245--253, 2004.

\bibitem{DieningHarjulehto2009}
L.~Diening, P.~Harjulehto, P.~H{\"a}st{\"o}, Y.~Mizuta, and T.~Shimomura.
\newblock Maximal functions in variable exponent spaces: limiting cases of the
  exponent.
\newblock {\em Ann. Acad. Sci. Fenn. Math.}, 34(2):503--522, 2009.

\bibitem{DieningHastoBuch2011}
L.~Diening, P.~Harjulehto, P.~H\"ast\"o, and M.~R{\ocirc{u}}{\v{z}}i{\v{c}}ka.
\newblock {\em Lebesgue and Sobolev Spaces with Variable Exponents}.
\newblock Number 2017 in Lecture Notes in Mathematics. Springer, 2011.

\bibitem{DieningHastoRoudenko2009}
L.~Diening, P.~H{\"a}st{\"o}, and S.~Roudenko.
\newblock Function spaces of variable smoothness and integrability.
\newblock {\em J. Funct. Anal.}, 256(6):1731--1768, 2009.

\bibitem{Fe83}
H.~G. {F}eichtinger.
\newblock {B}anach convolution algebras of {W}iener's type.
\newblock In {\em {F}unctions, {S}eries, {O}perators; {P}roc. {C}onf.
  {B}udapest 1980, {V}ol. {I}}, pages 509--524, 1983.

\bibitem{fe83-4}
H.~G. {F}eichtinger.
\newblock {M}odulation spaces on locally compact {A}belian groups.
\newblock Technical report, {J}anuary 1983.

\bibitem{fegr85}
H.~G. {F}eichtinger and P.~{G}r{\"o}bner.
\newblock {B}anach spaces of distributions defined by decomposition methods.
  {I}.
\newblock {\em {M}ath. {N}achr.}, 123:97--120, 1985.

\bibitem{FeGr86}
H.~G. Feichtinger and K.~Gr{\"o}chenig.
\newblock A unified approach to atomic decompositions via integrable group
  representations.
\newblock In {\em Function spaces and applications (Lund, 1986)}, volume 1302
  of {\em Lecture Notes in Math.}, pages 52--73. Springer, Berlin, 1988.

\bibitem{feGr89a}
H.~G. {F}eichtinger and K.~{G}r{\"o}chenig.
\newblock {B}anach spaces related to integrable group representations and their
  atomic decompositions, {I}.
\newblock {\em Journ. Funct. Anal.}, 21:307--340, 1989.

\bibitem{Fol84}
G.~B. Folland.
\newblock {\em Real analysis}.
\newblock Pure and Applied Mathematics (New York). John Wiley \& Sons, Inc.,
  New York, 1984.
\newblock Modern techniques and their applications, A Wiley-Interscience
  Publication.

\bibitem{fora05}
M.~{F}ornasier and H.~Rauhut.
\newblock {C}ontinuous frames, function spaces, and the discretization problem.
\newblock {\em {J}. {F}ourier {A}nal. {A}ppl.}, 11(3):245--287, 2005.


\bibitem{balhol10}
P.~{B}alazs and N.~{Holighaus}.
\newblock {D}iscretization in generalized coorbit spaces: extensions, annotations and errata for ``
{C}ontinuous frames, function spaces, and the discretization problem'' by M.~{F}ornasier and H.~{R}auhut. Online available: {\url{https://www.univie.ac.at/nonstatgab/warping/baho15.pdf}}

\bibitem{balholwies15}
P.~{B}alazs, N.~{Holighaus}, and C.~{W}iesmeyr.
\newblock {C}onstruction of warped time-frequency representations on nonuniform frequency scales, part ii: Integral transforms, functions spaces, atomic decompositions and Banach frames.
\newblock Online available: {arXiv: 1503.05439, 2005.}


\bibitem{fu13a}
H.~{F}{\"{u}}hr.
\newblock {C}oorbit spaces and wavelet coefficient decay over general dilation
  groups.
\newblock {\em {T}rans. {A}mer. {M}ath. {S}oc.}, 367:7373--7401, 2015.

\bibitem{fu13b}
H.~{F}{\"{u}}hr.
\newblock {V}anishing moment conditions for wavelet atoms in higher dimensions.
\newblock {\em {A}dv. {C}omput. {M}ath.}, 42(1):127--153, 2015.


\bibitem{furaitou15}
H.~{F}{\"{u}}hr and R.~{R}aisi{-T}ousi.
\newblock {S}implified vanishing moment criteria for wavelets over general
  dilation groups, with applications to abelian and shearlet dilation groups.
\newblock {\em ArXiv e-prints}, 2015.
\newblock arXiv:1407.0824 [math.FA].

\bibitem{fuvoigt14}
H.~{F}{\"{u}}hr and F.~{V}oigtlaender.
\newblock {W}avelet coorbit spaces viewed as decomposition spaces.
\newblock {\em J. Funct. Anal.}, 269:80-154, 2015.

\bibitem{Gr88}
K.~Gr{\"o}chenig.
\newblock Unconditional bases in translation and dilation invariant function
  spaces on {$\mathbb{R}\sp n$}.
\newblock In {\em Constructive theory of functions (Varna, 1987)}, pages
  174--183. Publ. House Bulgar. Acad. Sci., Sofia, 1988.

\bibitem{Gr91}
K.~{G}r{\"o}chenig.
\newblock {D}escribing functions: atomic decompositions versus frames.
\newblock {\em {M}onatsh. {M}athem.}, 112:1--41, 1991.

\bibitem{gr01}
K.~{G}r{\"o}chenig.
\newblock {\em {F}oundations of {T}ime-{F}requency {A}nalysis}.
\newblock {A}ppl. {N}umer. {H}armon. {A}nal. {B}irkh{\"a}user {B}oston, 2001.

\bibitem{Hol95}
M.~{H}olschneider.
\newblock {\em {W}avelets, {A}n {A}nalysis {T}ool}.
\newblock {O}xford {U}niversity {P}ress, 1995.

\bibitem{Jaffard}S. Jaffard, Y. Meyer.
\newblock {Wavelet methods for pointwise regularity and local oscillations of functions.}
\newblock {\em Memoirs of the AMS}, vol. {123}, 1996.

\bibitem{Ke09}
H.~{K}empka.
\newblock $2$-microlocal Besov and Triebel-Lizorkin spaces of variable integrability.
\newblock {\em {R}ev. {M}at. {C}omplut.}, 22(1):227--251, 2009.

\bibitem{Ke11}
H.~{K}empka.
\newblock Atomic, molecular and wavelet decomposition of 2-microlocal {B}esov
  and {T}riebel-{L}izorkin spaces with variable integrability.
\newblock {\em Funct. Appr.}, 43(2):171--208, 2010.

\bibitem{Ke10}
H.~{K}empka.
\newblock Atomic, molecular and wavelet decomposition of generalized
  2-microlocal {B}esov spaces.
\newblock {\em Journ. Function Spaces Appl.}, 8:129--165, 2010.

\bibitem{KeVybDiff}
H.~Kempka and J.~Vyb{\'{\i}}ral.
\newblock Spaces of variable smoothness and integrability: characterizations by
  local means and ball means of differences.
\newblock {\em J. Fourier Anal. Appl.}, 18(4):852--891, 2012.

\bibitem{KeVybNorm}
H.~Kempka and J.~Vyb{\'{\i}}ral.
\newblock A note on the spaces of variable integrability and summability of
  {A}lmeida and {H}\"ast\"o.
\newblock {\em Proc. Amer. Math. Soc.}, 141(9):3207--3212, 2013.

\bibitem{KovacikRakosnik91}
O.~Kov{\'a}{\v{c}}ik and J.~R{\'a}kosn{\'{\i}}k.
\newblock On spaces {$L^{p(x)}$} and {$W^{k,p(x)}$}.
\newblock {\em Czechoslovak Math. J.}, 41(116)(4):592--618, 1991.

\bibitem{LiSaUlYaYu11}
Y.~Liang, Y.~Sawano, T.~Ullrich, D.~Yang, and W.~Yuan.
\newblock {N}ew characterizations of {B}esov-{T}riebel-{L}izorkin-{H}ausdorff
  spaces including coorbits and wavelets.
\newblock {\em J. Fourier Anal. Appl.}, 18(5):1067--1111, 2012.

\bibitem{LiSaUlYaYu12}
Y.~Liang, Y.~Sawano, T.~Ullrich, D.~Yang, and W.~Yuan.
\newblock A new framework for generalized {B}esov-type and
  {T}riebel-{L}izorkin-type spaces.
\newblock {\em Diss. Math.}, 489, 2013.

\bibitem{Meyer86}
Y.~Meyer.
\newblock Principe d'incertitude, bases hilbertiennes et alg\`{e}bres
  d'op\'{e}rateurs.
\newblock {\em S\'{e}m. Bourbaki}, 28:209--223, 1985-1986.

\bibitem{Nekvinda2004}
A.~Nekvinda.
\newblock Hardy-{L}ittlewood maximal operator on {$L^{p(x)}(\mathbb R)$}.
\newblock {\em Math. Inequal. Appl.}, 7(2):255--265, 2004.

\bibitem{Orlicz31}
W.~Orlicz.
\newblock {\"U}ber konjugierte {E}xponentenfolgen.
\newblock {\em Studia Math.}, 3:200--212, 1931.

\bibitem{Pe75}
J.~{P}eetre.
\newblock {O}n spaces of {T}riebel-{L}izorkin type.
\newblock {\em {A}rk. {M}at.}, 13:123--130, 1975.

\bibitem{ra05-3}
H.~{R}auhut.
\newblock {C}oorbit space theory for quasi-{B}anach spaces.
\newblock {\em {S}tudia {M}ath.}, 180(3):237--253, 2007.

\bibitem{ra05-4}
H.~{R}auhut.
\newblock {W}iener amalgam spaces with respect to quasi-{B}anach spaces.
\newblock {\em {C}olloq. {M}ath.}, 109(2):345--362, 2007.

\bibitem{RaUl10}
{H.~{R}auhut and T.~{U}llrich.}
\newblock Generalized coorbit space theory and inhomogeneous function spaces of
  {B}esov-{L}izorkin-{T}riebel type.
\newblock {\em {J}. {F}unct. {A}nal.}, 260(11):3299-3362, 2011.

\bibitem{Ro57}
S.~Rolewicz.
\newblock On a certain class of linear metric spaces.
\newblock {\em Bull. Acad. Polon. Sci. Cl. III.}, 5:471--473, XL, 1957.

\bibitem{Ry99a}
V.~S. {R}ychkov.
\newblock {O}n a theorem of {B}ui, {P}aluszy\'nski and {T}aibleson.
\newblock {\em {P}roc. {S}teklov {I}nst.}, 227:280--292, 1999.

\bibitem{Sch12}
M.~{S}ch{\"{a}}fer.
\newblock Generalized coorbit space theory for quasi-banach spaces.
\newblock Diplomarbeit, Rheinische Friedrich-Wilhelms-Universität Bonn, Germany, 2012.

\bibitem{Voigt15}
F.~{V}oigtlaender.
\newblock Embedding theorems for decomposition spaces with applications to wavelet coorbit spaces.
\newblock PhD thesis, RWTH Aachen University, Germany, 2015. Online available: {\url{https://publications.rwth-aachen.de/record/564979/files/564979.pdf}}


\bibitem{Tr83}
H.~{T}riebel.
\newblock {\em {T}heory of {F}unction {S}paces}.
\newblock {B}irkh{\"a}user, Basel, 1983.

\bibitem{Tr88}
H.~{T}riebel.
\newblock {C}haracterizations of {B}esov-{H}ardy-{S}obolev spaces: a unified
  approach.
\newblock {\em {J}ourn. of {A}pprox. {T}heory}, 52:162--203, 1988.

\bibitem{Tr92}
H.~{T}riebel.
\newblock {\em {T}heory of {F}unction {S}paces II}.
\newblock {B}irkh{\"a}user, Basel, 1992.

\bibitem{Tu14}
A.~I. Tyulenev.
\newblock Some new function spaces of variable smoothness.
\newblock {\em {M}at. {S}b.}, 206(6): 85–128, 2015.

\bibitem{T10}
T.~{U}llrich.
\newblock Continuous characterizations of {B}esov-{L}izorkin-{T}riebel spaces
  and new interpretations as coorbits.
\newblock {\em Journ. Funct. Spaces Appl.}, Article ID 163213, 2012.




\bibitem{Wo97}
P.~{W}ojtaszczyk.
\newblock {\em {A} {M}athematical {I}ntroduction to {W}avelets}.
\newblock {C}ambridge {U}niversity {P}ress, 1997.

\end{thebibliography}
\def\ocirc#1{\ifmmode\setbox0=\hbox{$#1$}\dimen0=\ht0 \advance\dimen0
  by1pt\rlap{\hbox to\wd0{\hss\raise\dimen0
  \hbox{\hskip.2em$\scriptscriptstyle\circ$}\hss}}#1\else {\accent"17 #1}\fi}

\end{document}